%% file: ROM_GQR.tex
\documentclass[11pt,a4paper]{article}
\usepackage{pdfpages}
\usepackage{ifpdf}
\pdfoutput=1
\input{packages.tex}

\input{commands.tex}
\usepackage{verbatim}
\usepackage{graphicx}
\usepackage{epstopdf}

\begin{document}
\date{}
\title{A POD-Galerkin reduced order model for a LES filtering approach}


\author[1]{Michele Girfoglio\thanks{mgirfogl@sissa.it}}
\author[2]{Annalisa Quaini\thanks{quaini@math.uh.edu}}
\author[1]{Gianluigi Rozza\thanks{grozza@sissa.it}}
\affil[1]{SISSA, International School for Advanced Studies, Mathematics Area, mathLab, via Bonomea, Trieste 265 34136, Italy}
\affil[2]{Department of Mathematics, University of Houston, Houston TX 77204, USA}
\maketitle

\begin{abstract}
We propose a Proper Orthogonal Decomposition (POD)-Galerkin based Reduced Order Model (ROM) 
for a Leray model.  
For the implementation of the model, we combine a two-step algorithm called Evolve-Filter (EF) with
a computationally efficient finite volume method. The main novelty of the proposed approach relies 
in applying spatial filtering both for the collection of the snapshots and in the reduced order model, as well as in considering the pressure field at reduced level.
In both steps of the EF algorithm, velocity and pressure fields are approximated by using different POD 
basis and coefficients. For the reconstruction of the pressures fields, we use 
a pressure Poisson equation approach. We test our ROM on two benchmark problems: 2D and 3D
unsteady flow past a cylinder at Reynolds number $0 \leq Re \leq 100$. 
The accuracy of the reduced order model is assessed against
results obtained with the full order model. 
For the 2D case, a parametric study with respect to the filtering radius is also presented.  


\end{abstract}
\vspace*{0.5cm}

\section{Introduction}\label{sec:intro}

Reduced order models (ROMs) have been proposed as an efficient tool for the approximation of systems
of parametrized partial differential equations, as they significantly reducing the computational 
cost required by classical full order models (FOMs),
e.g.~finite element methods or finite volume methods.
The basic ROM framework consists of two steps. First, a database of several solutions is 
collected by solving the original full order model for different physical and/or geometrical 
configurations (\emph{offline phase}). 
Then, the information obtained during the offline phase is used to compute the solution for 
newly specified values of the parameters in a short amount of time (\emph{online phase}). 
For a comprehensive review on ROMs, we refer
to \cite{hesthaven2015certified, quarteroniRB2016, bennerParSys, Benner2015, Bader2016, ModelOrderReduction}. 
 
It is well known that the extension to turbulent flows present several challenges. 
One of the reasons is that projection based ROMs of turbulent flows are affected by energy stability problems \cite{Carlberg2013623}. 
This is related to the fact that proper orthogonal decomposition (POD) retains the modes
biased toward large, high-energy scales, while the turbulent kinetic energy is dissipated
at level of the small turbulent scales \cite{Moin1998}. A possible strategy to stabilize ROMs for turbulent flows is 
the introduction of dissipation via a closure model \cite{wang_turb, Aubry1988}. In \cite{couplet_sagaut_basdevant_2003},
it was shown theoretically and numerically that POD modes have a similar energy transfer mechanism to Fourier modes. 
Therefore, the use of Large Eddy Simulation (LES) could be beneficial. Data-driven ROMs for LES full order models have been successfully used for hydroacoustic analysis \cite{GadallaCianferraTezzeleStabileMolaRozza2020}. In addition, the efficiency of ROMs for LES-VMS stabilized finite elements has been proved through relevant numerical benchmarks \cite{StabileBallarinZuccarinoRozza2019}.

We focus on projection based ROMs (see, e.g., \cite{Amsallem_stab}) 
that have been successfully applied to several fluid dynamics problems. 
We propose a POD-Galerkin-based ROM for a LES filtering approach. 
We consider a variant of the so-called Leray model \cite{Leray34}, 
where the small-scale effects are described by a set of equations to be added to the discrete 
Navier-Stokes equations. This extra problem acts as a differential low-pass filter \cite{abigail_CMAME}.
For its  actual implementation, we use the Evolve-Filter (EF) algorithm 
\cite{Boyd1998283, Fischer2001265, Dunca2005,layton_CMAME}.
The Leray model has been extensively applied within a Finite Element framework,
while we focus on Finite Volume (FV) methods \cite{Girfoglio2019}. In the ROM context, the
Leray model has been applied to the 1D Kuramoto-Sivashinsky equations \cite{Sabetghadam2012}, 
stochastic Burgers equation \cite{Iliescu2016}, and Navier-Stokes equations \cite{Gunzburger2019b, Xie2018,Wells2017}. 
The EF algorithm has also been investigated in combination with regularized ROMs: 
applications include stochastic Burgers equation \cite{Xie2018_2} and 3D Navier-Stokes equations \cite{Wells2017}. 
In \cite{Gunzburger2019}, a relaxation step is added to 
the EF algorithm and applied to the 2D Navier-Stokes equations,
while in \cite{Xie2016} a ROM is applied to the approximate deconvolution model. 
We note that in \cite{Gunzburger2019b, Xie2018,Wells2017,Gunzburger2019,Xie2016}
LES filtering is used in the development of the ROM to address the inaccuracy and numerical
instability of the standard Galerkin ROM for convection-dominated problems.
Unlike those works, we use LES filtering also as FOM, i.e.~to generate the snapshot data. Such an approach ensures a greater consistency between FOM and ROM.
The strategy we proposed is mentioned in \cite{Xie2016} as a future perspective. However, 
to the best of our knowledge it has not been attempted so far.

The novelties of our approach include:
\begin{itemize}
\item[-] the computation of the pressure field at ROM level (through the pressure Poisson equation \cite{Stabile2017, Stabile2018, Star2019});
\item[-] the use of different POD coefficients and bases to approximate the two velocity fields in the Leray model.
\end{itemize}
We test our approach on two benchmarks: 
2D \cite{turek1996, John2004} and 3D \cite{turek1996} flow past a cylinder with time-dependent Reynolds number $0 \leq Re(t) \leq 100$. 
For the 2D test, we compare the evolution of velocity and pressure fields with the corresponding FOM quantities.
Moreover, we present a parametrization with respect to the filtering radius, which is a crucial model parameter. 
We are not aware of any other work that considers a parametric case. 
The 3D test aims to show that our methodology is not limited to 2D problems. 
While for the 2D test we used the EF method to implement the Leray model, 
for the 3D case we use a ``monolithic'' approach for the Leray model that seems
to introduce slightly less artificial diffusion. 

In \cite{Girfoglio2019}, we confirmed that the EF algorithm is over-diffusive also when combined with
a FV method. Therein, we showed that the Evolve-Filter-Relax algorithm \cite{layton_CMAME} 
with nonlinear differential filters is a better choice, especially for realist problems. 
This work on a ROM approach for the EF algorithm is an intermediate step towards the development of ROMs
for the Evolve-Filter-Relax algorithm with nonlinear filters, which are more challenging. 




The work is organized as follows. Sec.~\ref{sec:FOM} describes 
the full order model and the numerical method we use for it. 
Sec.~\ref{sec:ROM} presents the ingredients of the reduced order model. 
The numerical examples are reported in Sec. \ref{sec:results}. 
Sec. \ref{sec:conclusions} provides conclusions and perspectives.

\section{The full order model}\label{sec:FOM}

We consider a fixed domain $\Omega \subset \mathbb{R}^D$ with $D = 2, 3$ over a time 
interval of interest ($t_0$, $T$) $\subset \mathbb{R}^+$. Let $\bm{\pi} \in \mathcal{P} \subset \mathbb{R}^d$ 
be a parameter vector in a $d$-dimensional parameter space $\mathcal{P}$. 
The so-called \emph{Leray model} couples the Navier-Stokes equations (NSE) with a differential filter:
\begin{align}
& \rho\, \dt \u (\bm{x},t; \bm{\pi})+ \rho\,\div \left(\ubar(\bm{x},t; \bm{\pi}) \otimes \u(\bm{x},t; \bm{\pi})\right) - 2\mu \Delta\u(\bm{x},t; \bm{\pi}) + \nabla p(\bm{x},t; \bm{\pi}) = 0, \label{eq:leray1}  \\
& \div \u(\bm{x},t; \bm{\pi}) = 0,  \label{eq:leray2}  \\
& -2 \alpha^2\Delta\ubar(\bm{x},t; \bm{\pi}) +\ubar(\bm{x},t; \bm{\pi}) +\nabla \lambda(\bm{x},t; \bm{\pi}) = \u(\bm{x},t; \bm{\pi}), \label{eq:leray3}  \\
& \div \ubar(\bm{x},t; \bm{\pi}) = 0, \label{eq:leray4}
\end{align}
in $\Omega \times (t_0,T)$, endowed with the boundary conditions. In \eqref{eq:leray1}-\eqref{eq:leray4}, $\rho$ is the fluid density, 
$\mu$ is the dynamic viscosity, $\u$ is velocity, $p$ is the pressure, $\ubar$ is the \emph{filtered velocity}, and
variable $\lambda$ is a Lagrange multiplier to enforce the incompressibility constraint for $\ubar$. 
Parameter $\alpha$ can be interpreted as the \emph{filtering radius}. 
Problem \eqref{eq:leray1}-\eqref{eq:leray4} is endowed with 
suitable boundary
\begin{align}
&\u(\bm{x},t; \bm{\pi}) = \u_D(\bm{x},t; \bm{\pi})  \quad \mbox{on } \partial\Omega_D \times(t_0,T), \label{eq:bc-d} \\ 
&(2\mu \nabla\u(\bm{x},t; \bm{\pi}) - p(\bm{x},t; \bm{\pi})\mathbf{I})\n = \0  \quad\mbox{on } \partial\Omega_N \times(t_0,T), \label{eq:bc-n} \\ 
&\ubar(\bm{x},t; \bm{\pi}) = \u_D(\bm{x},t; \bm{\pi})  \quad \mbox{on } \partial\Omega_D \times(t_0,T), \label{eq:bc-filter-d} \\ 
&(2\alpha^2 \nabla\ubar(\bm{x},t; \bm{\pi}) - \lambda(\bm{x},t; \bm{\pi})\mathbf{I})\n = \0  \quad\mbox{on } \partial\Omega_N \times(t_0,T).\label{eq:bc-filter-n}
\end{align}
and the initial data $\u(\bm{x},t) = \u_0(\bm{x})$ in $\Omega \times\{t_0\}$. Here $\overline{\partial\Omega_D}\cup\overline{\partial\Omega_N}=\overline{\partial\Omega}$ and $\partial\Omega_D \cap\partial\Omega_N=\emptyset$. In addition, $\u_D$ and $\u_0$ are given.

For the sake of simplicity, from now on the dependance of the variables on $\bm{x}$, $t$, and 
parameter $\bm{\pi}$ will be omitted.

\subsection{The Evolve-Filter algorithm}\label{sec:FOM_2}

We start with the time discretization of the Leray model \eqref{eq:leray1}-\eqref{eq:leray4}.
Let $\Delta t \in \mathbb{R}$, $t^n = t_0 + n \Delta t$, with $n = 0, ..., N_T$ and $T = t_0 + N_T \Delta t$. Moreover, we denote by $y^n$ the approximation of a generic quantity $y$ at the time $t^n$. We adopt a Backward Differentiation Formula of order 2 (BDF2) \cite{quarteroni2007numerical}.

To decouple the Navier-Stokes system \eqref{eq:leray1}-\eqref{eq:leray2} from the filter system \eqref{eq:leray3}-\eqref{eq:leray4} 
at each time step, we consider the Evolve-Filter (EF) algorithm \cite{Boyd1998283, Fischer2001265, Dunca2005}. This algorithm 
reads as follows: given the velocities $\u^{n-1}$ and $\u^{n}$, at $t^{n+1}$:
\begin{enumerate}[i)]
\item \textit{evolve}: find intermediate velocity and pressure $(\v^{n+1},q^{n+1})$ such that
\begin{align}
&\rho\, \frac{3}{2\Delta t}\, \v^{n+1} + \rho\, \div \left(\u^* \otimes \v^{n+1}\right) - 2\mu\Delta\v^{n+1} +\nabla q^{n+1} = \b^{n+1},\label{eq:evolve-1.1}\\
& \div \v^{n+1} = 0\label{eq:evolve-1.2},
\end{align}
with boundary conditions
\begin{align}
&\v^{n+1} = \u_D^{n+1}  \quad \mbox{on } \partial\Omega_D \times(t_0,T), \label{eq:bc-d_1} \\ 
&(2\mu \nabla\v^{n+1} - q^{n+1}\mathbf{I})\n = \0  \quad\mbox{on } \partial\Omega_N \times(t_0,T), \label{eq:bc-n_1} 
\end{align}
and initial condition $\v^0 = \u_0$ in $\Omega \times\{t_0\}$. In eq.~\eqref{eq:evolve-1.1}, we set $\u^* = 2 \u^n-\u^{n-1}$ and $\b^{n+1} = (4\u^n - \u^{n-1})/(2\Delta t)$.
\item \textit{filter}: find $(\u^{n+1},\lambda^{n+1})$ such that
\begin{align}
&-\alpha^2 \Delta\u^{n+1} +\u^{n+1} +\nabla \lambda^{n+1} = \v^{n+1}, \label{eq:evolve-2.1}\\
& \div \u^{n+1}  = 0 \label{eq:filter-1.2},
\end{align}
with boundary conditions
\begin{align}
& \u^{n+1} = \u_D^{n+1}  \quad \mbox{on } \partial\Omega_D \times(t_0,T), \label{eq:bc-filter-d_1} \\ 
& (2\alpha^2 \nabla\u^{n+1} - \lambda^{n+1}\mathbf{I})\n = \0  \quad\mbox{on } \partial\Omega_N \times(t_0,T).\label{eq:bc-filter-n_1}
\end{align}
\end{enumerate}
We consider $\u^{n+1}$ and $q^{n+1}$ the approximation of the velocity and the pressure at the time $t^{n+1}$, respectively.

In this work, we consider only homogeneous Neumann boundary conditions. For the treatment of
non-homogeneous boundary conditions, we refer to \cite{BQV}.

\begin{rem}\label{rem:gen_Stokes}
Filter problem \eqref{eq:evolve-2.1}-\eqref{eq:filter-1.2} can be considered 
a generalized Stokes problem. In fact, if we divide eq.~\eqref{eq:evolve-2.1}
by $\Delta t$ and rearrange the terms we obtain:
\begin{align}
\frac{\rho}{\Delta t} \u^{n+1}  - \mubar \Delta\u^{n+1} + \nabla \qbar^{n+1} & = \frac{\rho}{\Delta t} \v^{n+1}, \quad \mubar = \rho \frac{\alpha^2}{\Delta t}, \label{eq:filter-1.1}
\end{align}
where the filtered pressure $\qbar^{n+1} = \rho \lambda^{n+1}/\Delta t$. Problem \eqref{eq:filter-1.1},\eqref{eq:filter-1.2} can be
seen as a time dependent Stokes problem with viscosity $\mubar$, discretized by the Backward Euler (or BDF1) scheme. 
A solver for problem \eqref{eq:filter-1.1},\eqref{eq:filter-1.2} can then be obtained 
by adapting a standard linearized Navier-Stokes solver. 
\end{rem}



\subsection{Space discrete problem: a Finite Volume approximation}\label{sec:FOM_3}
For the space discretization of problems \eqref{eq:evolve-1.1}-\eqref{eq:evolve-1.2} and \eqref{eq:filter-1.1},\eqref{eq:filter-1.2}, 
we adopt a Finite Volume (FV) method. We partition the computational domain $\Omega$ into cells or control volumes $\Omega_i$, with $i = 1, \dots, N_{c}$, where $N_{c}$ is the total number of cells in the mesh. 
Let  \textbf{A}$_j$ be the surface vector of each face of the control volume, 
with $j = 1, \dots, M$. 

The fully discretized form of problem \eqref{eq:evolve-1.1}-\eqref{eq:evolve-1.2} is given by
\begin{align}
&\rho\, \frac{3}{2\Delta t}\, \v^{n+1}_i + \rho\, \sum_j^{} \varphi^*_j \v^{n+1}_{i,j} - 2\mu \sum_j^{} (\nabla\v^{n+1}_i)_j \cdot \textbf{A}_j + \sum_j^{} q^{n+1}_{i,j} \textbf{A}_j  = {\bm b}^{n+1}_i, \label{eq:disc_evolve1} \\
&\sum_j^{} (\nabla q^{n+1})_j \cdot \textbf{A}_j = \sum_j^{} (\textbf{H}(\v_i^{n+1}))_j \cdot \textbf{A}_j, \label{eq:disc_evolve2}
\end{align}
where:
\begin{align}
\textbf{H}(\v^{n+1}_i) = -\rho \sum_j^{} \varphi^*_j \v^{n+1}_{i,j} + 2\mu \sum_j^{} (\nabla\v^{n+1}_i)_j \cdot \textbf{A}_j + {\bm b}^{n+1}_i \quad \text{with} \quad \varphi^*_j = \u^{*}_j \cdot \textbf{A}_j. \label{eq:H}
\end{align}
In \eqref{eq:disc_evolve1}-\eqref{eq:H}, $\v^{n+1}_i$ and ${\bm b}^{n+1}_i$ denotes the average velocity and source term in control volume $\Omega_i$, respectively. Moreover, we denote with $\v^{n+1}_{i,j}$ and $q^{n+1}_{i,j}$ the velocity and pressure
associated to the centroid of face $j$ normalized by the volume of $\Omega_i$.

The fully discrete problem associated to the filter \eqref{eq:filter-1.1},\eqref{eq:filter-1.2} is given by 
\begin{align}
&\frac{\rho}{\Delta t} \u^{n+1}_i - \sum_j^{} \mubar_j(\nabla\u^{n+1}_i)_j \cdot \textbf{A}_j +  \sum_j^{} \qbar^{n+1}_{i,j} \textbf{A}_j = \frac{\rho}{\Delta t} \v^{n+1}_i, \label{eq:disc_filter1} \\
&\sum_j^{} (\nabla \qbar^{n+1}_i)_j \cdot \textbf{A}_j = \sum_j^{} (\overline{\textbf{H}}(\u^{n+1}_i))_j \cdot \textbf{A}_j, \label{eq:disc_filter2}
\end{align}
with
\begin{align}
\overline{\textbf{H}}(\u^{n+1}_i) =  \sum_j^{} \mubar_j(\nabla\u^{n+1}_i)_j \cdot \textbf{A}_j + \dfrac{\rho}{\Delta t}\v^{n+1}. \label{eq:Hbar}
\end{align}
In \eqref{eq:disc_filter1}-\eqref{eq:Hbar}, we denoted with $\u^{n+1}_i$ the average filtered velocity 
in control volume $\Omega_i$, while $\qbar^{n+1}_{i,j}$ is the auxiliary pressure
at the centroid of face $j$ normalized by the volume of $\Omega_i$.
For more details, we refer the reader to \cite{Girfoglio2019}.

We have implemented the EF algorithm within the C++ finite volume library OpenFOAM\textsuperscript{\textregistered} \cite{Weller1998}.
For the solution of the linear system associated with \eqref{eq:disc_evolve1}-\eqref{eq:disc_evolve2} we used the PISO algorithm \cite{PISO}, 
while for problem \eqref{eq:disc_filter1}-\eqref{eq:disc_filter2} we chose a slightly modified version of the SIMPLE algorithm
\cite{SIMPLE}, called SIMPLEC algorithm \cite{Doormaal1984}. Both PISO and SIMPLEC are partitioned algorithms
that decouple the computation of the pressure from the computation of the velocity. 

The approach described in this section represents our Full Order Model (FOM).  

\section{The reduced order model}\label{sec:ROM}
The Reduced Order Model (ROM) we propose is an extension of the model introduced in \cite{Stabile2017, Stabile2018}. In Sec~\ref{sec:ROM_1} we introduce the procedure we use to construct a POD-Galerkin ROM
and in Sec.~\ref{sec:ROM_2} we present the strategy we choose for pressure stabilization at reduced order level. 
Finally, Sec.~\ref{sec:ROM_3} describes the lifting function method  we apply to
enforce non-homogeneous Dirichlet boundary conditions for the velocity field at the reduced order level. 
The ROM computations are carried out using ITHACA-FV \cite{RoSta17}, 
an in-house open source C++ library. 

\subsection{A POD-Galerkin projection method}\label{sec:ROM_1}
The main idea of reduced order modeling for parametrized PDEs is the assumption that 
solutions live in a low dimensional manifold. Thus, any solution can be approximated as a 
linear combination of a reduced number of global basis functions. 

We approximate velocity fields $\v$ and $\u$ and pressure fields $q$ and $\overline{q}$ 
as linear combinations of the dominant modes (basis functions), which are assumed to be dependent on space variables only,
multiplied by scalar coefficients that depend on the time and/or the parameters:
\begin{align}
\v \approx \v_r = \sum_{i=1}^{N_v^r} \beta_i(\bm{\pi}, t) \bm{\varphi}_i(\bm{x}), \quad 
q \approx q_r = \sum_{i=1}^{N_q^r} \gamma_i(\bm{\pi}, t) \psi_i(\bm{x}), \label{eq:ROM_1} \\
\u \approx \u_r = \sum_{i=1}^{N_u^r} \overline{\beta_i}(\bm{\pi}, t) \bm{\overline{\varphi}}_i(\bm{x}), \quad
\overline{q} \approx \overline{q}_r = \sum_{i=1}^{N_{\overline{q}}^r} \overline{\gamma_i}(\bm{\pi}, t) \overline{\psi}_i(\bm{x}). \label{eq:ROM_2}
\end{align}
In \eqref{eq:ROM_1}-\eqref{eq:ROM_2}, $N_{\Phi}^r$ denotes the cardinality of a reduced basis for the space field $\Phi$ belongs to.

We point out that, unlike previous works \cite{Gunzburger2019b, Wells2017,Xie2018,Gunzburger2019}, 
we compute also pressure fields $q$ and $\overline{q}$.
Furthermore, we consider different basis and different coefficients for the approximation of the velocity fields $\v$ and $\u$. 
This follows from the fact that we apply the filtering step for both FOM and ROM. 

In the literature, one can find several techniques to generate the reduced basis spaces, e.g.~Proper Orthogonal Decomposition (POD), the Proper Generalized Decomposition (PGD) and the Reduced Basis (RB) with a greedy sampling strategy.
See, e.g., \cite{Rozza2008, ChinestaEnc2017, Kalashnikova_ROMcomprohtua, quarteroniRB2016, Chinesta2011, Dumon20111387, Huerta2020, ModelOrderReduction}. 
We find the reduced basis by using the method of snapshots.

Let $\mathcal{K}= \{\bm{\pi}^1, \dots, \bm{\pi}^{N_k}\}$ be a finite dimensional training set of samples chosen inside the parameter space $\mathcal{P}$ and for each time instance $t^k \in \{t^1, \dots, t^{N_t}\} \subset (t_0, T]$. We solve the FOM described in Sec.~\ref{sec:FOM}  
for each $\bm{\pi}^k \in \mathcal{K}  \subset \mathcal{P}$. The total number of snapshots $N_s$ is given by $N_s = N_k \cdot N_t$. 
The snapshots matrices are obtained from the full-order snapshots:  
\begin{align}\label{eq:space}
\bm{\mathcal{S}}_{{{\Phi}}} = [{{\Phi}}(\bm{\pi}^1, t^1), \dots, {{\Phi}}(\bm{\pi}^{N_k}, t^{N_t})] \in \mathbb{R}^{N_{\Phi}^h \times N_s} \quad
\text{for} \quad {{\Phi}} = \{\v_h, \u_h, q_h, \overline{q}_h\},
\end{align}
where the subscript $h$ denotes a solution computed with the FOM and $N_{\Phi}^h$
is the dimension of the space $\Phi$ belong to in the FOM. Note that ${\Phi}$ could be either a scalar or
a vector field. 
The POD problem consists in finding, for each value of the dimension of the POD space $N_{POD} = 1, \dots, N_s$, the scalar coefficients $a_1^1, \dots, a_1^{N_s}, \dots, a_{N_s}^1, \dots, a_{N_s}^{N_s}$ and functions ${\zeta}_1, \dots, {\zeta}_{N_s}$, that minimize the error between the snapshots and their projection onto the POD basis. In the $L^2$-norm, we have
\begin{align}
E_{N_{POD}} = \text{arg min} \sum_{i=1}^{N_s} ||{{\Phi}_i} - \sum_{k=1}^{N_{POD}} a_i^k {\zeta}_k || \quad \forall N_{POD} = 1, \dots, N_s    \cl
\text{with} \quad ({\zeta}_i, {\zeta}_j)_{L_2(\Omega)} = \delta_{i,j} \quad \forall i,j = 1, \dots, N_s. \label{eq:min_prob}
\end{align}

It can be shown \cite{Kunisch2002492} that eq.~\eqref{eq:min_prob} is equivalent to the following eigenvalue problem
\begin{align}
\bm{\mathcal{C}}^{{\Phi}} \bm{Q}^{{\Phi}} &= \bm{Q}^{{\Phi}} \bm{\Lambda}^{{\Phi}}, \label{eq:eigen_prob} \\
\mathcal{C}_{ij}^\Phi &= ({\Phi}_i, {\Phi}_j)_{L_2(\Omega)} \quad \text{for} \quad i,j = 1, \dots, N_s,
\end{align}
where $\bm{\mathcal{C}}^{{\Phi}}$ is the correlation matrix computed from the snapshot matrix $\bm{\mathcal{S}}_{{{\Phi}}}$, $\bm{Q}^{{\Phi}}$ is the matrix of eigenvectors and $\bm{\Lambda}^{{\Phi}}$ is a diagonal matrix whose diagonal entries are the
eigenvalues of $\bm{\mathcal{C}}^{{\Phi}}$. 
Then, the basis functions are obtained as follows:
\begin{align}\label{eq:basis_func}
{\zeta}_i = \dfrac{1}{N_s \Lambda_i^\Phi} \sum_{j=1}^{N_s} {\Phi}_j Q_{ij}^\Phi.
\end{align}
The POD modes resulting from the aforementioned methodology are:
\begin{align}\label{eq:spaces}
L_\Phi = [\bm{\zeta}_1, \dots, \bm{\zeta}_{N_\Phi^r}] \in \mathbb{R}^{N_\Phi^h \times N_\Phi^r},
\end{align}
where $N_\Phi^r < N_s$ are chosen according to the eigenvalue decay of the vectors of eigenvalues \bm{\Lambda}.
The reduced order model can be obtained through a Galerkin projection of the governing equations onto the POD spaces. 

Let
\begin{align}
&M_{r_{ij}} = (\bm{\varphi}_i, \bm{\varphi}_j)_{L_2(\Omega)}, \quad \widetilde{M}_{r_{ij}} = (\bm{\varphi}_i, \overline{\bm{\varphi}_j})_{L_2(\Omega)}, \quad A_{r_{ij}} = (\bm{\varphi}_i, \Delta \bm{\varphi}_j)_{L_2(\Omega)}, \label{eq:matrices_evolve1} \\
&B_{r_{ij}} = (\bm{\varphi}_i, \nabla \psi_j)_{L_2(\Omega)}, \quad P_{r_{ij}} = (\psi_i, \nabla \cdot \bm{\varphi}_j)_{L_2(\Omega)},
\label{eq:matrices_evolve2}
\end{align}
where $\bm{\varphi}_i$ and $\psi_i$ are the basis functions in \eqref{eq:ROM_1}. The reduced algebraic system at time $t^{n+1}$
for problem \eqref{eq:evolve-1.1}-\eqref{eq:evolve-1.2} is: 
\begin{align}
&\rho\, \frac{3}{2\Delta t} \bm{M}_r \bm{\beta}^{n+1} + \rho \bm{G}_r(\overline{\bm{\beta}}^{n}, \overline{\bm{\beta}}^{n-1}) \bm{\beta}^{n+1} - 2\mu \bm{A}_r \bm{\beta}^{n+1} + \bm{B}_r \bm{\gamma}^{n+1} = \dfrac{\rho}{\Delta t} \widetilde{\bm{M}}_r \left(2\overline{\bm{\beta}}^{n} - \dfrac{1}{2}\overline{\bm{\beta}}^{n-1}\right), \label{eq:reduced_1} \\
&\bm{P}_r \bm{\beta}^{n+1} = 0, \label{eq:reduced_2}
\end{align}
where vectors $\bm{\beta}^{n+1}$ and $\bm{\gamma}^{n+1}$ contain the values of coefficients $\beta_i$ and $\gamma_i$ 
in \eqref{eq:ROM_1} at time $t^{n+1}$.
The term $\bm{G_r}(\overline{\bm{\beta}}^{n}, \overline{\bm{\beta}}^{n-1}) \bm{\beta}^{n+1}$ in \eqref{eq:reduced_1}
is related to the non-linear convective term: 
\begin{align}\label{eq:convective_matrix}
\left(\bm{G_r}(\overline{\bm{\beta}}^{n}, \overline{\bm{\beta}}^{n-1}) \bm{\beta}^{n+1}\right)_i = (2\overline{\bm{\beta}}^{n} - \overline{\bm{\beta}}^{n-1})^T  \bm{\mathcal{G}}_{r_{i\bm{..}}} \bm{\beta}^{n+1}
\end{align}
where $\bm{\mathcal{G}_r}$ is a third-order tensor defined as follows \cite{quarteroni2007numerical, Rozza2009}
\begin{align}\label{eq:C_tensor}
\mathcal{G}_{r_{ijk}} = (\bm{\varphi_i}, \nabla \cdot (\bm{\varphi_j} \otimes \overline{\bm{\varphi}_k}))_{L_2(\Omega)}.
\end{align}

Next, let 
\begin{align}
&{\overline{M}}_{r_{ij}} = (\bm{\overline{\varphi}}_i, \bm{\overline{\varphi}}_j)_{L_2(\Omega)}, \quad
\overline{A}_{r_{ij}} = (\bm{\overline{\varphi}}_i, \Delta \bm{\overline{\varphi}}_j)_{L_2(\Omega)}, \label{eq:matrices_filter1} \\
&\overline{B}_{r_{ij}} = (\bm{\overline{\varphi}}_i, \nabla \overline{\psi}_j)_{L_2(\Omega)}, \quad
\overline{P}_{r_{ij}} = (\overline{\psi}_i, \nabla \cdot \bm{\overline{\varphi}}_j)_{L_2(\Omega)},    \label{eq:matrices_filter2}
\end{align}
where $\overline{\bm{\varphi}}_i$ and $\overline{\psi}_i$ are the basis functions in \eqref{eq:ROM_2}. 
The reduced algebraic system at time $t^{n+1}$
for problem \eqref{eq:filter-1.1},\eqref{eq:filter-1.2} is
\begin{align}
&\frac{\rho}{\Delta t} \bm{{\overline{M}}}_r \bm{\overline{\beta}}^{n+1}  - \mubar \bm{\overline{A}}_r \bm{\overline{\beta}}^{n+1} + \bm{\overline{B}}_r\bm{\overline{\gamma}}^{n+1} = \frac{\rho}{\Delta t} \bm{\widetilde{M}}^T_r \bm{\beta}^{n+1} \label{eq:reduced2_1}, \\
&\bm{\overline{P}}_r \bm{\overline{\beta}}^{n+1} = 0. \label{eq:reduced2_2}
\end{align}
where vectors $\overline{\bm{\beta}}^{n+1}$ and $\overline{\bm{\gamma}}^{n+1}$ contain the values of coefficients $\overline{\beta}_i$ and $\overline{\gamma}_i$ in \eqref{eq:ROM_2} at time $t^{n+1}$.

The complete reduced algebraic system at time $t^{n+1}$ is given by \eqref{eq:reduced_1}-\eqref{eq:reduced_2},\eqref{eq:reduced2_1}-\eqref{eq:reduced2_2}.
Finally, the initial conditions for the ROM algebraic system \eqref{eq:reduced_1}-\eqref{eq:reduced_2}
are obtained performing a Galerkin projection of the initial full order condition onto the POD basis spaces:
\begin{align}
{\beta^0}_i = (\v(\bm{x},\bm{\pi},t_0), \bm{\varphi}_i)_{L_2(\Omega)}, \cl
{\overline{\beta}^0}_i = (\u(\bm{x},\bm{\pi},t_0), \bm{\overline{\varphi}}_i)_{L_2(\Omega)}. \el
\end{align}

\subsection{Pressure fields reconstruction and pressure stability}\label{sec:ROM_2}

For the accurate reconstruction of the pressure field at the reduced level, different approaches have been proposed. 
One option is to use a global POD basis for both pressure and velocity field and same temporal coefficients
\cite{Bergmann2009Enablers,Lorenzi2016}. Another option to satisfy the inf-sup (or Ladyzhenskaya-Brezzi-Babuska) 
condition \cite{BREZZI199027, boffi_mixed} is through the supremizer enrichment method \cite{Rozza2007,Ballarin2014,Stabile2018}.
Finally, one can take the divergence of the momentum equation to obtain a Poisson equation for pressure that is projected onto a POD basis
\cite{akhtar2009stability}. This method, called Poisson pressure equation (PPE), has been recently extended to a 
finite volume context \cite{Stabile2017, Stabile2018, StabileZancanaroRozza2020, Star2019}. 

We choose to adopt and extend the PPE method used in \cite{Stabile2017, Stabile2018, Star2019} for both pressure fields in the EF algorithm.
We take the divergence of eq.~\eqref{eq:evolve-1.1} and \eqref{eq:filter-1.1}
and account for divergence free conditions~\eqref{eq:evolve-1.2} and \eqref{eq:filter-1.2}
to obtain the Poisson pressure equation:
\begin{align}
&\Delta q^{n+1} = -\rho\, \nabla \cdot \left(\div \left(\u^* \otimes \v^{n+1}\right)\right), \label{eq:system_def1_2}\\ 
&\Delta \overline{q}^{n+1} = 0, \label{eq:system_def1_4}
\end{align}
with boundary conditions \eqref{eq:bc-d_1}, \eqref{eq:bc-filter-d_1}, and:
\begin{align}
& \partial_{{n}} q^{n+1} = -2 \mu \bm{n} \cdot \left(\nabla \times \nabla \times \v^{n+1} \right) - \bm{n} \cdot \left(\rho\dfrac{3}{2\Delta t}\v^{n+1} - \bm{b}^{n+1}\right) \quad \mbox{on } \partial\Omega_N \times(t_0,T), \label{eq:system_def2_2}\\
& \partial_{{n}} \overline{q}^{n+1} = -2 \overline{\mu} \bm{n} \cdot \left(\nabla \times \nabla \times \u^{n+1} \right)  \quad \mbox{on } \partial\Omega_N \times(t_0,T), \label{eq:system_def2_4}
\end{align}
where $\partial_{{n}}$ denotes the derivative with respect to the normal vector \bm{n}. 
So, at reduced level, instead of solving \eqref{eq:evolve-1.1}-\eqref{eq:evolve-1.2} and \eqref{eq:filter-1.1},\eqref{eq:filter-1.2}, 
we solve modified systems \eqref{eq:evolve-1.1}, \eqref{eq:system_def1_2} and \eqref{eq:filter-1.1}, \eqref{eq:system_def1_4}.

For the enforcement of non-homogeneous Neumann conditions for the pressure fields, we
refer to \cite{Orszag1986, JOHNSTON2004221}.

\begin{rem}\label{rem2}
Systems \eqref{eq:evolve-1.1}-\eqref{eq:evolve-1.2} and \eqref{eq:filter-1.1},\eqref{eq:filter-1.2}
are not equivalent to systems \eqref{eq:evolve-1.1},\eqref{eq:system_def1_2} and \eqref{eq:filter-1.1},\eqref{eq:system_def1_4}
for steady flows \cite{Li2020, Orszag1986, JOHNSTON2004221}. However, as discussed in Remark \ref{rem:gen_Stokes}, 
the filter problem can be considered as a time-dependent Stokes problem. 
\end{rem}

The matrix form of eq.~\eqref{eq:system_def1_2} and \eqref{eq:system_def1_4} is obtained, as usual, after integrating 
by parts the Laplacian terms in the weak formulation and accounting for the boundary conditions. We obtain:
\begin{align}
&\bm{D}_r \bm{\gamma}^{n+1} + \rho \bm{J}_r(\overline{\bm{\beta}}^{n}, \overline{\bm{\beta}}^{n-1}) \bm{\beta}^{n+1} - 2\mu \bm{N}_r\bm{\beta}^{n+1} - \dfrac{\rho}{2\Delta t}\bm{F}_r = 0 \label{eq:reduced_q}, \\
&\overline{\bm{D}}_r \overline{\bm{\gamma}}^{n+1} - 2\overline{\mu} \overline{\bm{N}}_r \overline{\bm{\beta}}^{n+1} = 0, \label{eq:reduced_q_bar}
\end{align}
where
\begin{align}
&D_{r_{ij}} = (\nabla \psi_i, \nabla \psi_j)_{L_2(\Omega)}, \label{eq:pressure_matrices1} \\
& N_{r_{ij}} = (\bm{n} \times \nabla \psi_i, \nabla \times \bm{\phi}_j) _{L_2(\partial \Omega)}, \quad F_{r_{ij}} = (\psi_i, \bm{n} \cdot (3\bm{\varphi}_j - 4\overline{\bm{\varphi}}_j+ \overline{\bm{\varphi}}_j))_{L_2(\partial\Omega)}, \label{eq:pressure_matrices2} \\
& \overline{D}_{r_{ij}} =(\nabla \overline{\psi}_i, \nabla \overline{\psi}_j)_{L_2(\Omega)},  \quad \overline{N}_{r_{ij}} = (\bm{n} \times \nabla \overline{\psi}_i, \nabla \times \bm{\overline{\phi}_j}) _{L_2(\partial \Omega)}. \label{eq:pressure_matrices3}
\end{align}

The residual associated with the non-linear term in the equation \eqref{eq:reduced_q} is evaluated using the same strategy proposed 
for eq.~\eqref{eq:reduced_1}. We have
\begin{align}\label{eq:convective_matrix_2}
\left(\bm{J_r}(\overline{\bm{\beta}}^{n}, \overline{\bm{\beta}}^{n-1}) \bm{\beta}^{n+1}\right)_i = (2\overline{\bm{\beta}}^{n} - \overline{\bm{\beta}}^{n-1})^T  \bm{\mathcal{J}}_{r_{i\bm{..}}} \bm{\beta}^{n+1},
\end{align}
where $\bm{\mathcal{J}}_{r}$ is a third-order tensor defined as follows
\begin{align}\label{eq:J_tensor}
\mathcal{J}_{r_{ijk}} = (\nabla \psi_i, \nabla \cdot (\bm{\varphi}_j \otimes \overline{\bm{\varphi}}_k))_{L_2(\Omega)}.
\end{align}

Finally, the ROM algebraic system that has to be solved at every time step is \eqref{eq:reduced_1}, \eqref{eq:reduced_q}, \eqref{eq:reduced2_1}, \eqref{eq:reduced_q_bar}.

\subsection{Treatment of the Dirichlet boundary conditions: the lifting function method}\label{sec:ROM_3}

In order to homogeneize the velocity fields snapshots and make them independent on the boundary conditions,
we use the lifting function method \cite{Stabile2017}. The lifting functions, also called control functions,  
are problem-dependent: they have to be divergence free in order to retain 
the divergence-free property of the basis functions and they have to satisfy the boundary conditions of the FOM. The velocity snapshots are then modified according to:
\begin{align}\label{eq:homog1}
\v'_h = \v_h - \sum_{j=1}^{N_{BC}} v_{{BC}_j}(\bm{\pi}, t) \bm{\chi}_j(\bm{x}), \\
\u'_h = \u_h - \sum_{j=1}^{N_{BC}} u_{{BC}_j}(\bm{\pi}, t)\bm{\chi}_j(\bm{x}), 
\end{align}
where $N_{BC}$ is the number of non-homogeneous Dirichlet boundary conditions, $\bm{\chi} (\bm{x})$ 
are the control functions, and $v_{{BC}_j}$ and $u_{{BC}_j}$ are suitable temporal coefficients. 
The POD is applied to the snapshots satisfying the  homogeneous boundary conditions 
and then the boundary value is added back:
\begin{align}\label{eq:homog2}
\v_r= \sum_{j=1}^{N_{BC}} v_{{BC}_j}(\bm{\pi}, t)\bm{\chi}_j(\bm{x})  + \sum_{i=1}^{N_v^r} \beta_i(\bm{\pi}, t) \bm{\varphi}_i(\bm{x}), \\
\u_r= \sum_{j=1}^{N_{BC}} u_{{BC}_j}(\bm{\pi}, t)\bm{\chi}_j(\bm{x})  + \sum_{i=1}^{N_u^r} \overline{\beta_i}(\bm{\pi}, t) \overline{\bm{\varphi}_i}(\bm{x}).
\end{align}

\section{Numerical results}\label{sec:results}

We consider two well-known test cases \cite{John2004,turek1996}: 2D and 3D flow past a cylinder. 
We investigate the performance of the ROM model in the reconstruction of the time evolution of the flow field. 
For the 2D case, parametrization of the filtering radius is also introduced. 

\subsection{2D flow past a cylinder}
This benchmark has been thoroughly investigated at FOM level in a finite volume environment in \cite{Girfoglio2019}. 
To the best of our knowledge, it is the first time that this benchmark is considered within a ROM framework. 

\begin{figure}[h]
\centering
\includegraphics[height=0.115\textwidth]{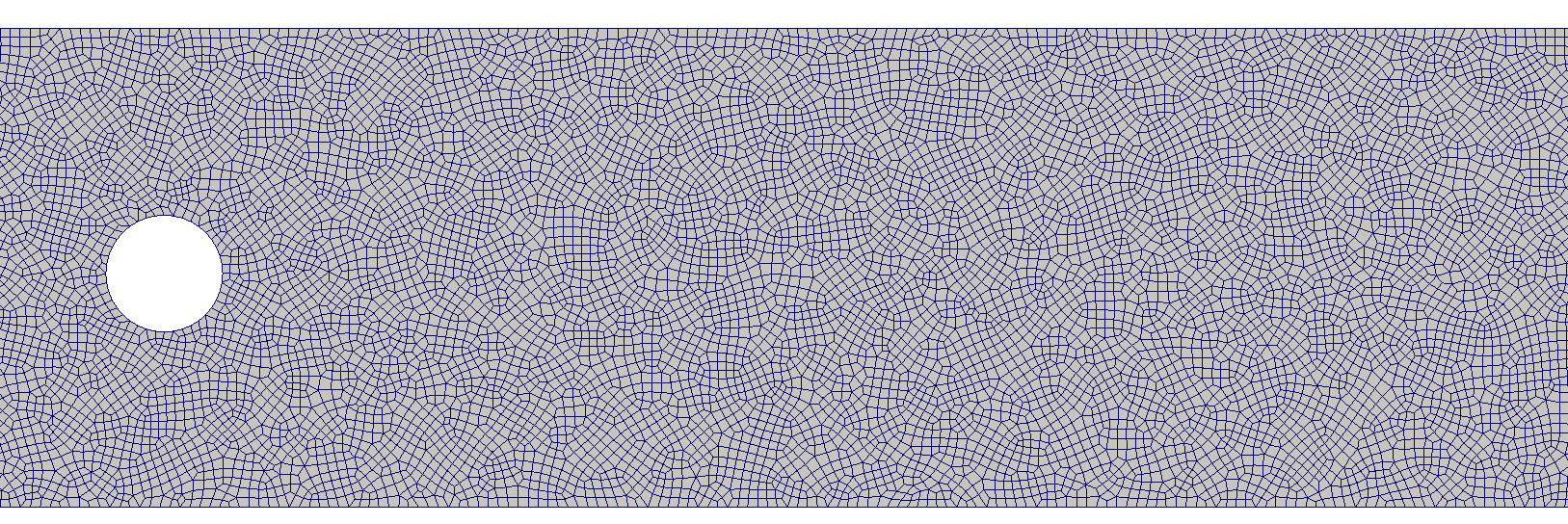}
\includegraphics[height=0.11\textwidth]{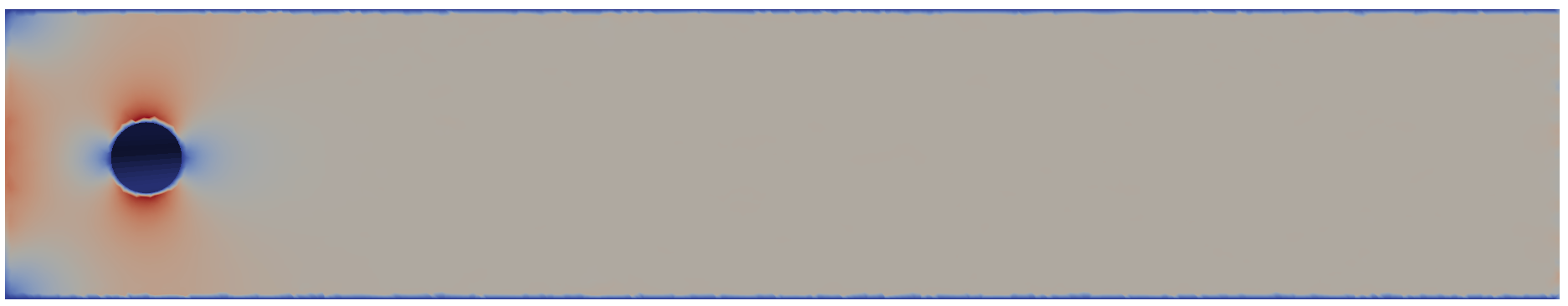}
\includegraphics[height=0.11\textwidth]{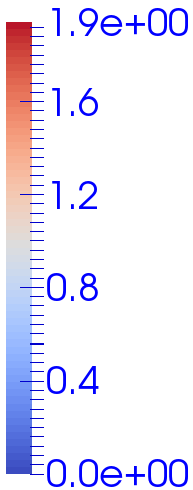}
\caption{2D flow past a cylinder: (left) part of the mesh under consideration and (right) the lifting function for velocity.}
\label{fig:example_cyl}
\end{figure}

The computational domain is a 2.2 $\times$ 0.41 rectangular channel with a cylinder of radius 0.05 centered at (0.2, 0.2), 
when taking the bottom left corner of the channel as the origin of the axes. 
Fig.~\ref{fig:example_cyl} (left) shows part of the computational domain. 
The channel is filled with fluid with density $\rho = 1$ and viscosity $\mu = 10^{-3}$.
We impose a no slip boundary condition on the upper and lower wall and on the cylinder. At the inflow, we prescribe
the following velocity profile:
\begin{align}\label{eq:cyl_bc}
\v(0,y,t) = \left(\dfrac{6}{0.41^2} \sin\left(\pi t/8 \right) y \left(0.41 - y \right), 0\right), \quad y \in [0, 0.41], \quad t \in (0, 8],
\end{align}
and ${\partial q}/{\partial \n} = {\partial \overline{q}}/{\partial \n} = 0$. At the outflow we prescribe $\nabla \v \cdot  \n = 0$ and $q = \qbar = 0$. 
We start the simulations from fluid at rest. Note that the Reynolds number is time dependent, with $0 \leq Re \leq 100$ \cite{turek1996}.

We consider a hexaedral computational grid with $h_{min} =  4.2e-3$ and $h_{max} = 1.1e-2$
for a total of $1.59e4$ cells. The quality of this mesh is high: it features
very low values of maximum non-orthogonality (36$^\circ$), 
average non-orthogonality (4$^\circ$), 
skewnwss (0.7), and maximum aspect ratio (2). 
Fig.~\ref{fig:example_cyl} (left) shows a part of the mesh. 
We chose this mesh because it is the coarsest among all the meshes considered in \cite{Girfoglio2019}
and thus the most challenging for our filtering approach.

While in \cite{Girfoglio2019} the choice of the time step depended on the 
Courant-Friedrichs-Lewy number ($CFL_{max}$) set to $CFL_{max} = 0.2$, 
in this work we consider a fixed time step for sake of consistency with the ROM. 
In order to obtain comparable solutions, we set $\Delta t = 4e-4$, which allows 
to obtain $CFL_{max} = 0.2$ at the time when the velocity reaches its maximum value. 
The same time step is used for the temporal resolution of 
the ROM dynamical system. 
For the convective term, we use a second-order accurate Central Differencing (CD) scheme \cite{Lax1960}. In this
way, we avoid introducing stabilization and are able to assess the effects produced by the
filter. 

As far as we know, in a ROM-FV framework only test cases with steady 
\cite{Stabile2017, Stabile2018, Star2019} or time-dependent \cite{Star2019} uniform Dirichlet 
boundary conditions have been considered. For this benchmark, we need to enforce time-dependent 
non-uniform Dirichlet boundary conditions. 
For this purpose, we consider a divergence free function with the following non-uniform velocity distribution at the inflow:
\begin{equation}\label{eq:BC_lift}
\bm{\chi} (0, y) = \left(\dfrac{6}{0.41^2} y \left(0.41 - y \right), 0\right), \quad y \in [0, 2.2],
\end{equation}
and uniform null values the rest of the boundary. See Figure \ref{fig:example_cyl} (right).


We will compare our findings mainly with \cite{Stabile2017, Stabile2018, Star2019}, because these references
developed a NSE-ROM finite volume framework with a PPE stabilization method for the reconstruction of the pressure field. 
Moreover, the benchmarks presented in \cite{Stabile2017, Stabile2018, Star2019} share some features with the ones
we consider: a 2D flow past cylinder at $Re = \mathcal{O} (100)$ \cite{Stabile2017, Stabile2018}
(although with a steady and uniform inflow condition), and time-dependent uniform inflow boundary condition \cite{Star2019}
(although in a Y-junction flow problem).

We are going to investigate the accuracy of our reduced order model with respect to the time history of the velocity and pressure fields. 
We set $\alpha = 0.0032$ and refer to \cite{Girfoglio2019} for details on this choice.
The snapshots are collected every 0.01 s using an equispaced grid method in time. Therefore, the dimension of the correlation matrix 
$\bm{\mathcal{C}}^{{\Phi}}$ in \eqref{eq:eigen_prob} is $800 \times 800$ and $N_s^u = N_s^v = N_s^q = N_s^{\qbar}= 800$. 
Following \cite{GeorgakaStabileRozzaBluck2019}, we performed a convergence test as the number of snapshots increases. 
We multiplied the frequency at which the snapshots are collected
by 5 and 10, leading to $N_s^u = N_s^v = N_s^q = N_s^{\qbar} = 160$ and $N_s^u = N_s^v = N_s^q = N_s^{\qbar} = 80$, respectively.
We calculated the $L^2$ relative error: 
\begin{equation}\label{eq:error1}
E_{\Phi} = \dfrac{||\Phi_h(t, \pi) - \Phi_r(t, \pi)||_{L^2(\Omega)}}{||{\Phi_h}(t, \pi)||_{L^2(\Omega)}},
\end{equation}
where $\Phi_h$ and $\Phi_r$ are a particular field computed with the FOM (i.e., $\v_h$, $\u_h$, $q_h$ or $\qbar_h$) 
and the ROM (i.e., $\v_r$, $\u_r$, $q_r$ or $\qbar_r$), respectively. 
Moreover, we evaluate the relative error of the kinetic energy
\begin{equation}\label{eq:error2}
E_{K_{\Phi_h}} = \dfrac{|K_{\Phi_h}(t, \pi) - K_{\Phi_r}(t, \pi)|}{K_{\Phi_h}(t, \pi)},
\end{equation}
where $K_{{\Phi_h}}$ and $K_{{\Phi_r}}$ are the values of the kinetic energy computed by the FOM (for $\v_h$ or $\u_h$)
and by the ROM (for $\v_r$ or $\u_r$), respectively.

\begin{figure}
\centering
 \begin{overpic}[width=0.45\textwidth]{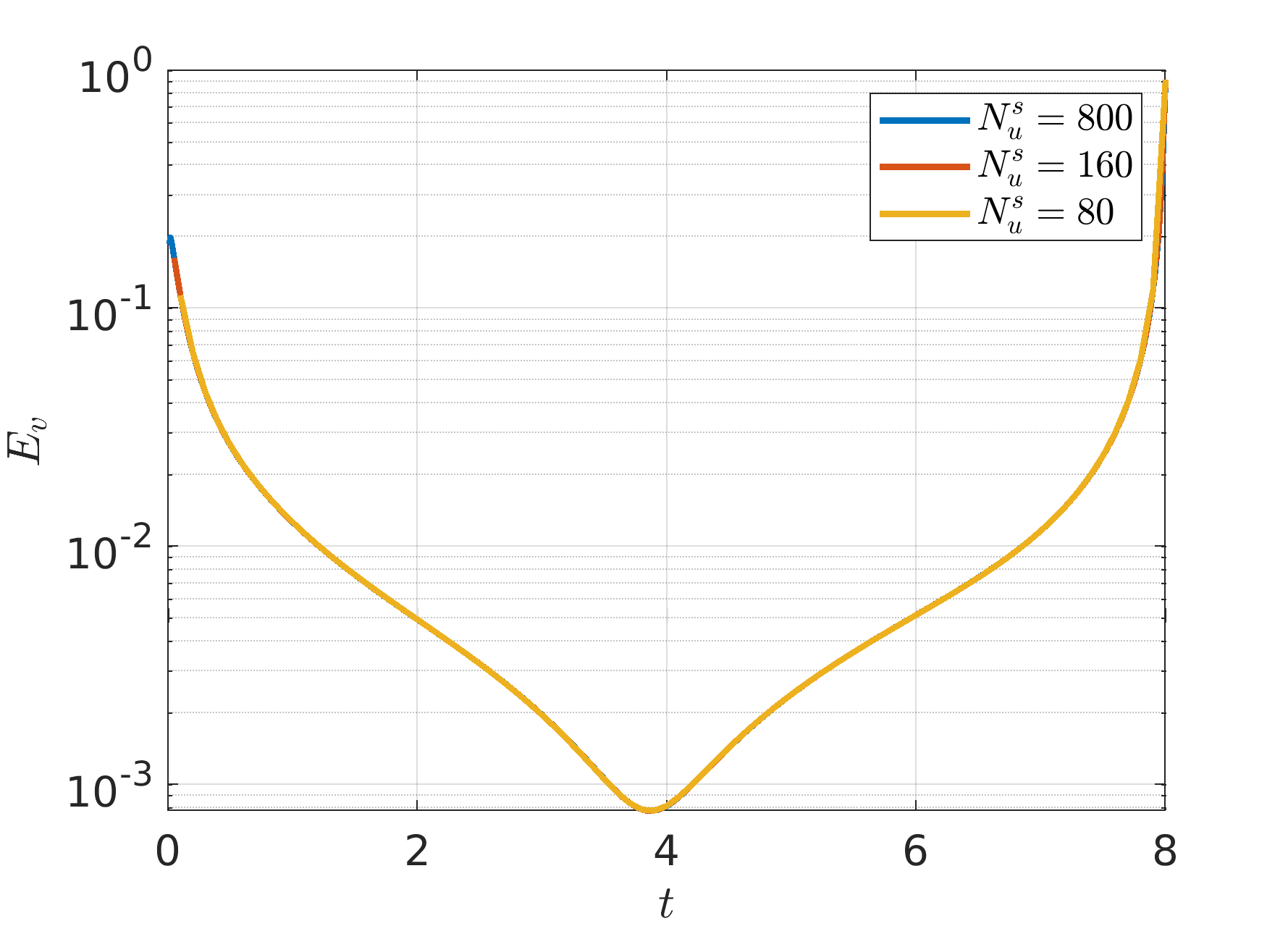}
      \end{overpic}
\begin{overpic}[width=0.45\textwidth]{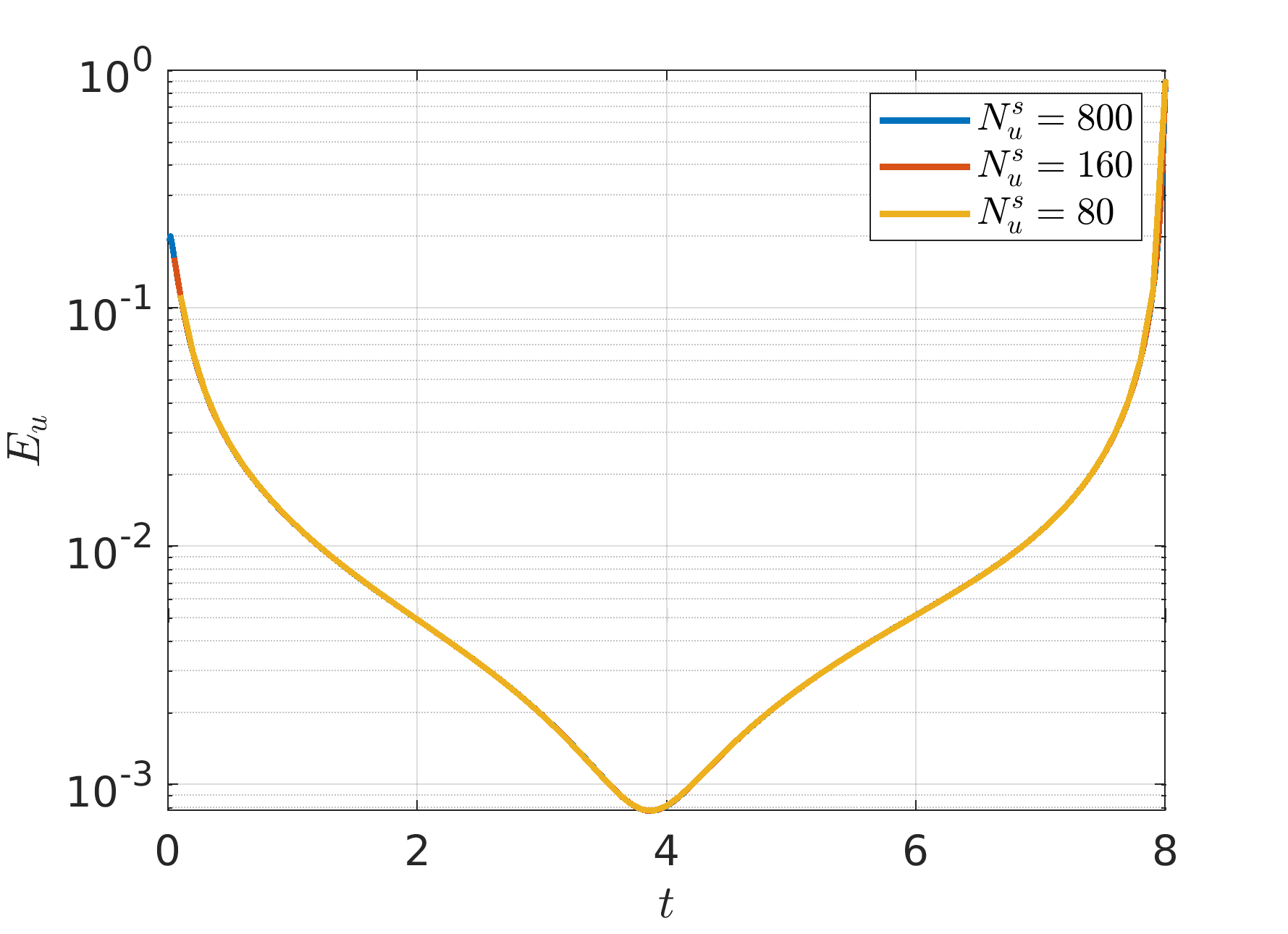}
      \end{overpic}\\
\begin{overpic}[width=0.45\textwidth]{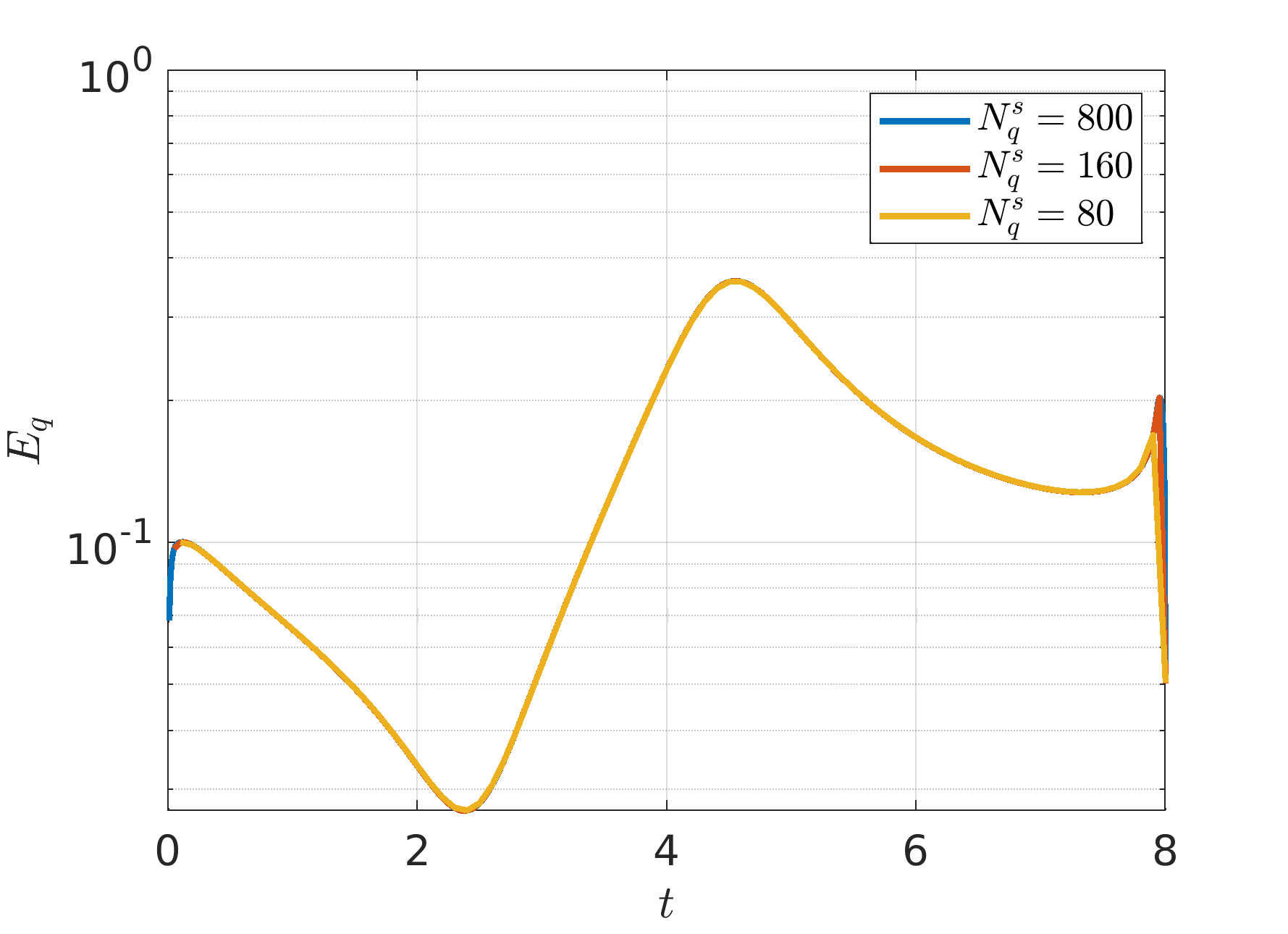}
      \end{overpic}
\begin{overpic}[width=0.45\textwidth]{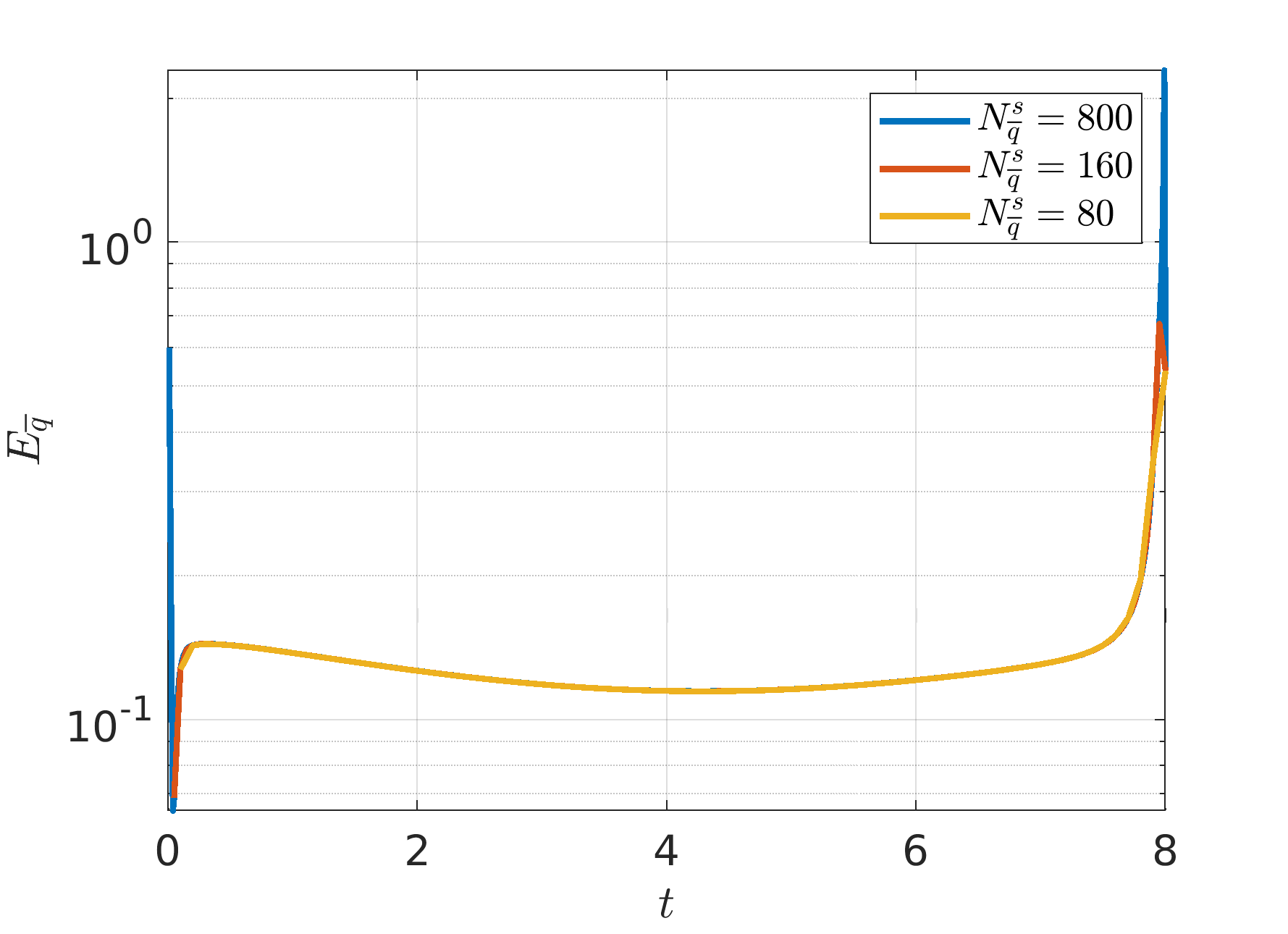}
      \end{overpic}\\
\begin{overpic}[width=0.45\textwidth]{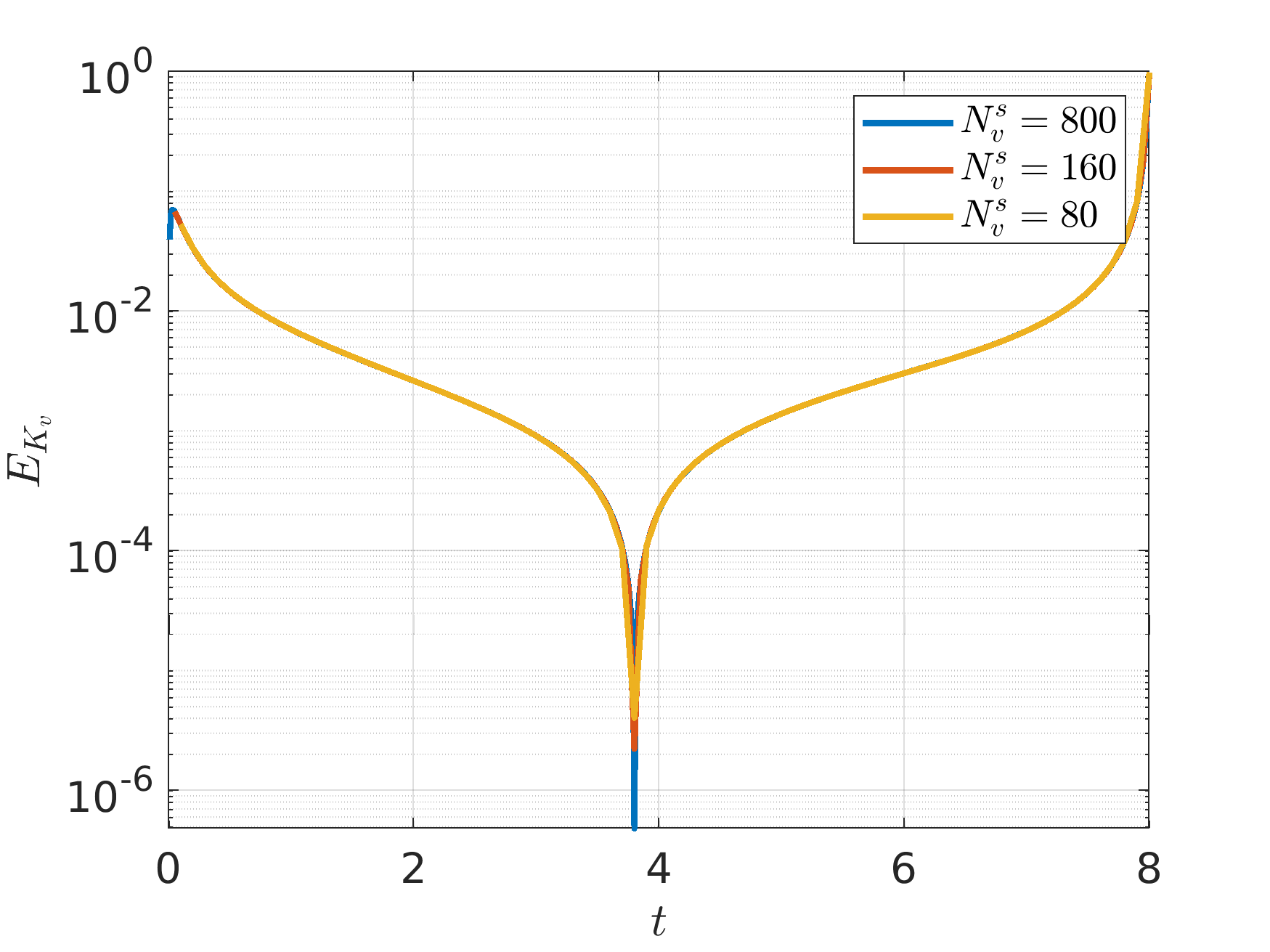}
      \end{overpic}
\begin{overpic}[width=0.45\textwidth]{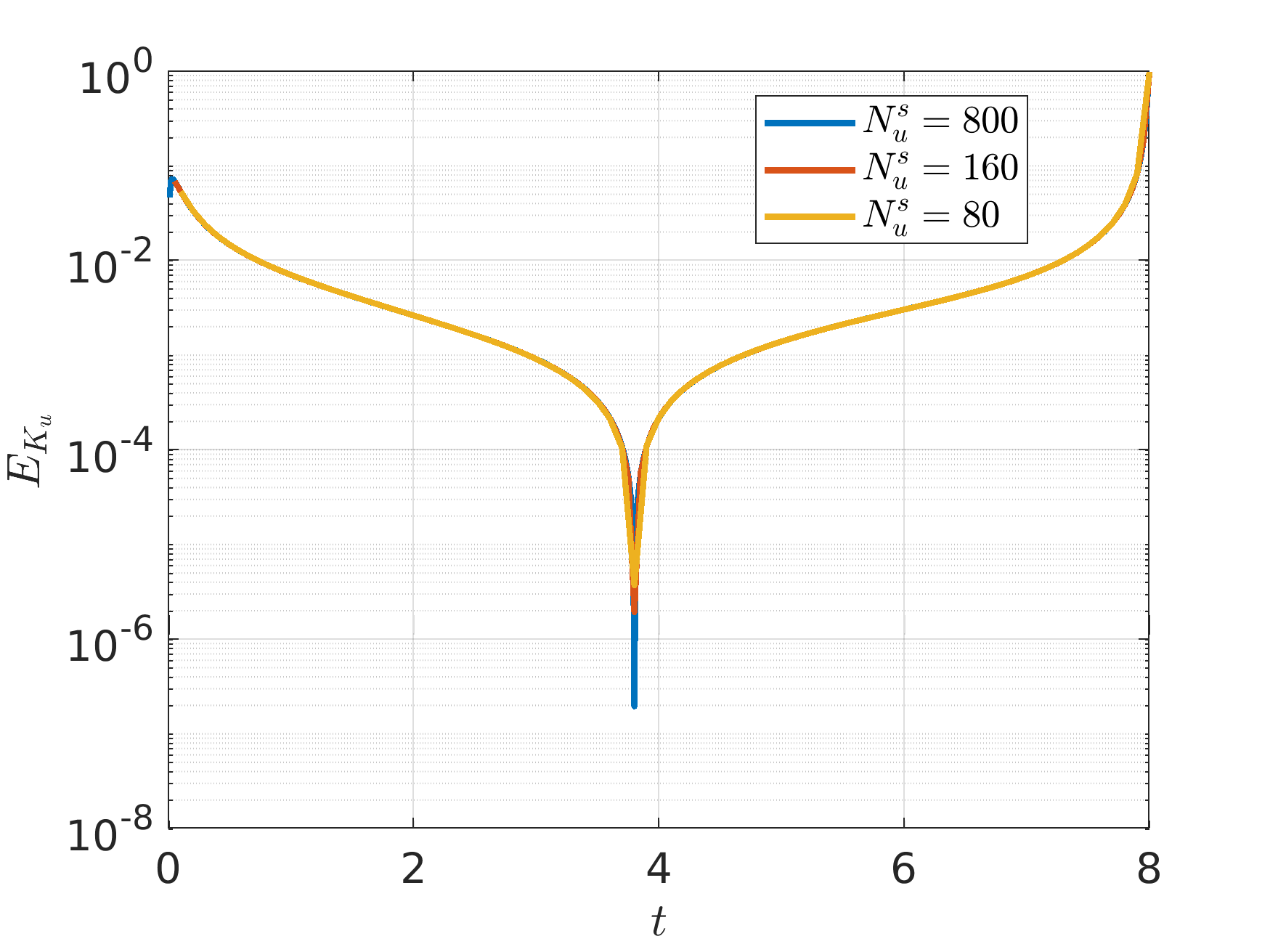}
      \end{overpic}\\
\caption{2D flow past a cylinder: time history of $L^2$ norm of the relative error for the velocity fields (top), pressure fields (center), and for kinetic energies of the system (bottom) for different numbers of snapshots. We consider 2 modes for $\v$, $\u$ and $q$, and 1 mode for $\bar{q}$.}
\label{fig:err_t}
\end{figure}

Fig.~\ref{fig:err_t} shows error \eqref{eq:error1} for all the velocity and pressure fields and error \eqref{eq:error2}
over time for the three different sampling frequencies. 
Fig.~\ref{fig:err_t} shows that, except for the initial and the final time of the simulation, 
the relative error for velocities fields is significantly lower than $10^{-1}$ over most of the time interval 
and reaches a minimum value of $7.8 \cdot 10^{-4}$.
Moreover, we notice that there is no substantial difference in the errors for the different sampling frequencies. Thus, 
to reduce the computational cost of the offline phase, we will consider the largest sampling frequency (i.e., 0.1) 
for the results presented from here on.
 
We report minimum, maximum and average (over time) relative errors
in Table \ref{tab:errors_t}.
We remark that the average errors for the velocity fields are comparable to the values 
(between $1.2 \cdot 10^{-2}$ and $3 \cdot 10^{-2}$) obtained for steady flow past cylinder 
$Re = 100$ \cite{Stabile2018}. Larger velocity errors during the first few 
time steps might be due to the transient nature of the flow,
as also noted for lid driven cavity flow studied in \cite{Star2019}. 
As for pressure  $q$, the relative error reaches its maximum value far from 
the end points 
of the time interval, while the average error for  $q$ is comparable to
the errors in \cite{Stabile2018} and \cite{Star2019}. 
The error for $\qbar$ stays close to the average during most of the time interval.
Finally, we observe that the smallest relative errors are achieved for the kinetic energies of the system.

\begin{table}
\centering
\begin{tabular}{lcccccc}
\multicolumn{2}{c}{} \\
\cline{1-7}
 & $\u$ & $\v$ & $q$ & $\overline{q}$ & $K_u$ & $K_v$  \\
\hline
Maximum $E_\Phi$ & 9.2e-1  & 9.1e-1 & 3.6e-1 & 5.4e-1 & 9.8e-1 & 9.8e-1 \\
Miminum $E_\Phi$ & 7.8e-4 & 7.8e-4 & 2.7e-2 & 1.1e-1 & 4e-6 & 3.7e-6\\
Average $E_\Phi$ & 2.4e-2 & 2.3e-2 & 1.4e-1 & 1.3e-1 & 1.9e-2 & 1.9e-2\\
\hline
\end{tabular}
\caption{2D flow past a cylinder: maximum, minimum and average (over time) 
relative errors for the velocity and pressure fields, and for the kinetic energies of the system. 
The sampling frequency of the snapshots is 0.1.}
\label{tab:errors_t}
\end{table}

\begin{figure}
\centering
 \begin{overpic}[width=0.45\textwidth]{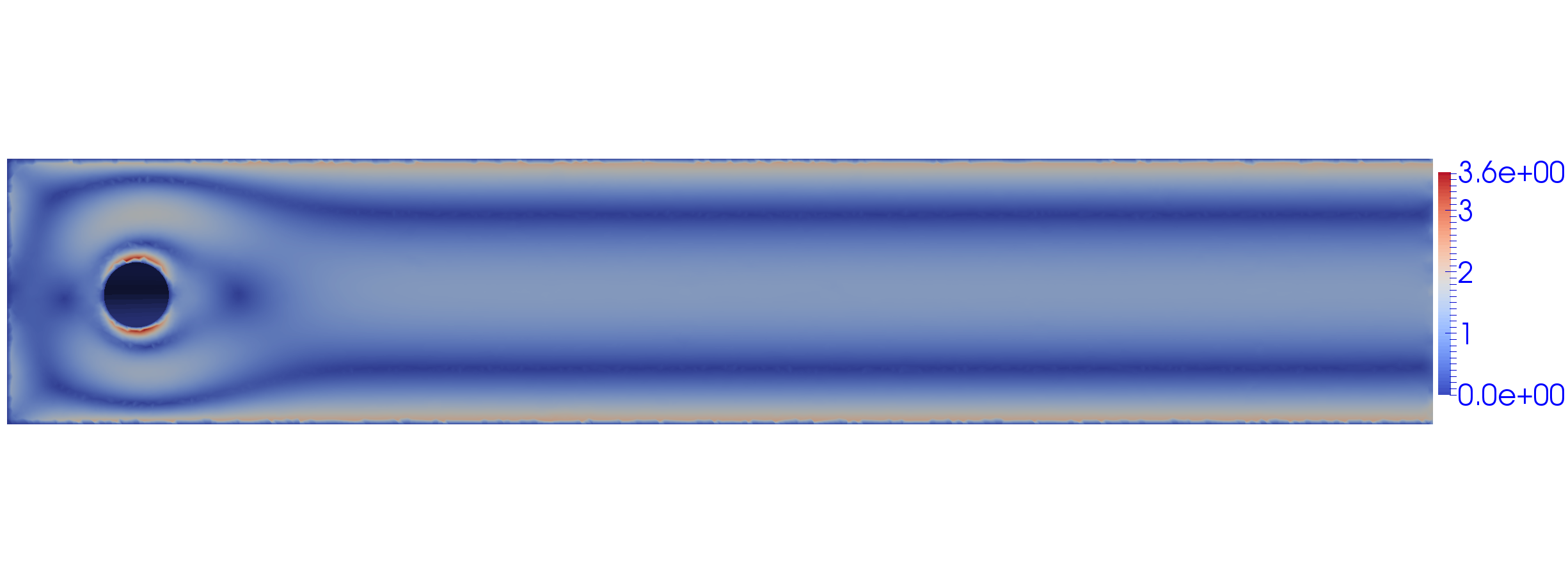}
        \put(30,30){\small{First mode for $\v$}}
      \end{overpic}
       \begin{overpic}[width=0.45\textwidth]{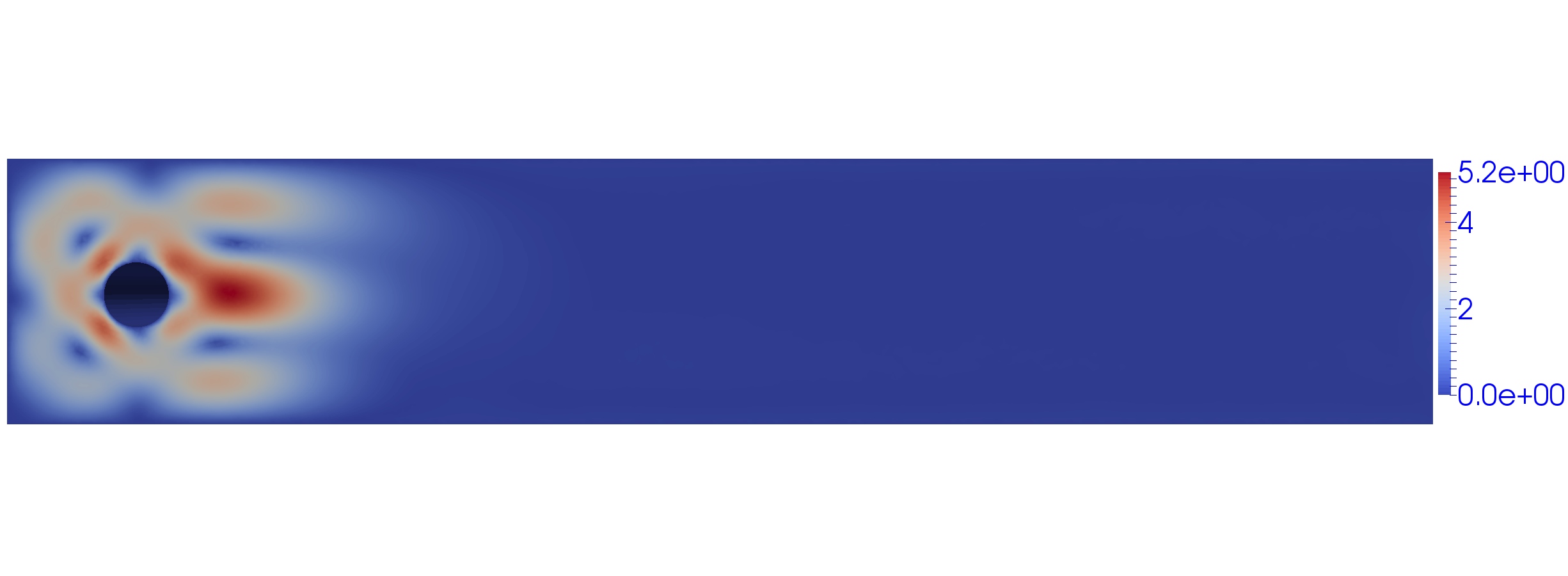}
        \put(30,30){\small{Second mode for $\v$}}
      \end{overpic}\\
 \begin{overpic}[width=0.45\textwidth]{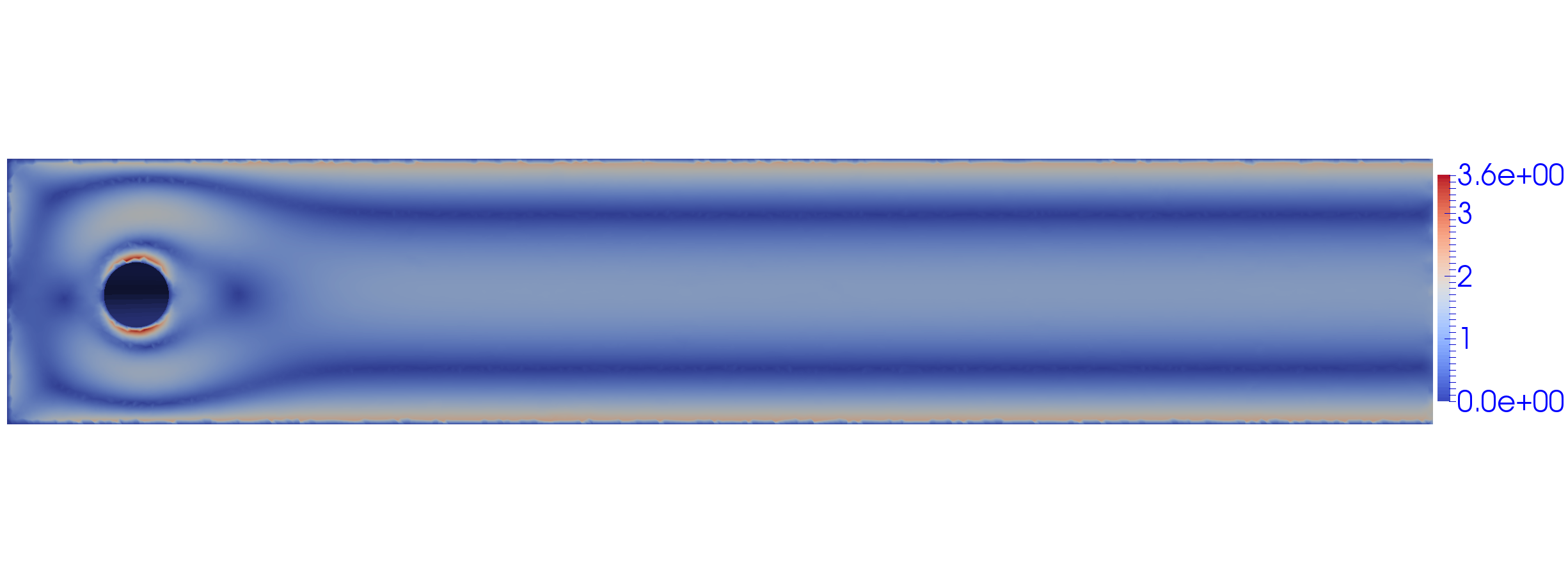}
        \put(30,30){\small{First mode for $\u$}}
      \end{overpic}
 \begin{overpic}[width=0.45\textwidth]{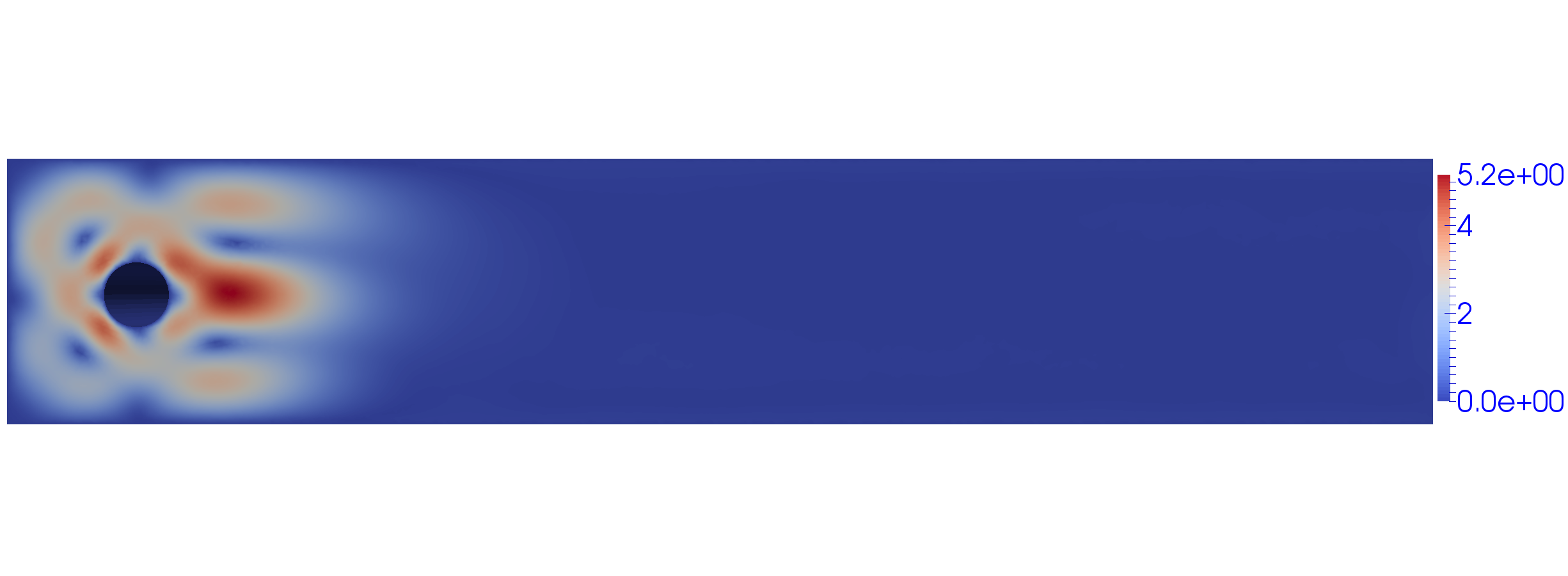}
        \put(30,30){\small{Second mode for $\u$}}
      \end{overpic}\\
 \begin{overpic}[width=0.45\textwidth]{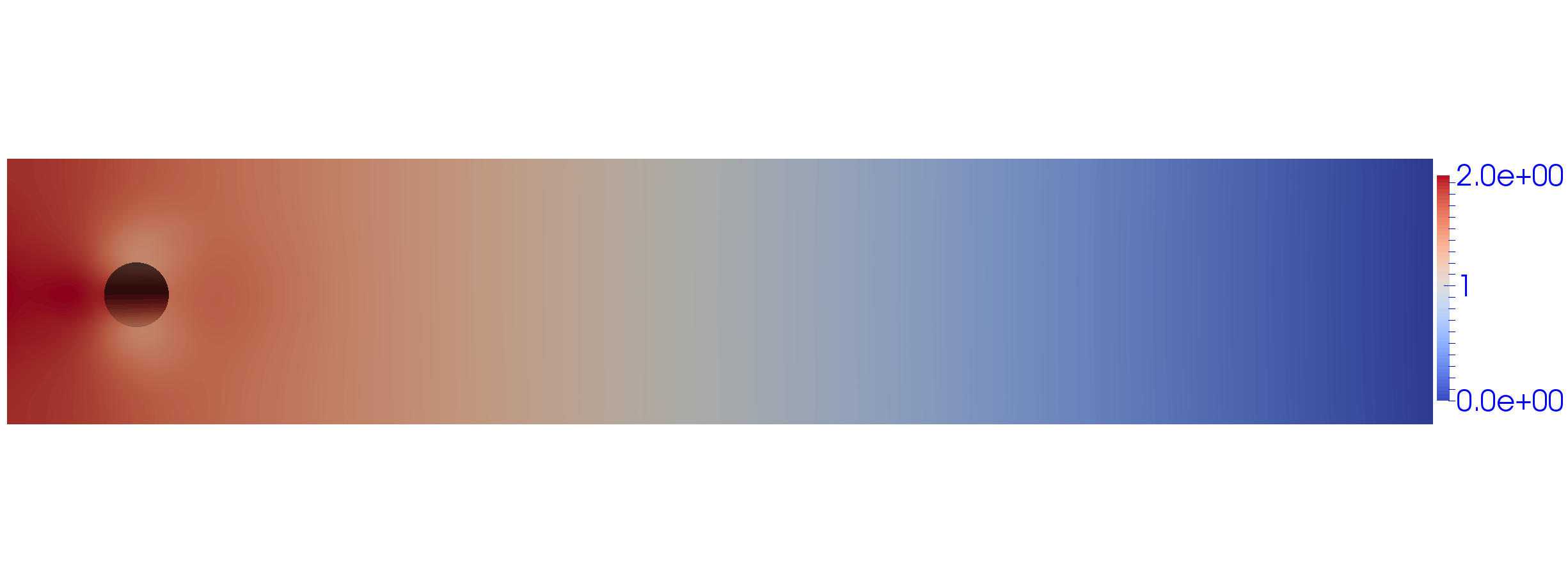}
\put(30,30){\small{First mode for $q$}}
      \end{overpic}
 \begin{overpic}[width=0.45\textwidth]{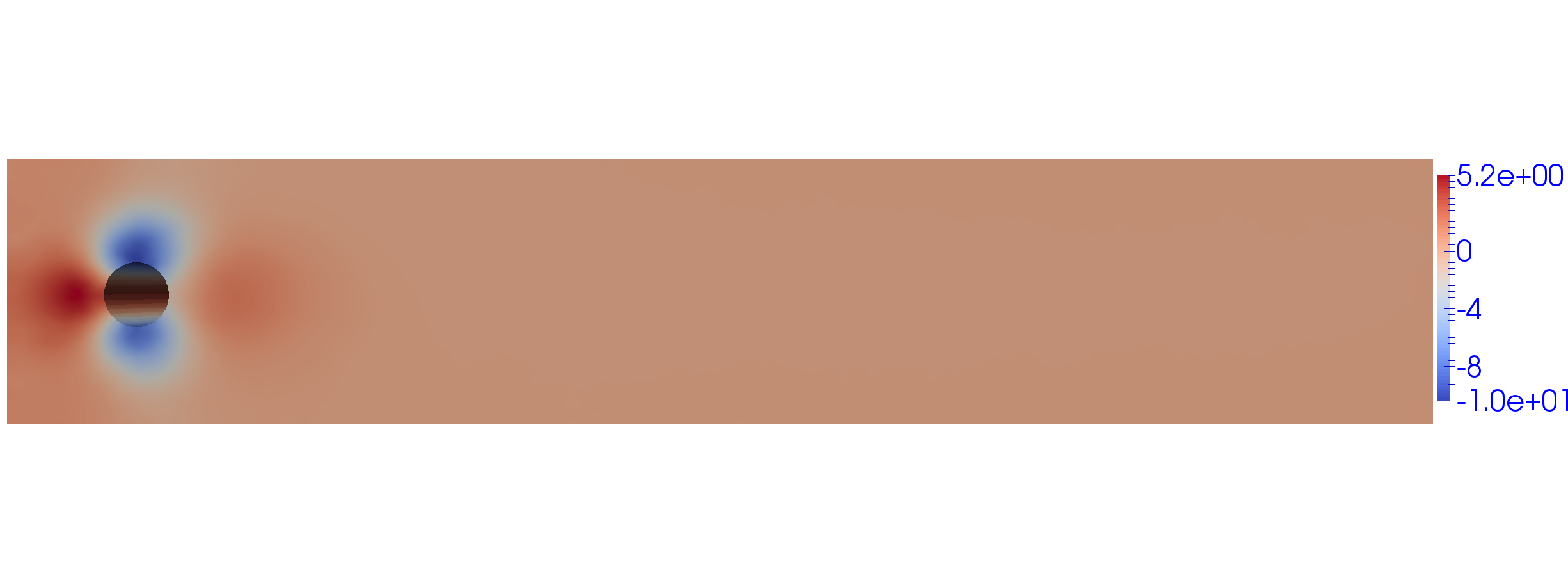}
\put(30,30){\small{Second mode for $q$}}
      \end{overpic}\\
 \begin{overpic}[width=0.45\textwidth]{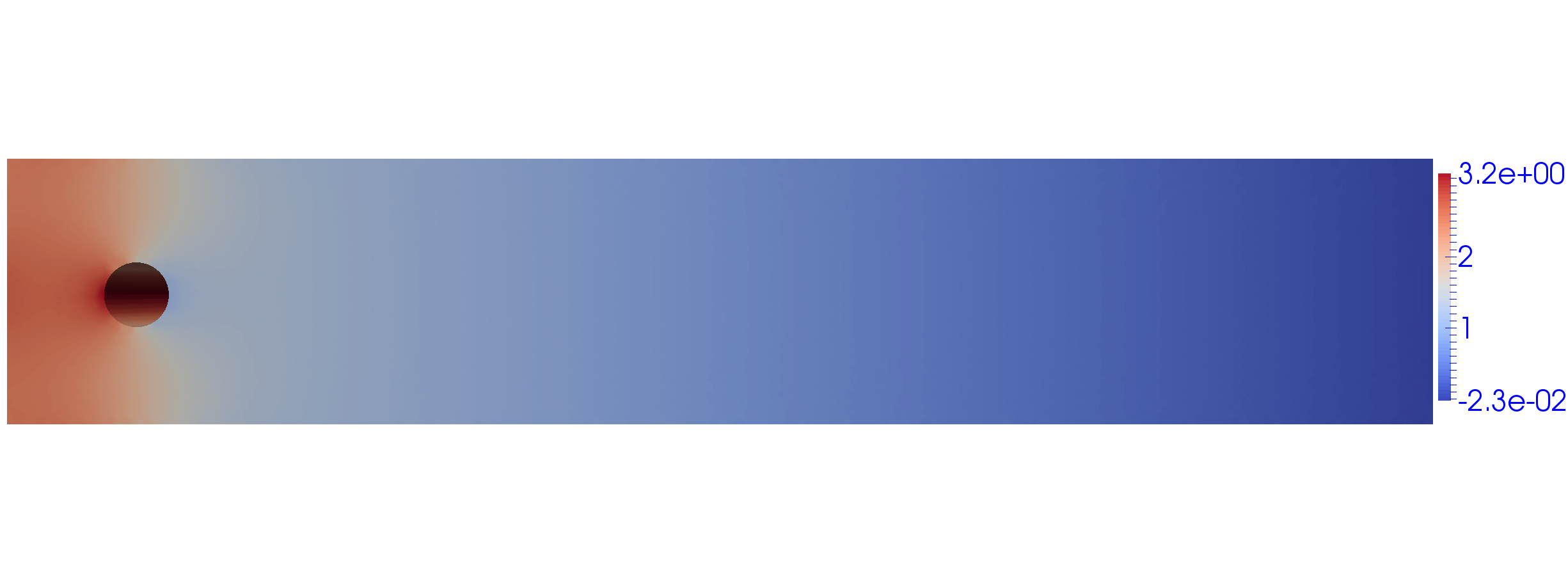}
\put(30,30){\small{First mode for $\overline{q}$}}
      \end{overpic}\\
\caption{2D flow past a cylinder: first 2 POD modes for velocity $\v$ (1nd raw), velocity $\u$ (2nd raw), pressure $q$ (3nd raw) and first POD mode for pressure $\overline{q}$ (4nd raw).}
\label{fig:modes2D}
\end{figure}

Table \ref{tab:cum2D} lists the first 3 cumulative eigenvalues, based on the first 15 most energetic POD modes. We see
that 2 modes for $\u$, $\v$ and $q$, and 1 mode for $\overline{q},$ are sufficient to retain 99.99\% of the energy contained in the snapshots. It has been verified that adding more modes does not increase the accuracy of the ROM results.
These first (homogenized) velocity and pressure modes are plotted in Figure \ref{fig:modes2D}.  
 
 \begin{table}
\centering
\begin{tabular}{lcccc}
\multicolumn{2}{c}{} \\
\cline{1-5}
N modes & $\u$ & $\v$ & $q$ & $\overline{q}$ \\
\hline
 1 & 0.999588 & 0.999582 & 0.967431 & 0.999985 \\
 2 & 0.999924 & 0.999924 & 0.999916 & 0.999997 \\
 3 & 0.999998 & 0.999998 & 0.999995 & 0.999999 \\
\hline
\end{tabular}
\caption{2D flow past a cylinder: first 3 cumulative eigenvalues for the velocity and pressure fields.}
\label{tab:cum2D}
\end{table}

\begin{figure}
\centering
       \begin{overpic}[width=0.45\textwidth]{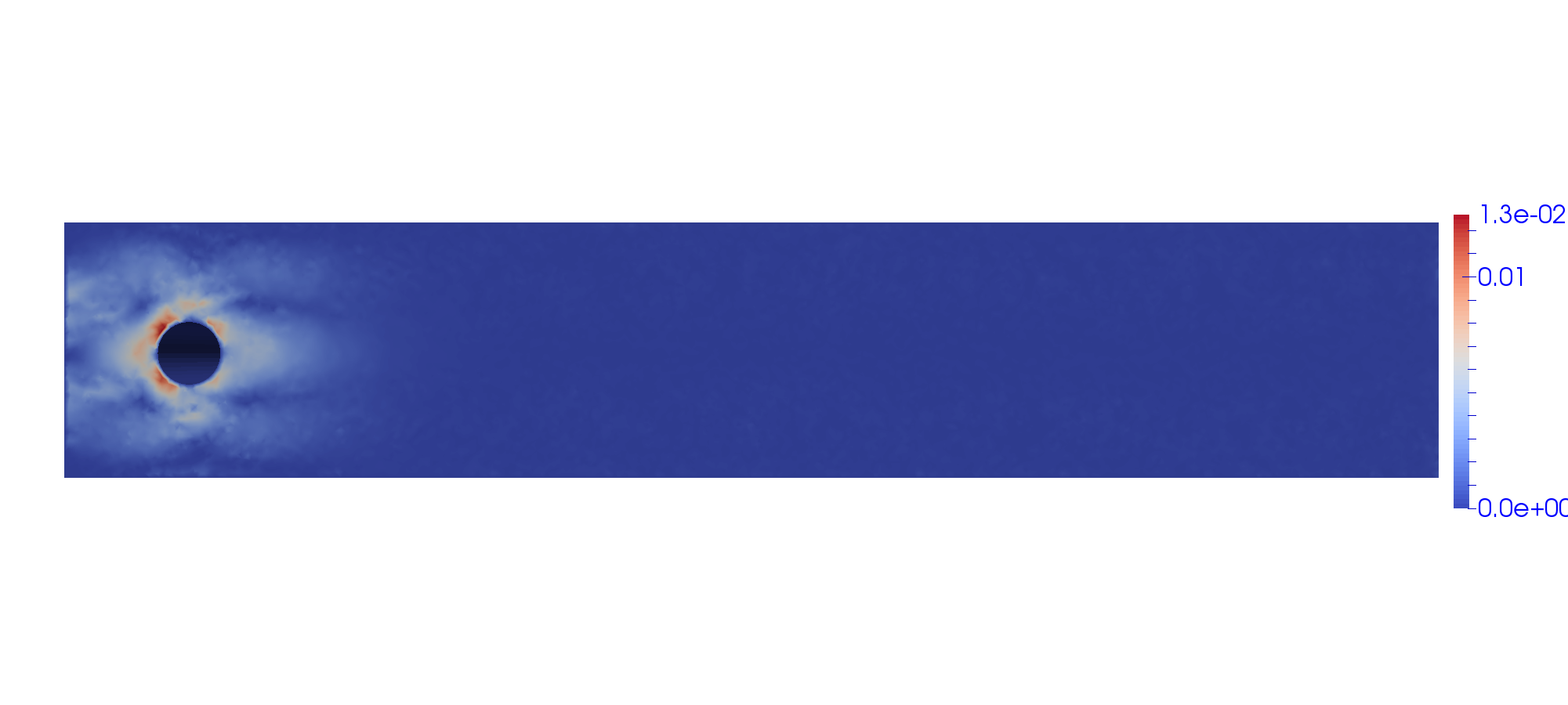}
      \end{overpic}
 \begin{overpic}[width=0.45\textwidth]{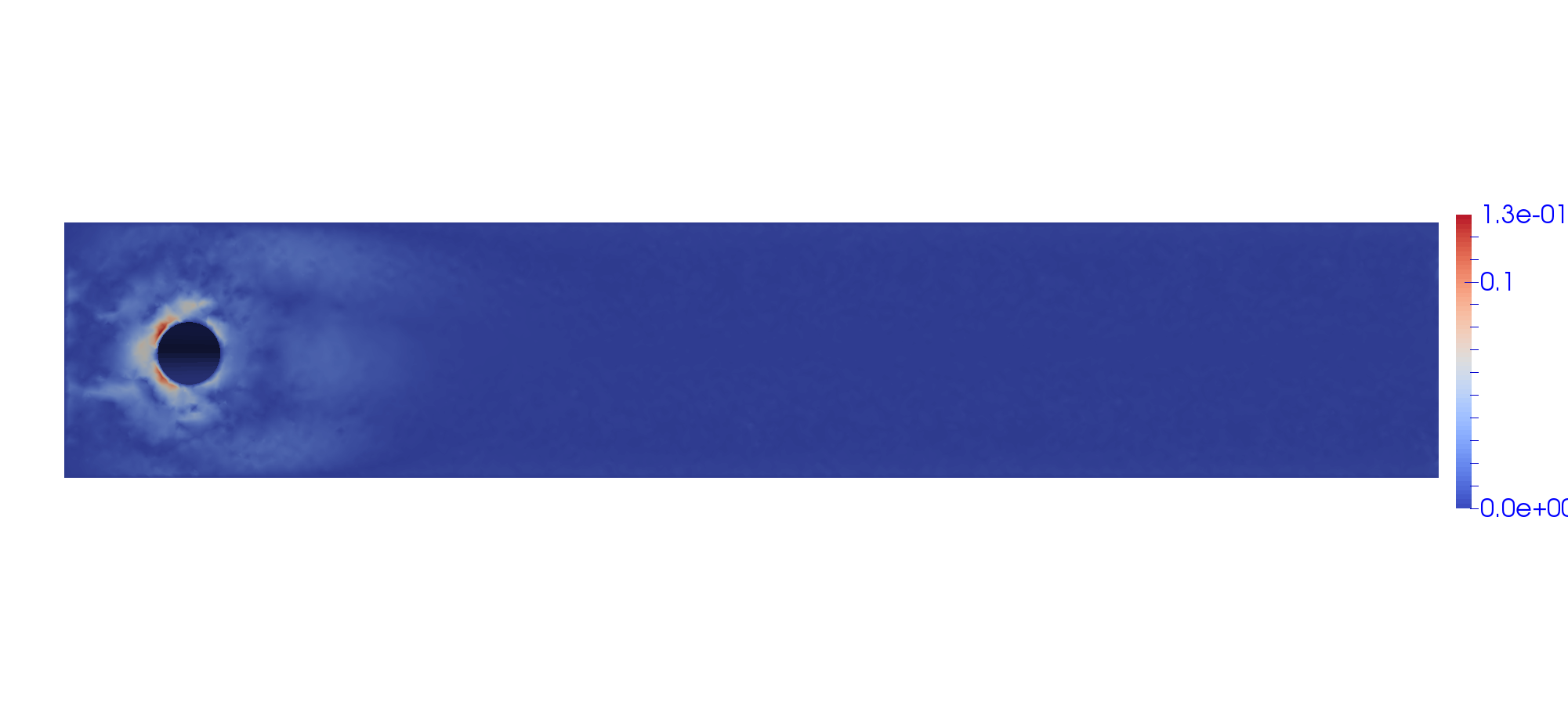}
      \end{overpic} \\
\caption{2D flow past a cylinder: difference between the first (left) and the second (right) POD modes related to the velocities $\u$ and $\v$. 
}
\label{fig:compUUev}
\end{figure}

\begin{rem}
Fig.~\ref{fig:compUUev} displays the difference between the first two POD modes for velocities $\u$ and $\v$. We note that there are significant differences right next to the cylinder where the filtering stabilization plays a key role. This justifies the fact to consider different bases to approximate the two velocity fields. 
\end{rem}

\begin{figure}
\centering
 \begin{overpic}[width=0.45\textwidth]{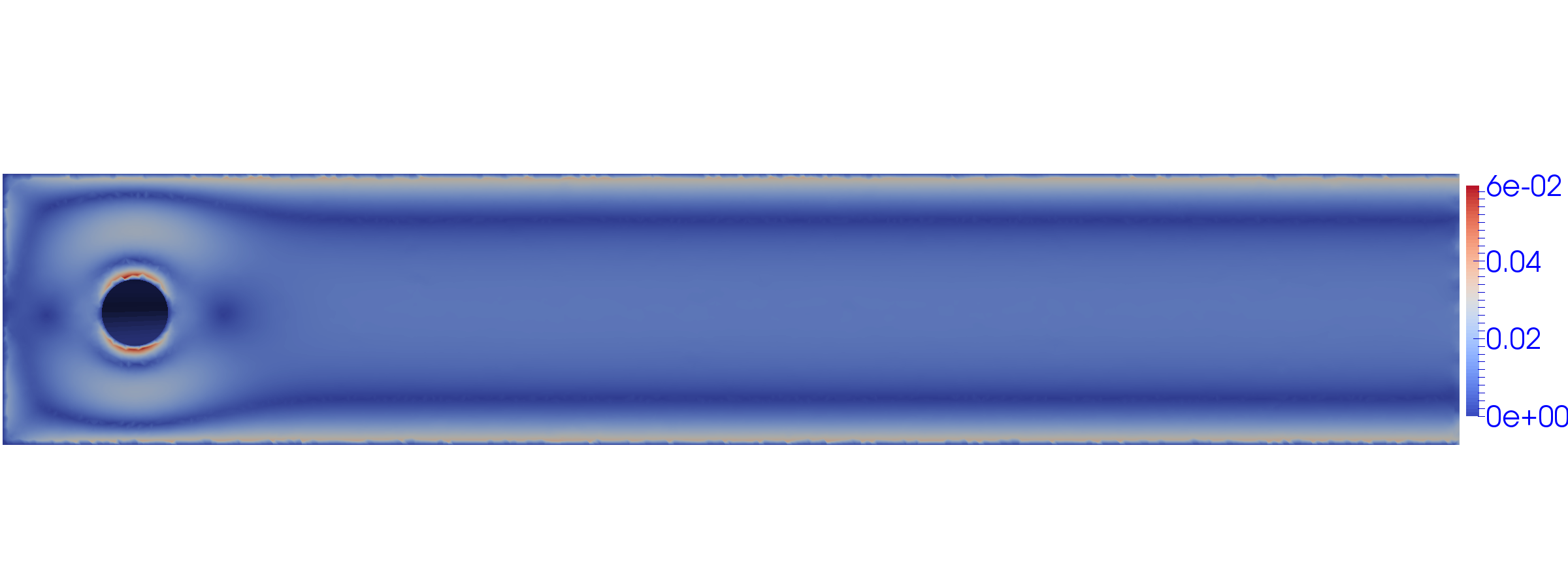}
        \put(40,30){$\u_{FOM}$}
      \end{overpic}
 \begin{overpic}[width=0.45\textwidth]{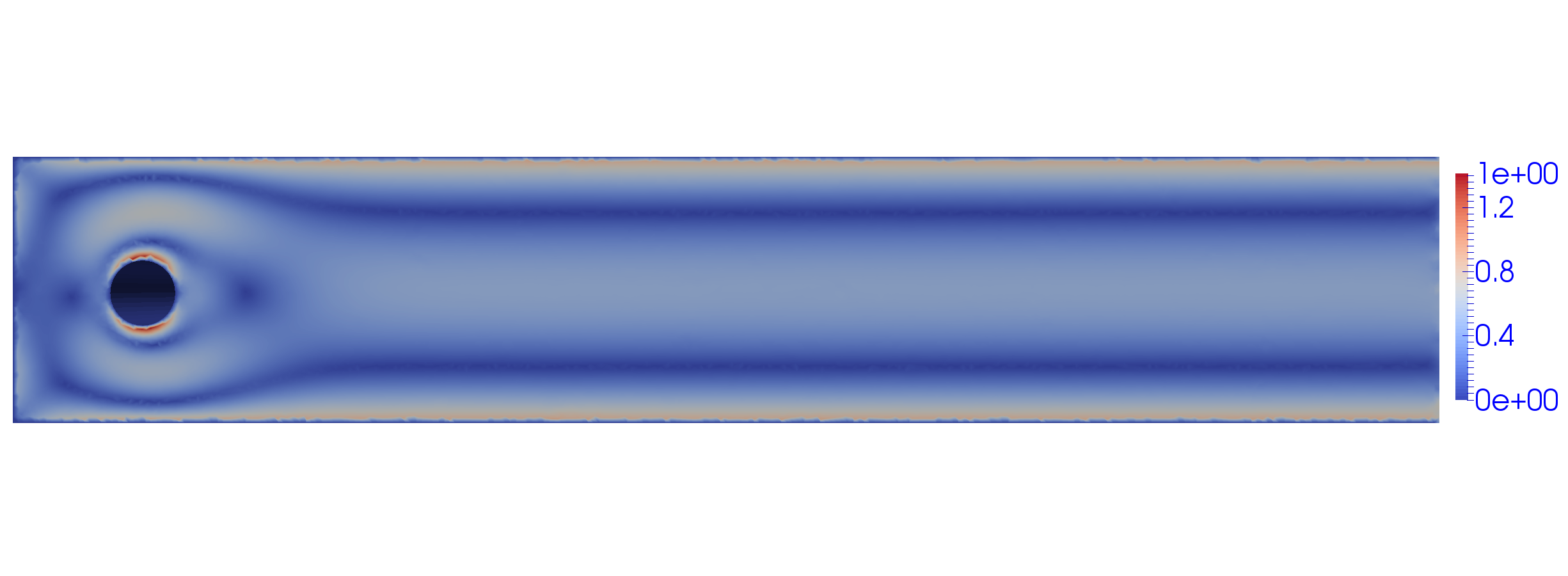}
        \put(40,30){$\u_{FOM}$}
      \end{overpic}\\
 \begin{overpic}[width=0.45\textwidth]{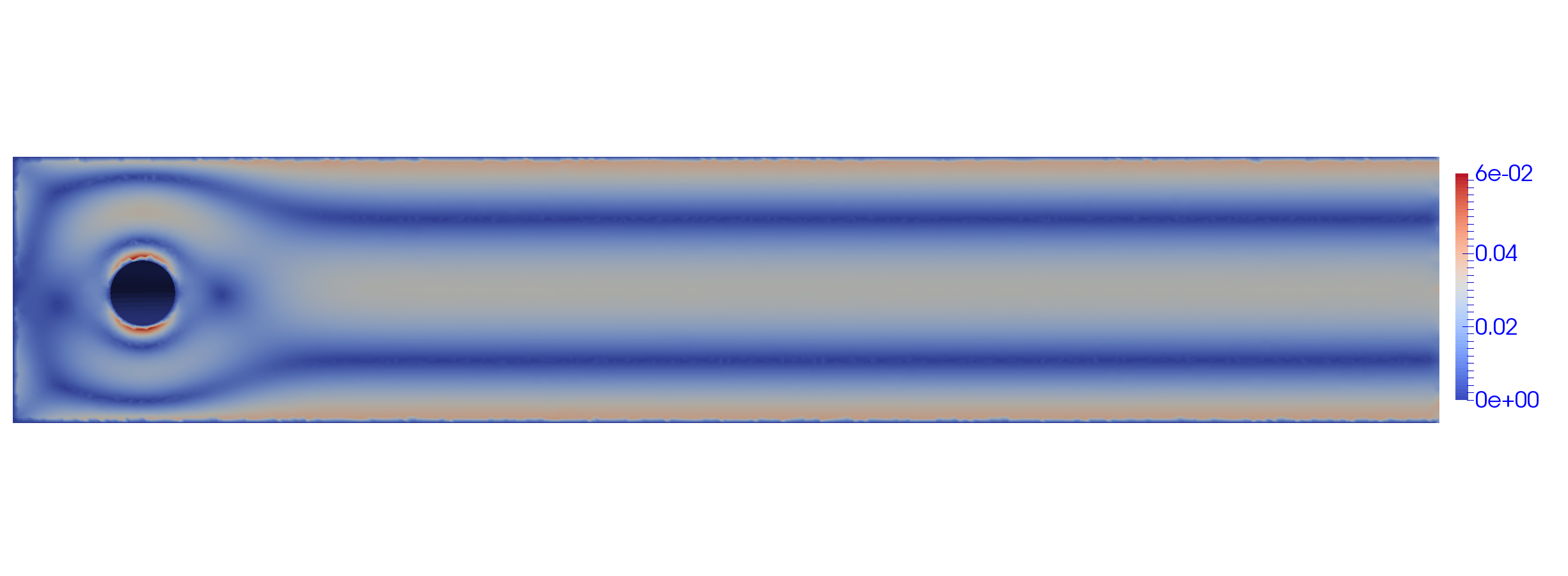}
        \put(40,30){$\u_{FOM}$}
      \end{overpic}
 \begin{overpic}[width=0.45\textwidth]{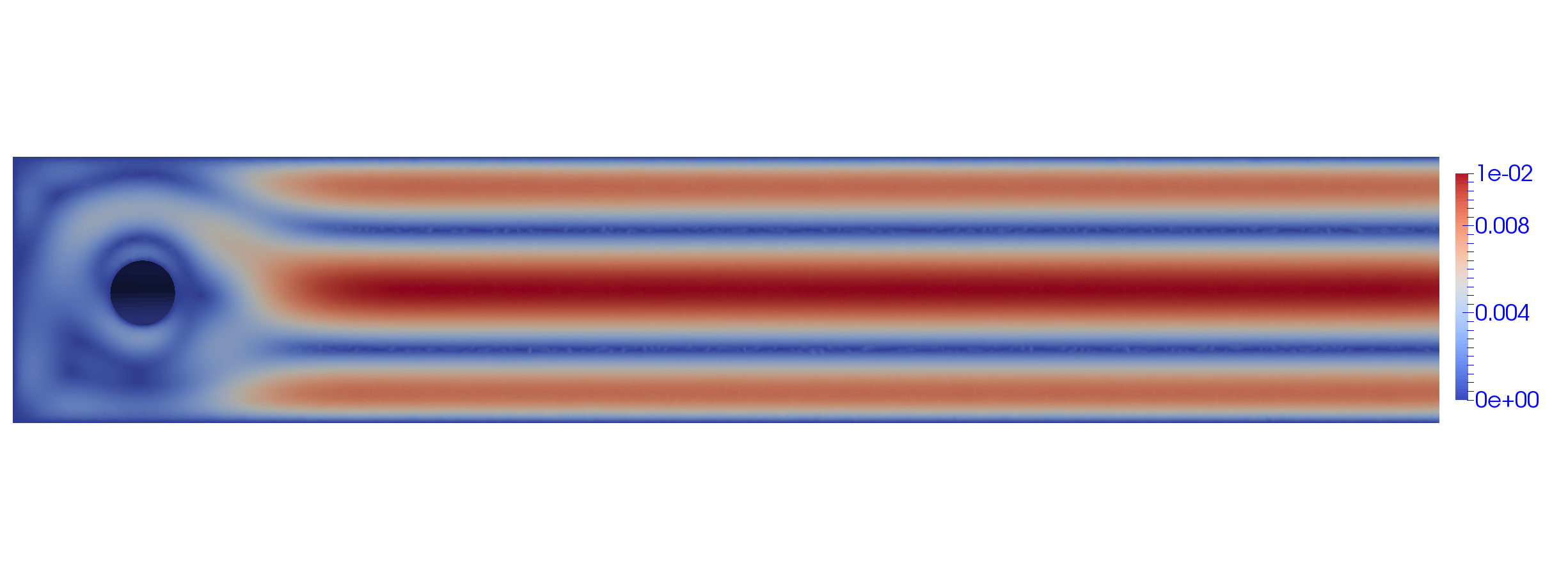}
        \put(40,30){$\u_{FOM}$}
      \end{overpic}\\
\caption{2D flow past a cylinder: FOM (homogeneized) velocity field $\u$ at $t = 0.1$, $5$ (first raw, from left to right), $7.9$ and $8$ (second raw, from left to right).}
\label{fig:2D_hom}
\end{figure}

We report in Fig.~\ref{fig:2D_hom} the homogeneized velocity $\u$ computed at times $t = 0.1, 5, 7.9, 8$. 
We observe that up to $t = 7.9$ the spatial structure of the flow field is very similar to the most energetic POD velocity mode (see Fig.~\ref{fig:modes2D}, first panel in the second row) whilst for $t = 8$ it is not comparable with either the first or the second POD mode. 
This could explain the larger relative errors for $\u$ towards the end of the time interval shown in Fig.~\ref{fig:err_t}.

\begin{figure}
\centering
       \begin{overpic}[width=0.45\textwidth]{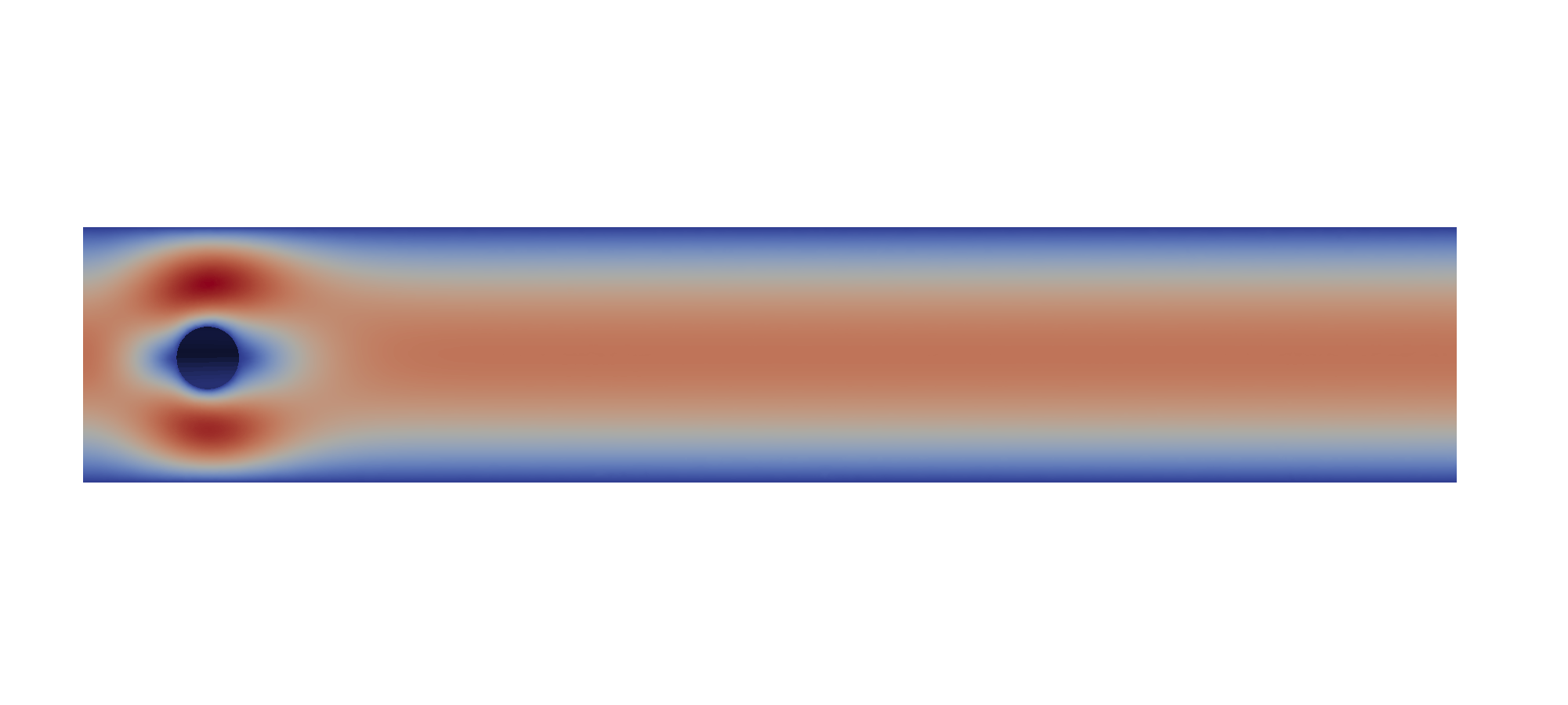}
        \put(40,35){$\u_{FOM}$}
      \end{overpic}
       \begin{overpic}[width=0.45\textwidth]{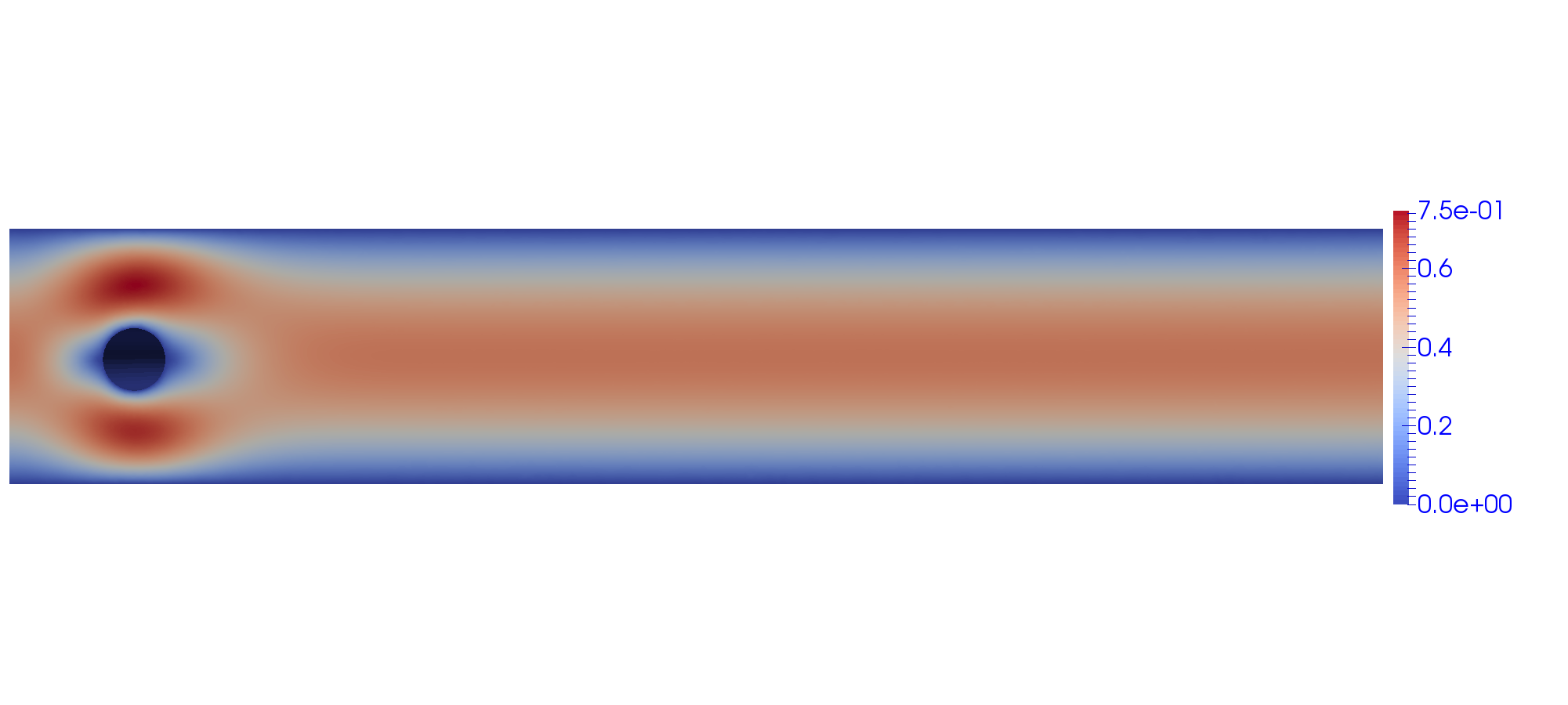}
        \put(40,35){$\u_{ROM}$}
      \end{overpic} \\
       \begin{overpic}[width=0.45\textwidth]{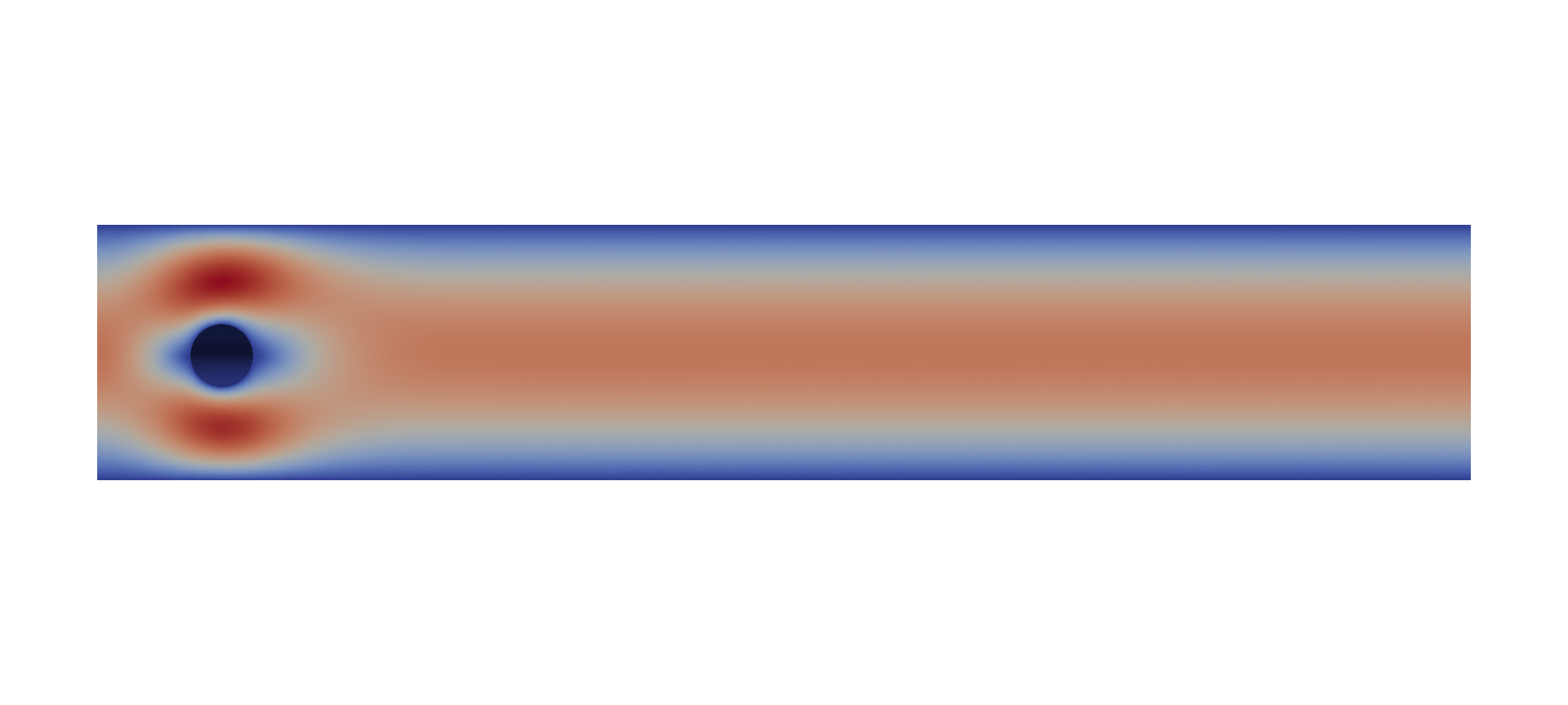}
        \put(40,35){$\v_{FOM}$}
      \end{overpic}
       \begin{overpic}[width=0.45\textwidth]{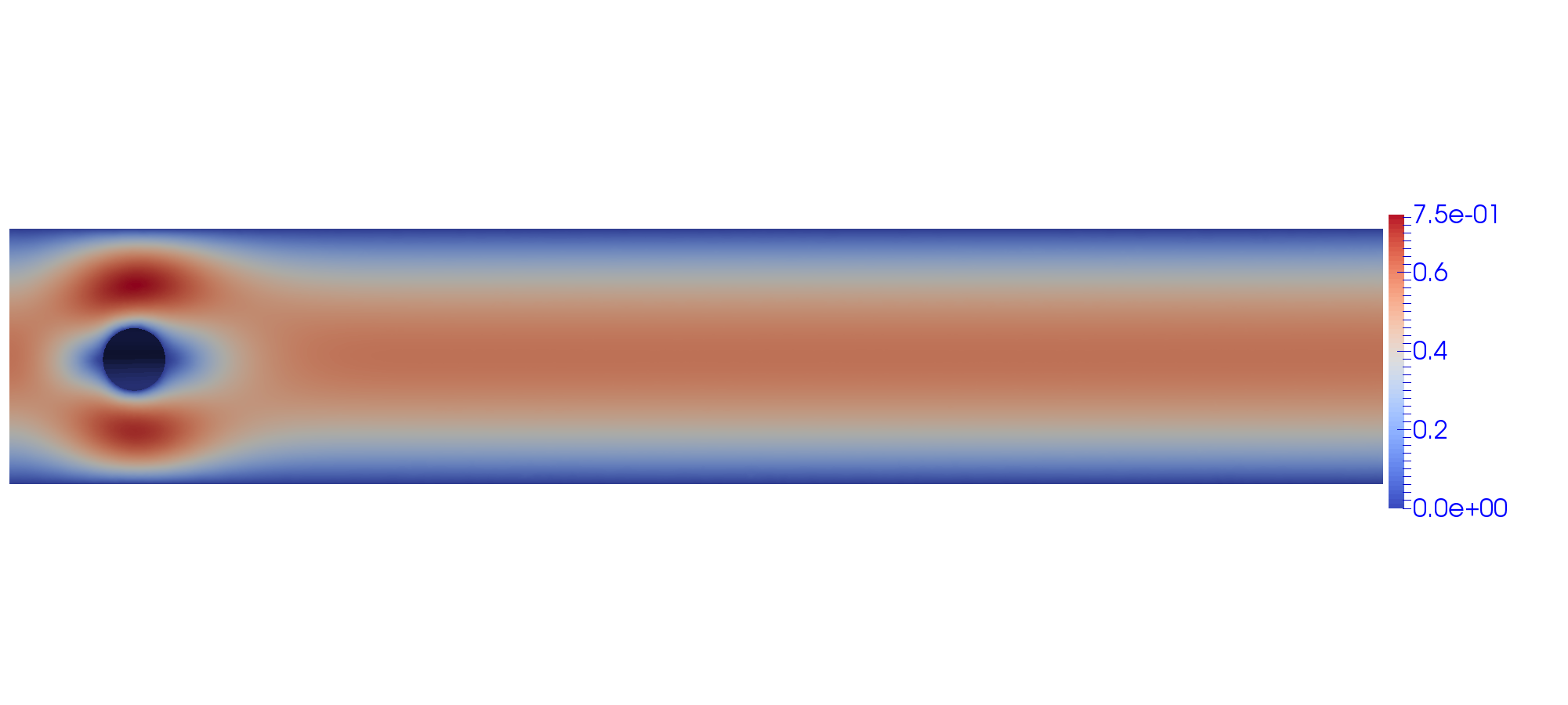}
        \put(40,35){$\v_{ROM}$}
      \end{overpic}\\
       \begin{overpic}[width=0.45\textwidth]{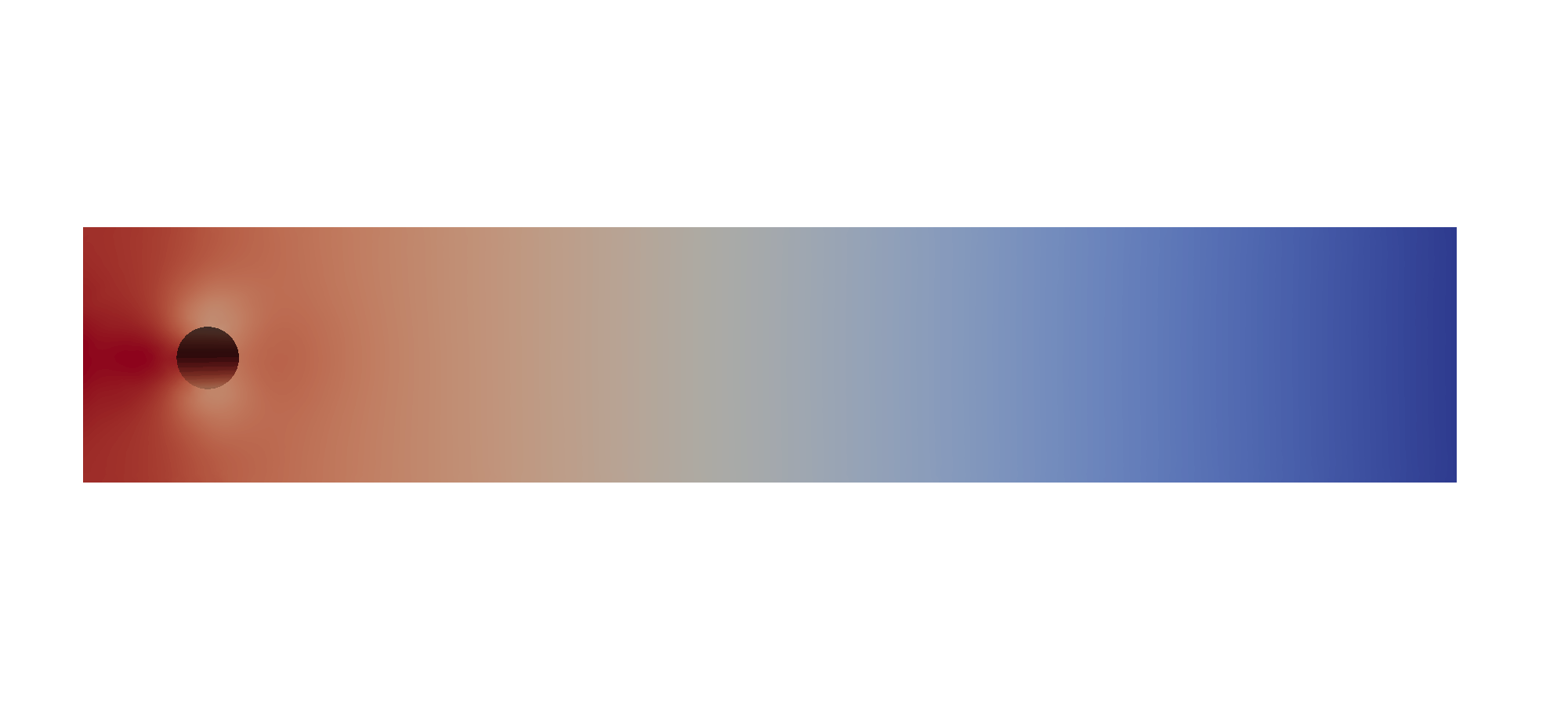}
        \put(40,35){$q_{FOM}$}
      \end{overpic}
       \begin{overpic}[width=0.45\textwidth]{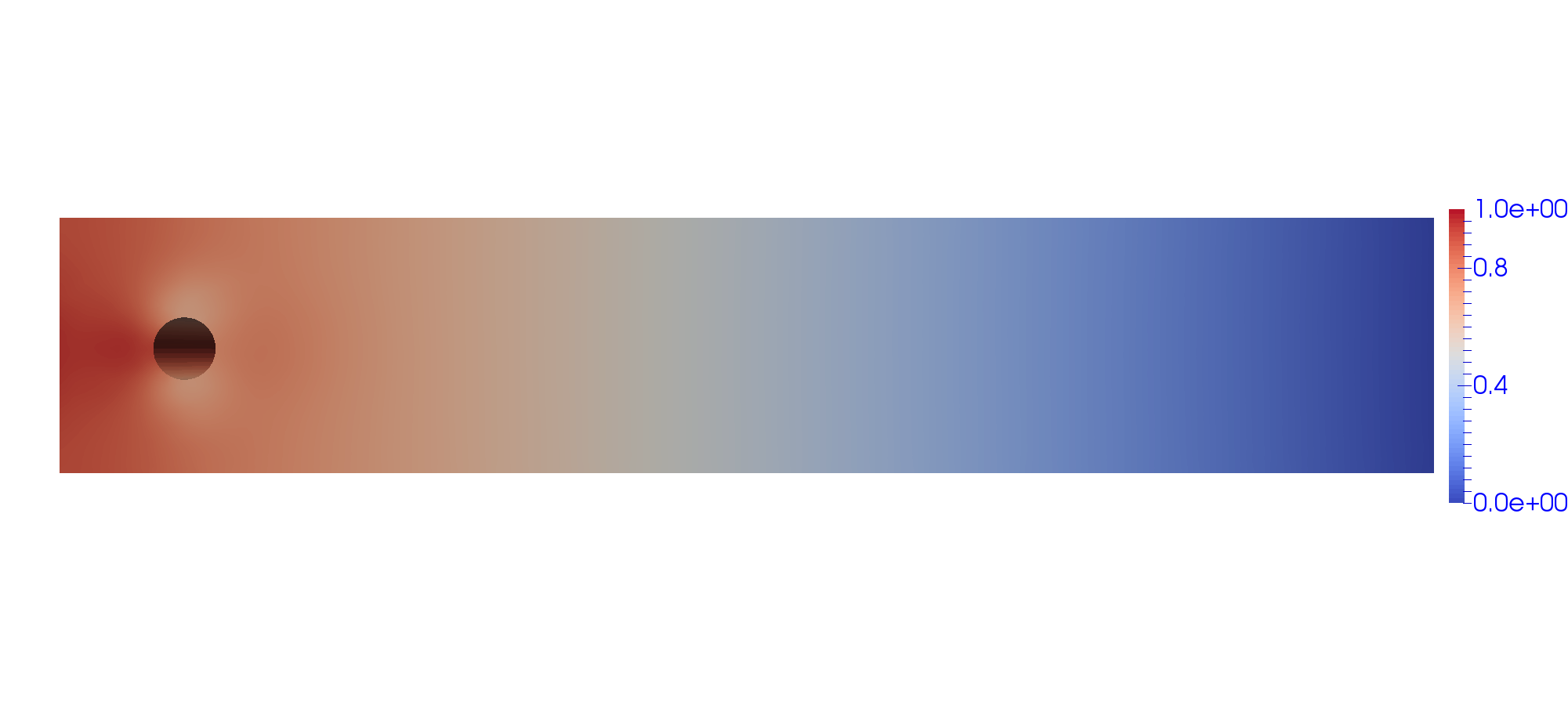}
        \put(40,35){$q_{ROM}$}
      \end{overpic}\\
       \begin{overpic}[width=0.45\textwidth]{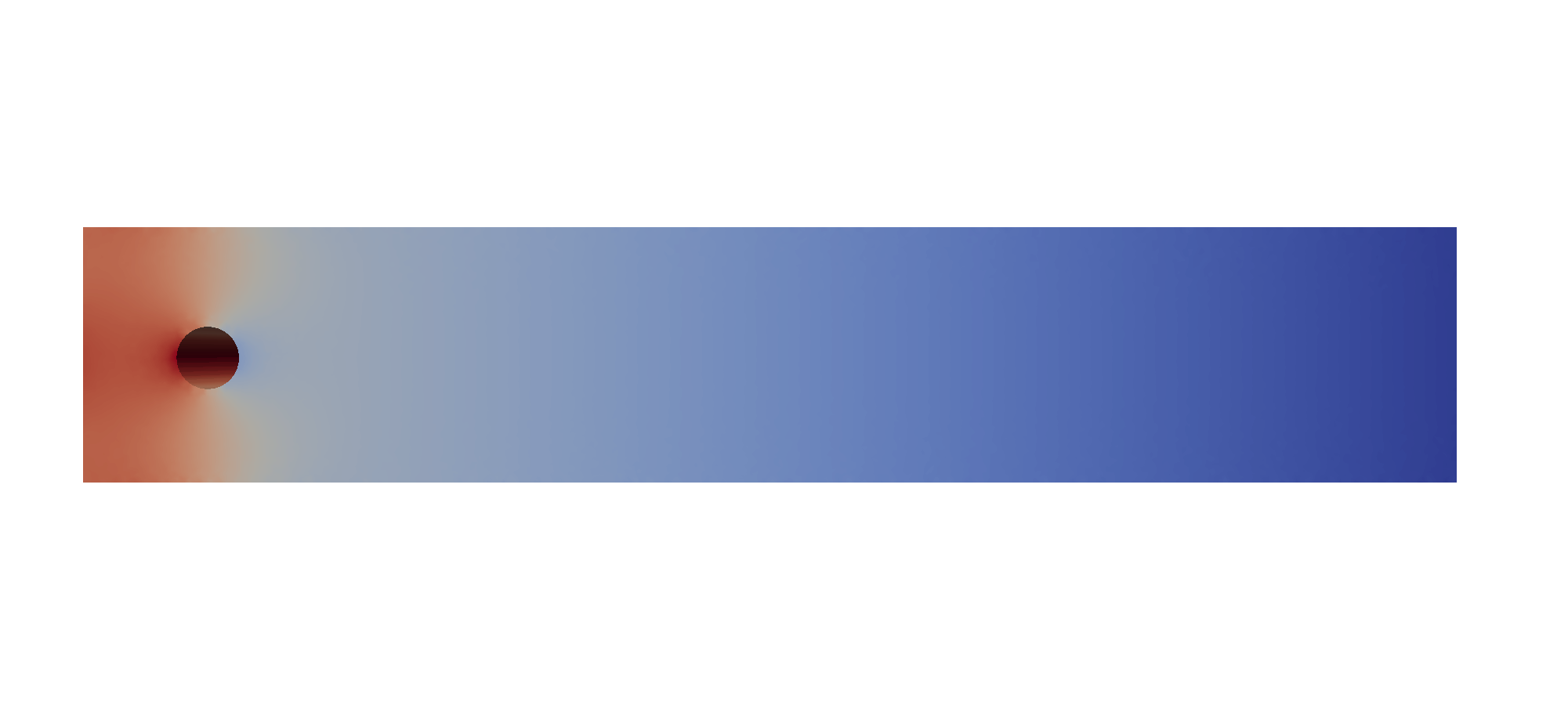}
        \put(40,35){$\overline{q}_{FOM}$}
      \end{overpic}
       \begin{overpic}[width=0.45\textwidth]{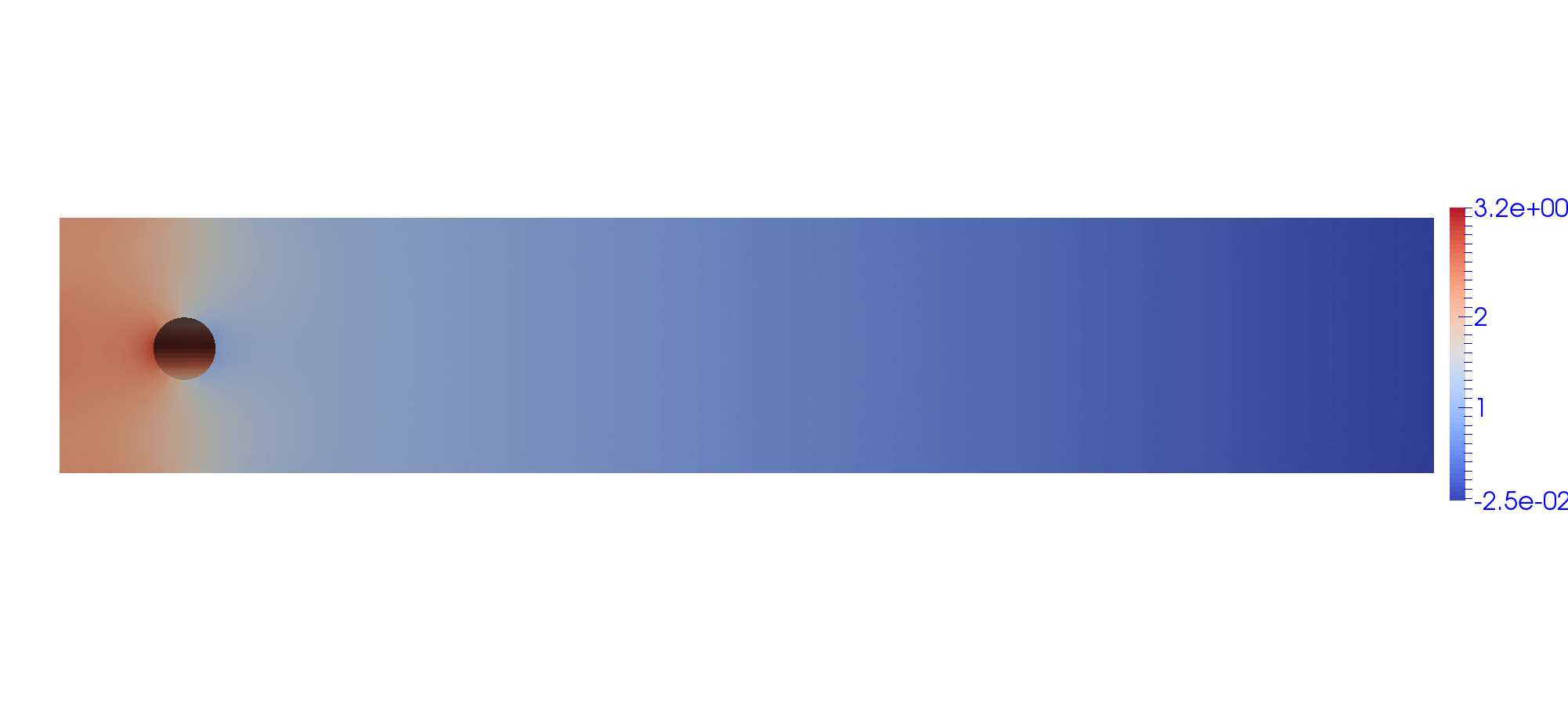}
        \put(40,35){$\overline{q}_{ROM}$}
      \end{overpic}
\caption{2D flow past a cylinder: comparison between FOM and ROM $\u$ (first row), $\v$ (second row), $q$ (third row),
and $\overline{q}$ (fourth row) at time $t = 1$. We consider 2 modes for $\v$, $\u$ and $q$, and 1 mode for $\bar{q}$. 
}
\label{fig:comp_FOM_ROM_t1}
\end{figure}

\begin{figure}
\centering
       \begin{overpic}[width=0.45\textwidth]{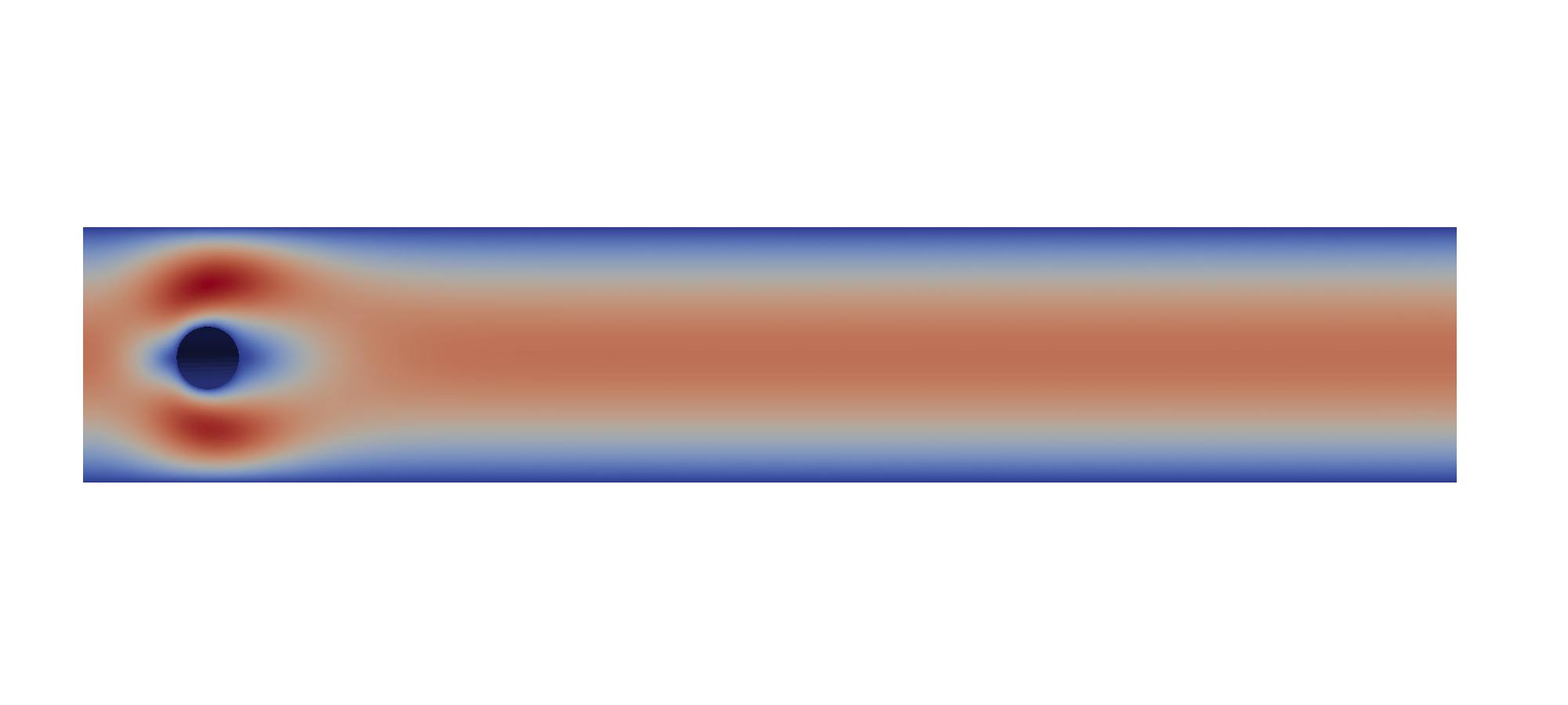}
        \put(40,35){$\u_{FOM}$}
      \end{overpic}
       \begin{overpic}[width=0.45\textwidth]{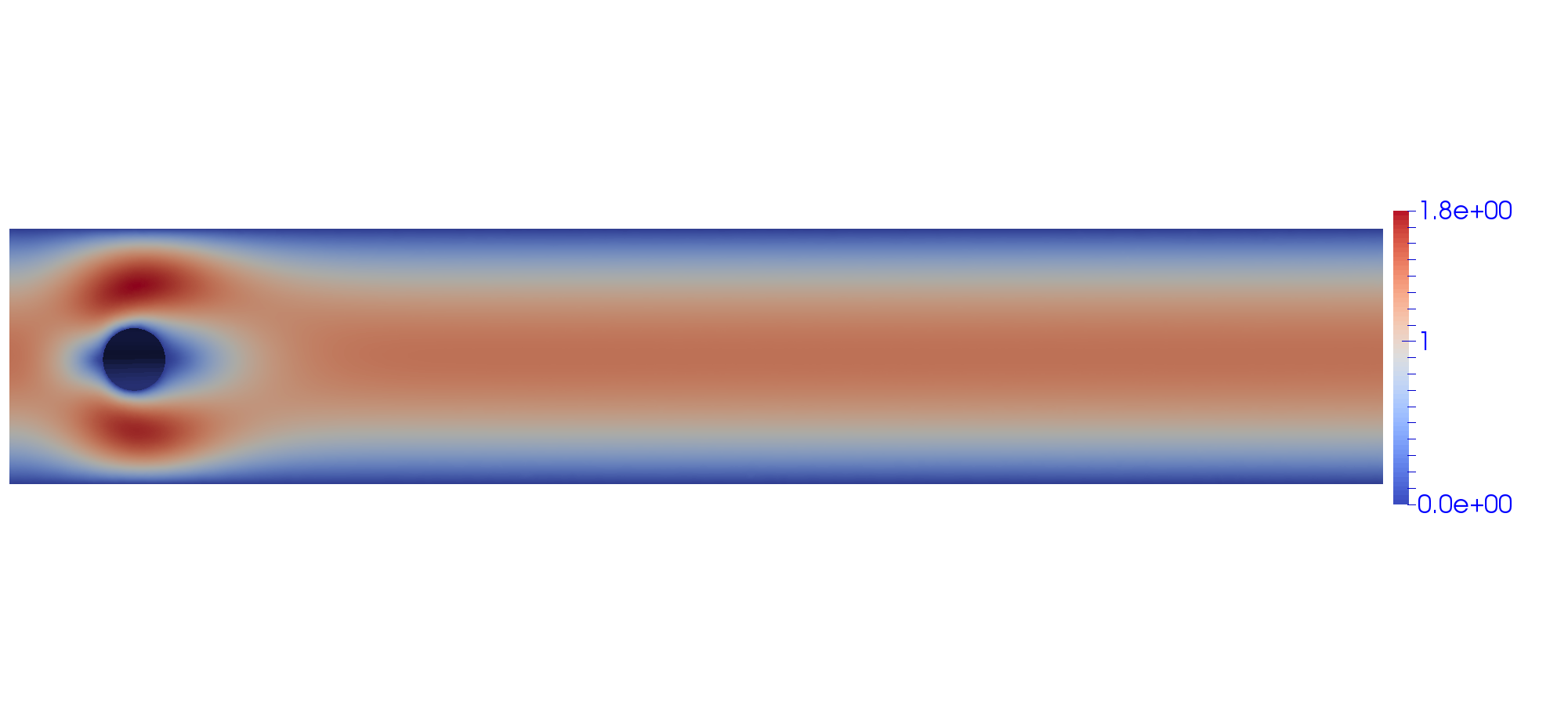}
        \put(40,35){$\u_{ROM}$}
      \end{overpic} \\
       \begin{overpic}[width=0.45\textwidth]{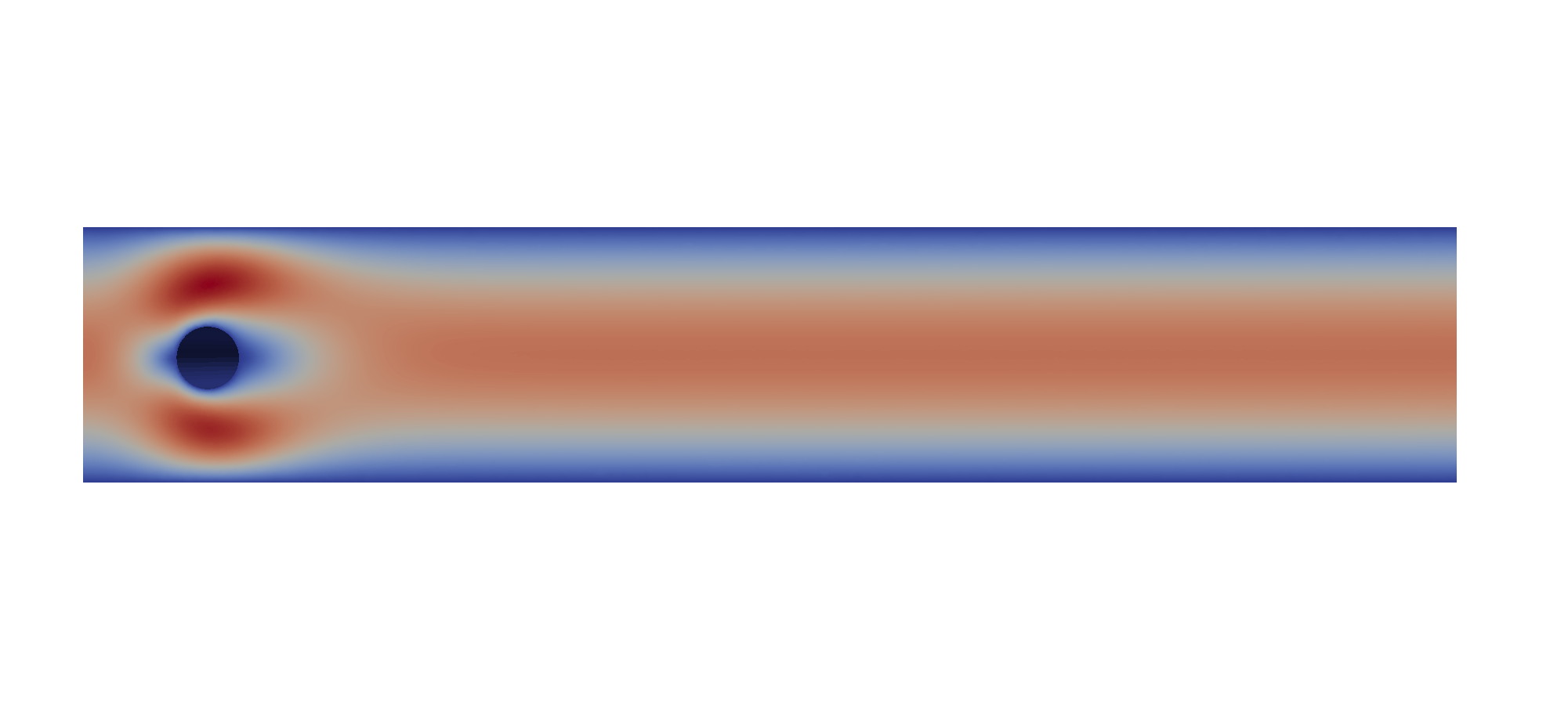}
        \put(40,35){$\v_{FOM}$}
      \end{overpic}
       \begin{overpic}[width=0.45\textwidth]{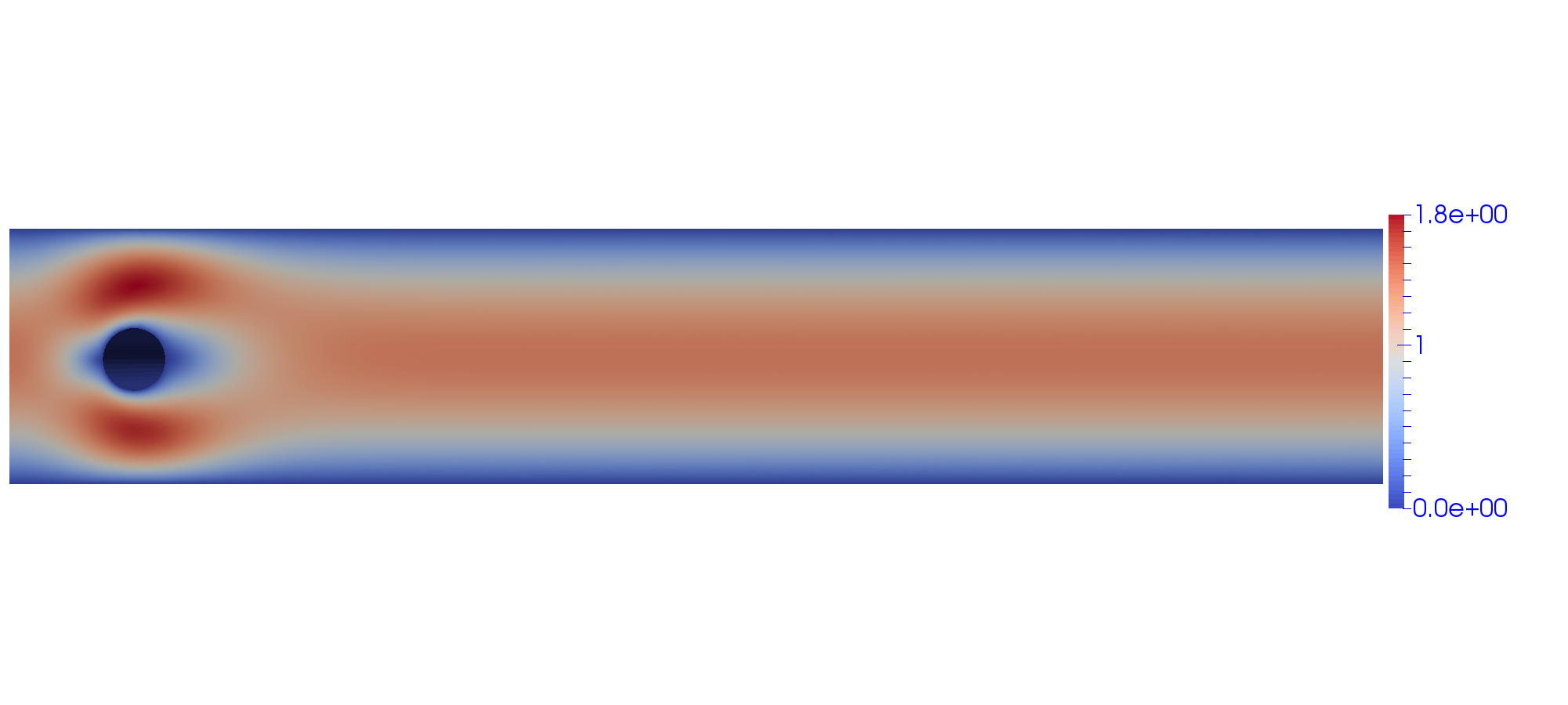}
        \put(40,35){$\v_{ROM}$}
      \end{overpic}\\
       \begin{overpic}[width=0.45\textwidth]{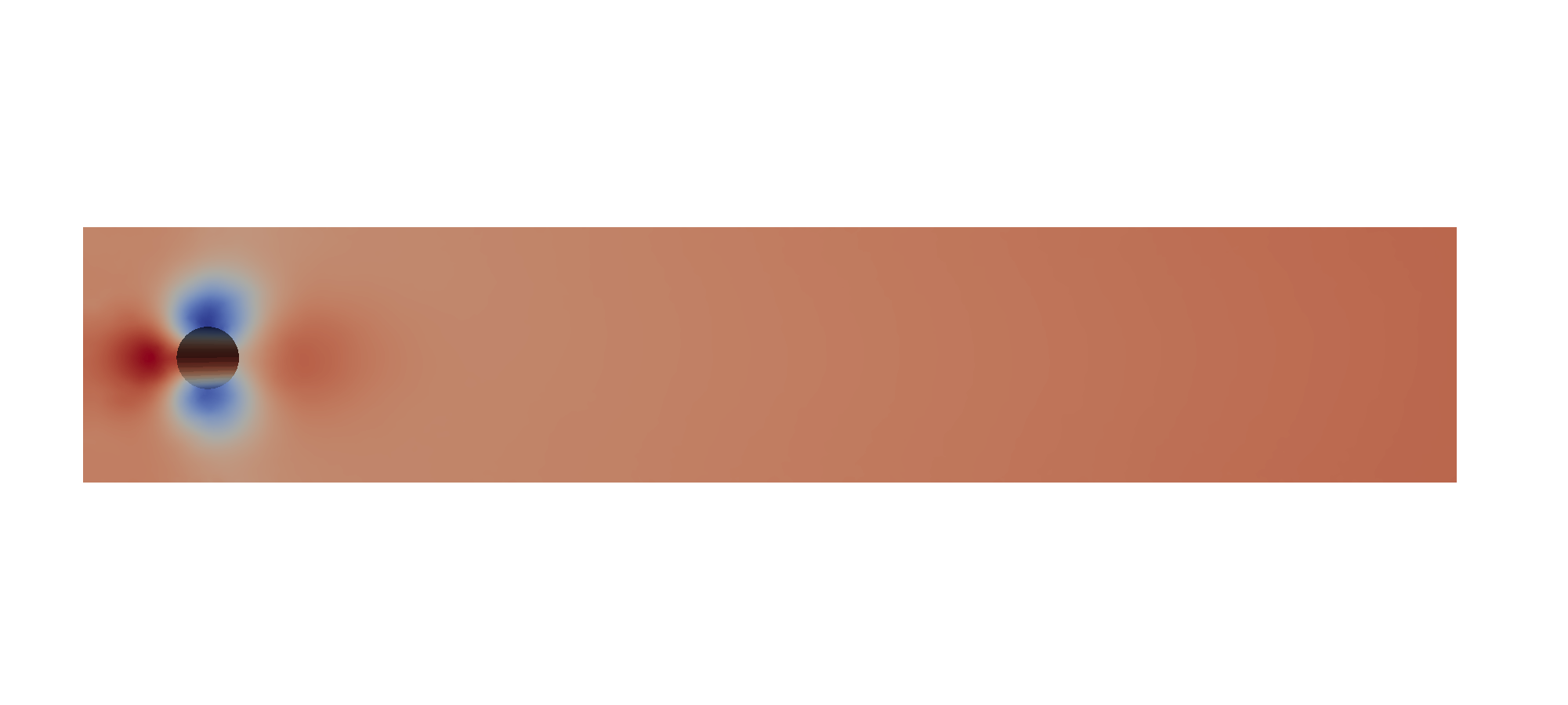}
        \put(40,35){$q_{FOM}$}
      \end{overpic}
       \begin{overpic}[width=0.45\textwidth]{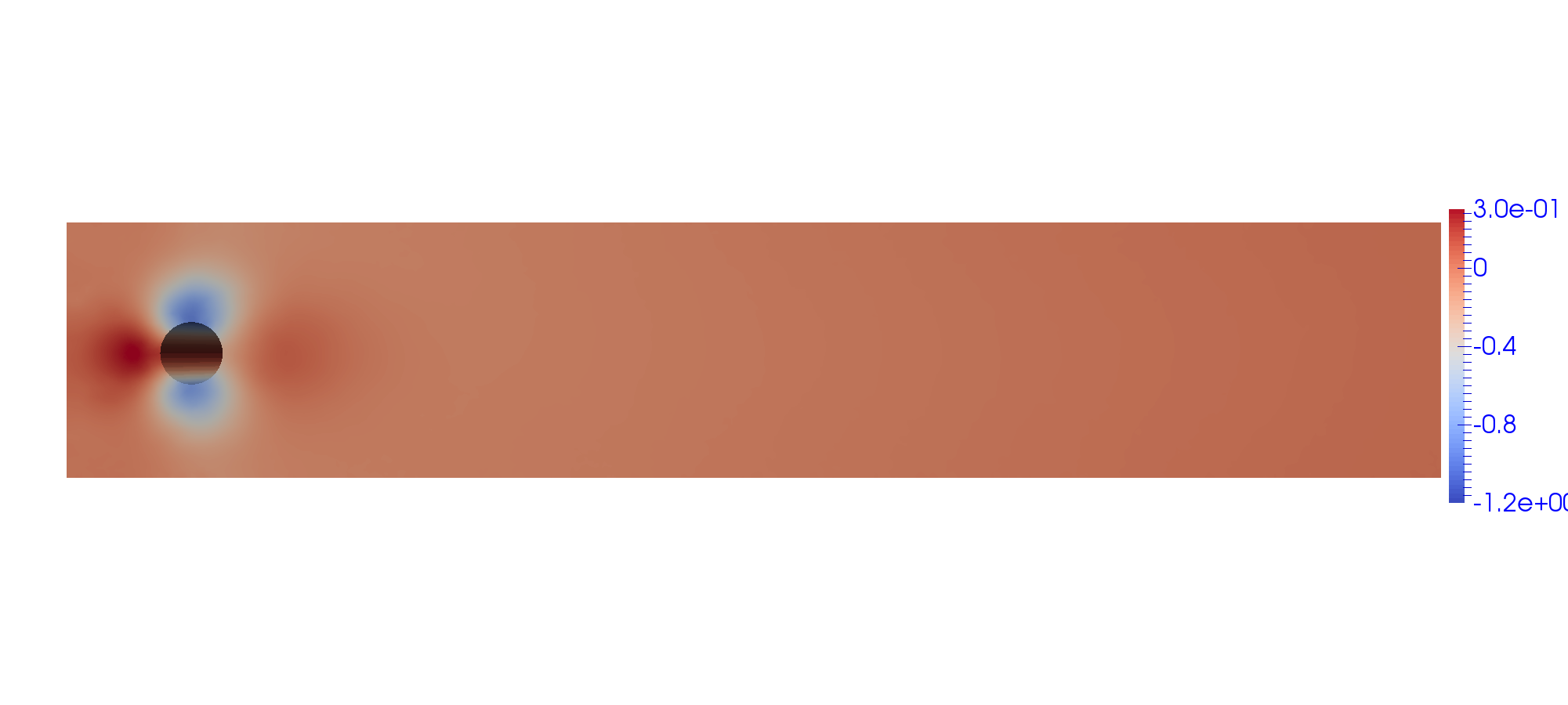}
        \put(40,35){$q_{ROM}$}
      \end{overpic}\\
       \begin{overpic}[width=0.45\textwidth]{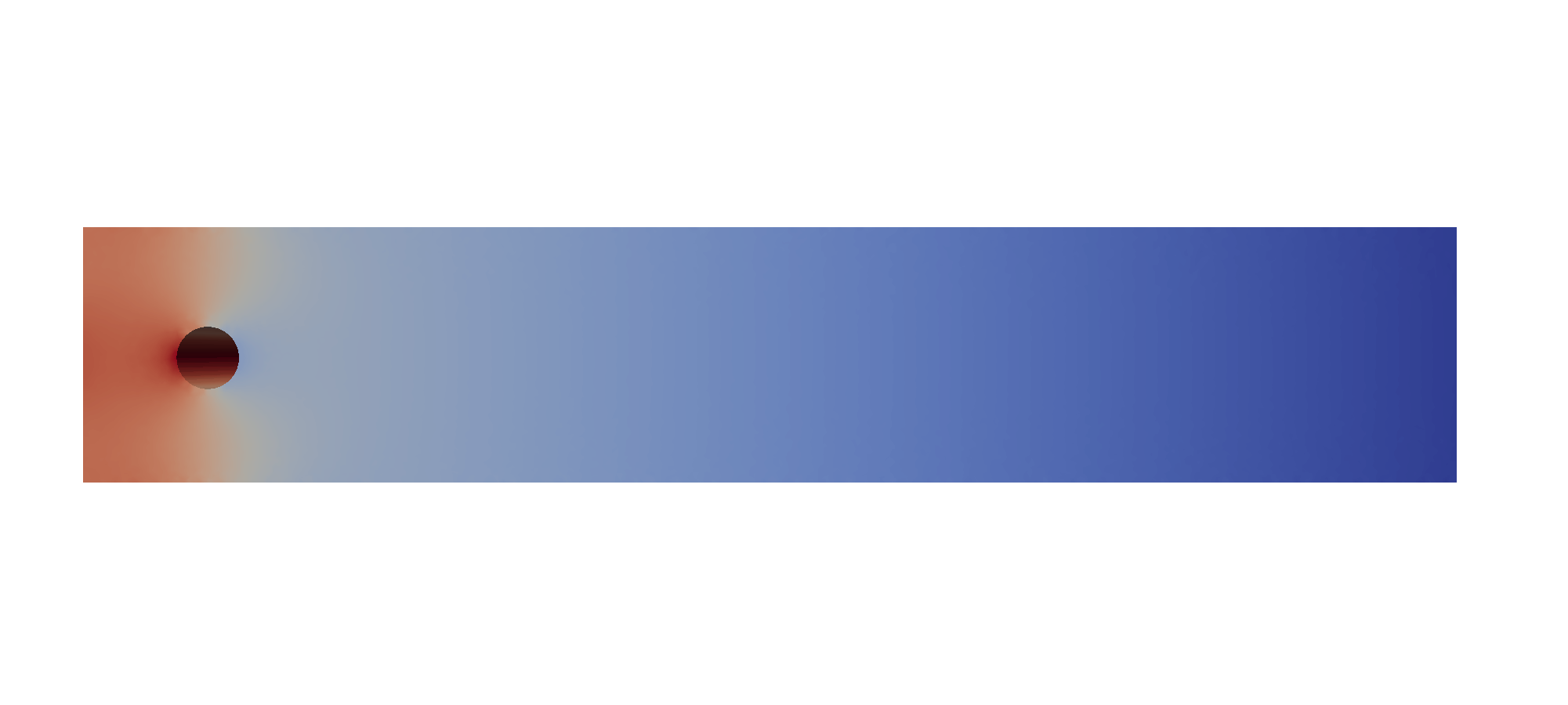}
        \put(40,35){$\overline{q}_{FOM}$}
      \end{overpic}
       \begin{overpic}[width=0.45\textwidth]{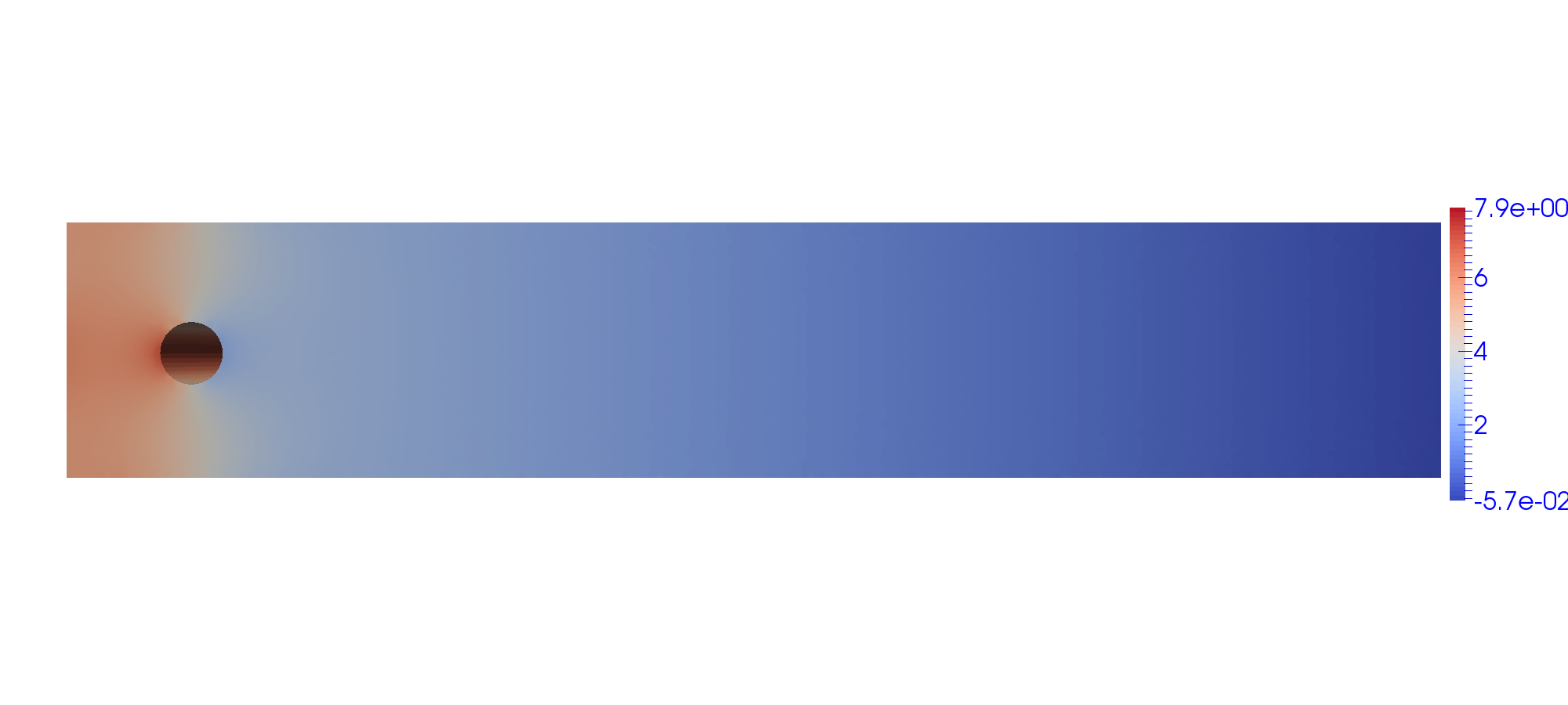}
        \put(40,35){$\overline{q}_{ROM}$}
      \end{overpic}
\caption{2D flow past a cylinder: comparison between FOM and ROM $\u$ (first row), $\v$ (second row), $q$ (third row),
and $\overline{q}$ (fourth row) at time $t = 5$. We consider 2 modes for $\v$, $\u$ and $q$, and 1 mode for $\bar{q}$.
}
\label{fig:comp_FOM_ROM_t5}
\end{figure}

\begin{figure}
\centering
       \begin{overpic}[width=0.45\textwidth]{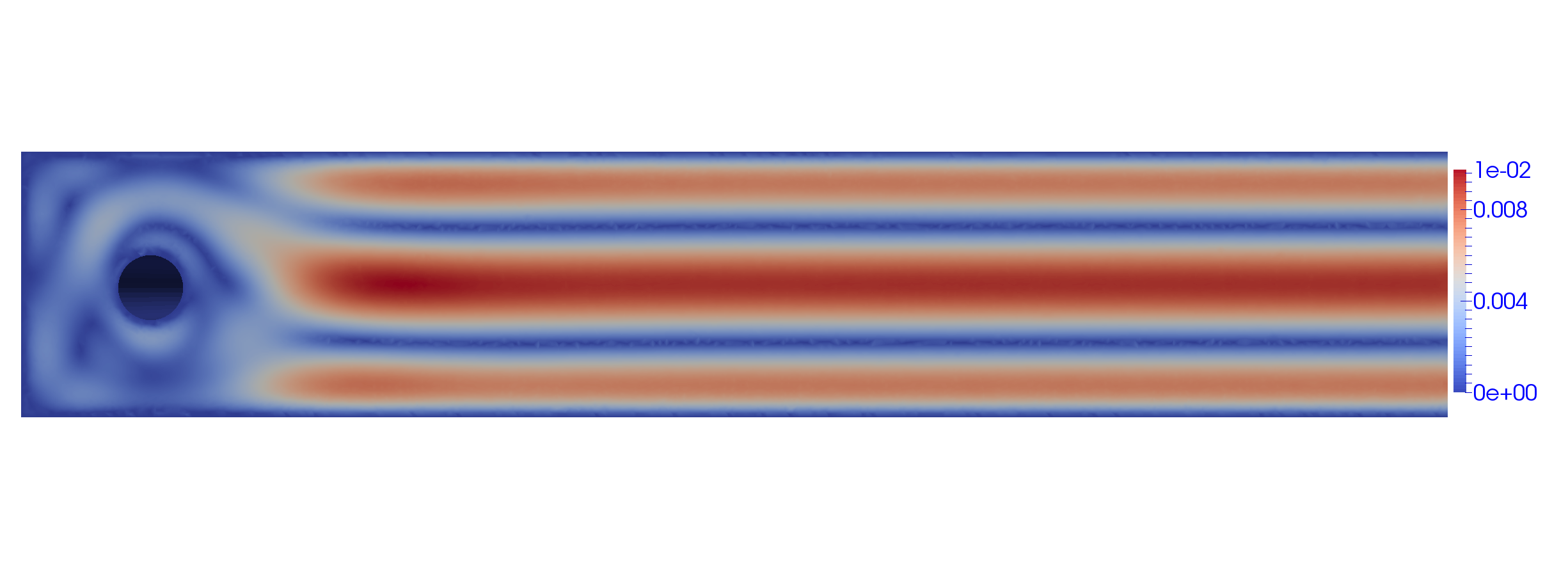}
        \put(20,30){$\v_{FOM} - \v_{ROM}$, $t = 1$}
      \end{overpic}
 \begin{overpic}[width=0.45\textwidth]{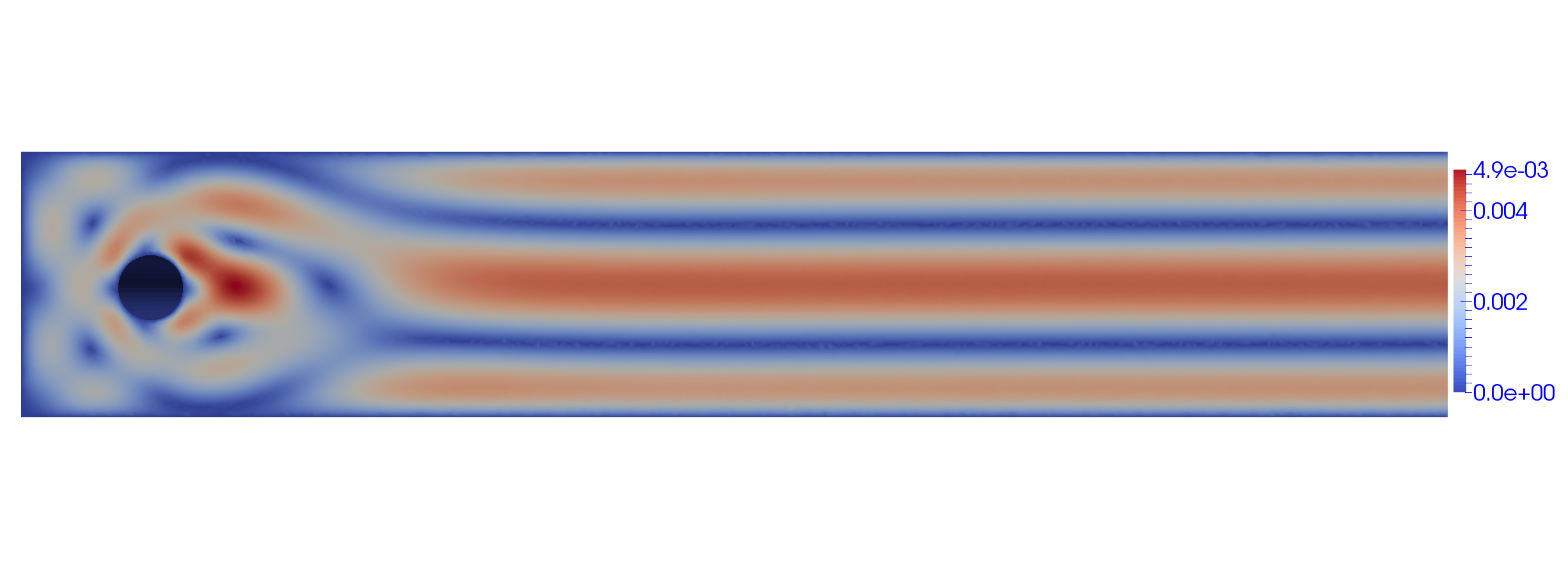}
        \put(20,30){$\v_{FOM} - \v_{ROM}$, $t = 5$}
      \end{overpic} \\
       \begin{overpic}[width=0.45\textwidth]{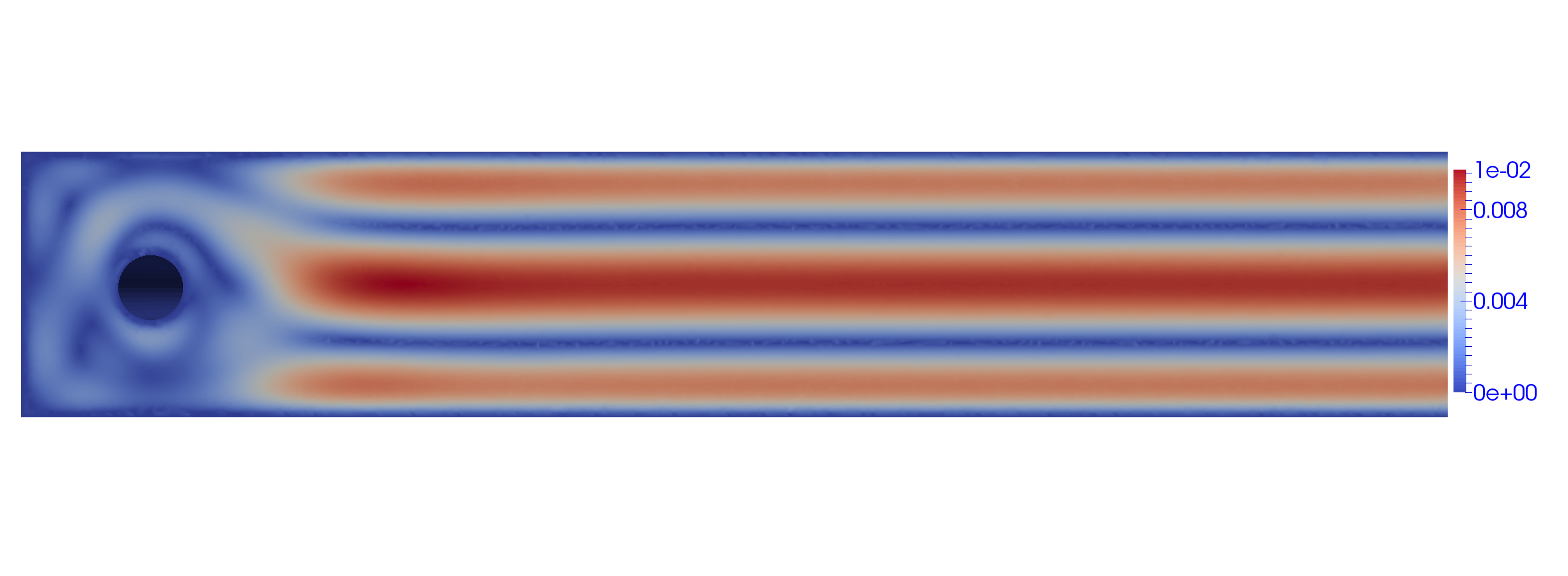}
        \put(20,30){$\u_{FOM} - \u_{ROM}$, $t = 1$}
      \end{overpic}
 \begin{overpic}[width=0.45\textwidth]{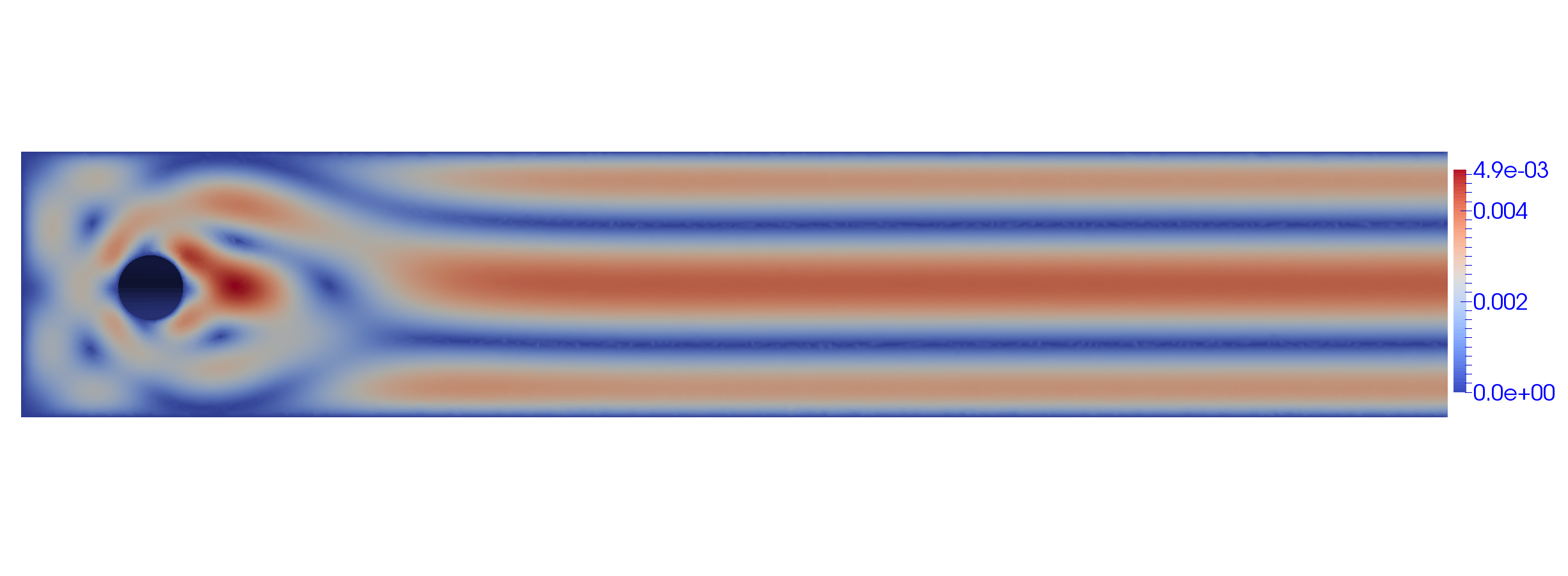}
        \put(20,30){$\u_{FOM} - \u_{ROM}$, $t = 5$}
      \end{overpic}\\
       \begin{overpic}[width=0.45\textwidth]{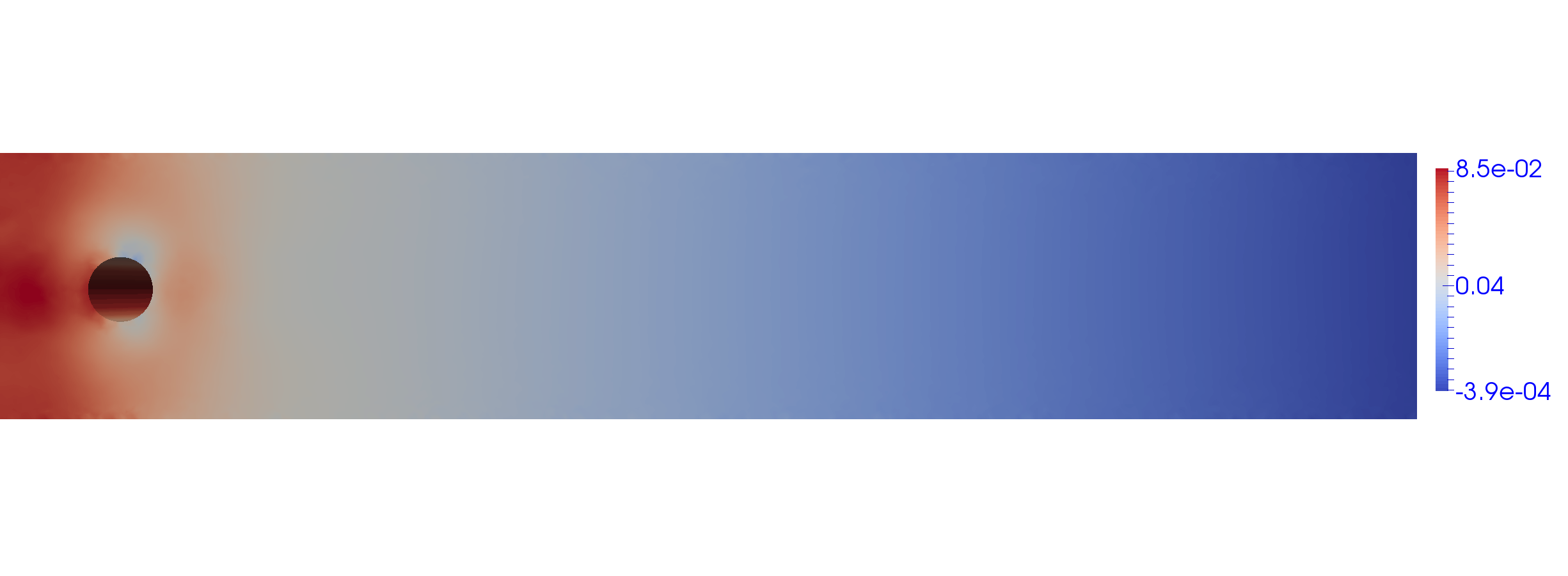}
        \put(20,30){$q_{FOM} - q_{ROM}$, $t = 1$}
      \end{overpic}
      \begin{overpic}[width=0.45\textwidth]{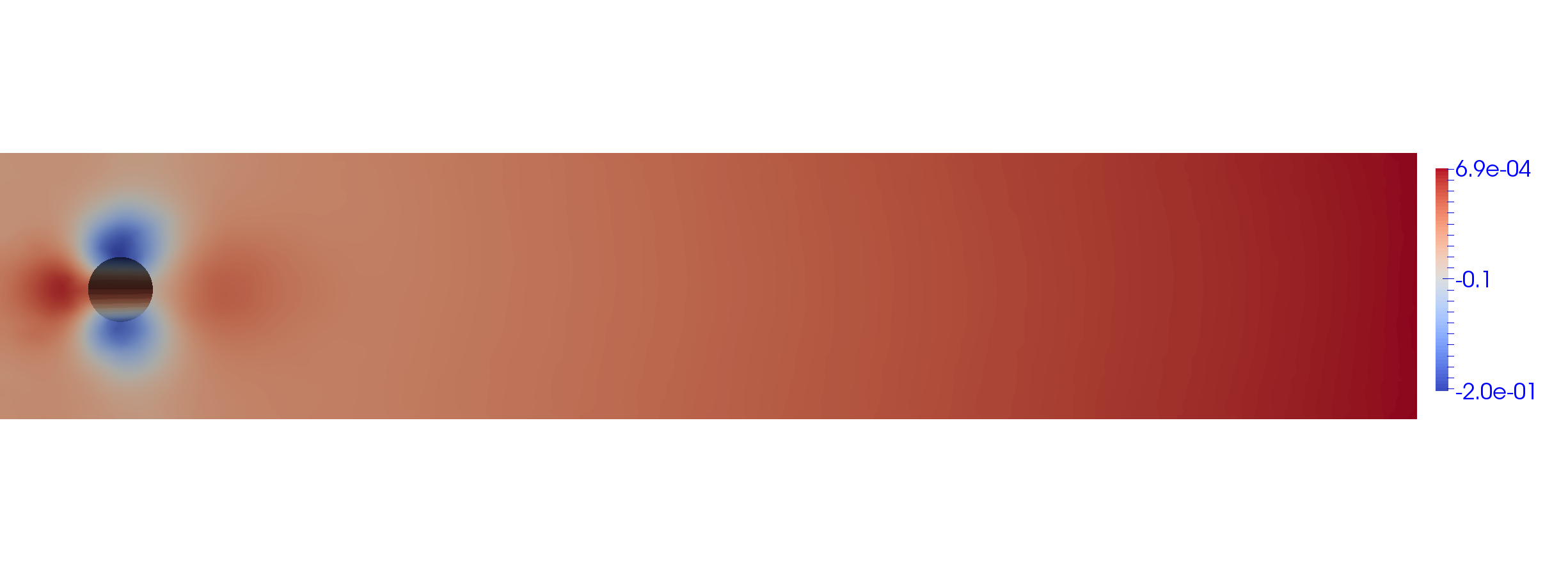}
        \put(20,30){$q_{FOM} - q_{ROM}$, $t = 5$}
      \end{overpic}\\
       \begin{overpic}[width=0.45\textwidth]{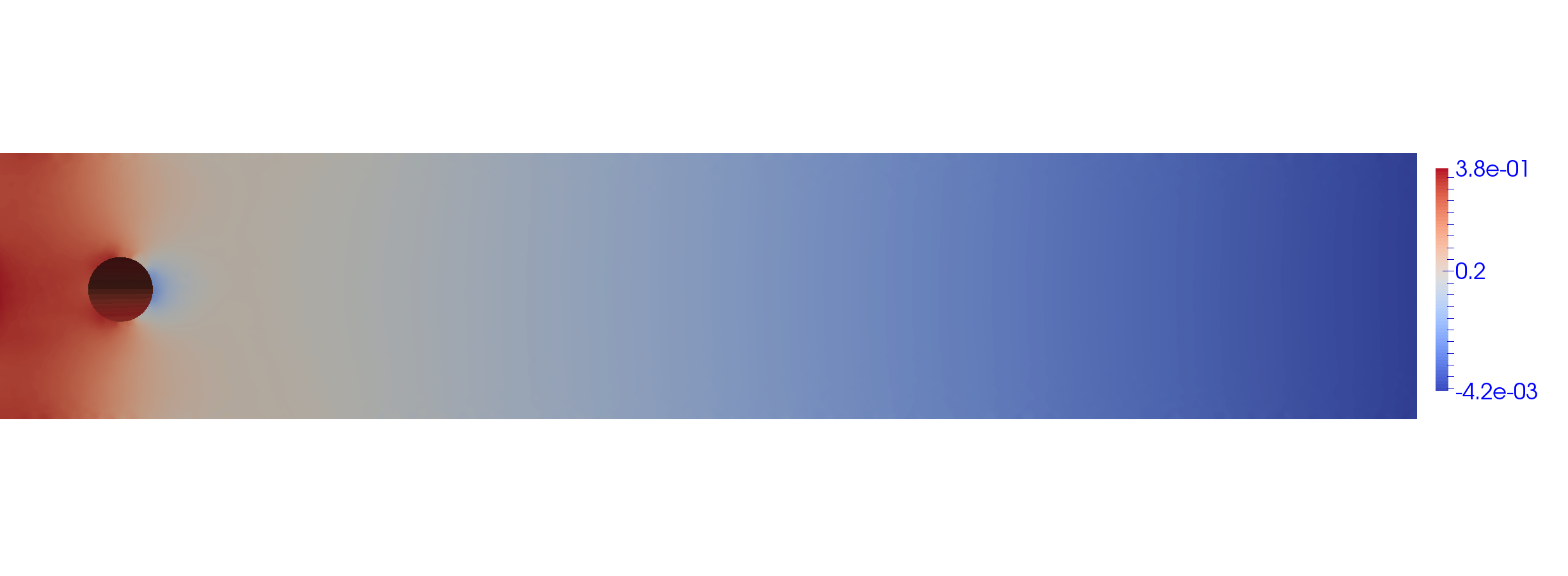}
        \put(20,30){$\overline{q}_{FOM} - \overline{q}_{ROM}$, $t = 1$}
      \end{overpic}
 \begin{overpic}[width=0.45\textwidth]{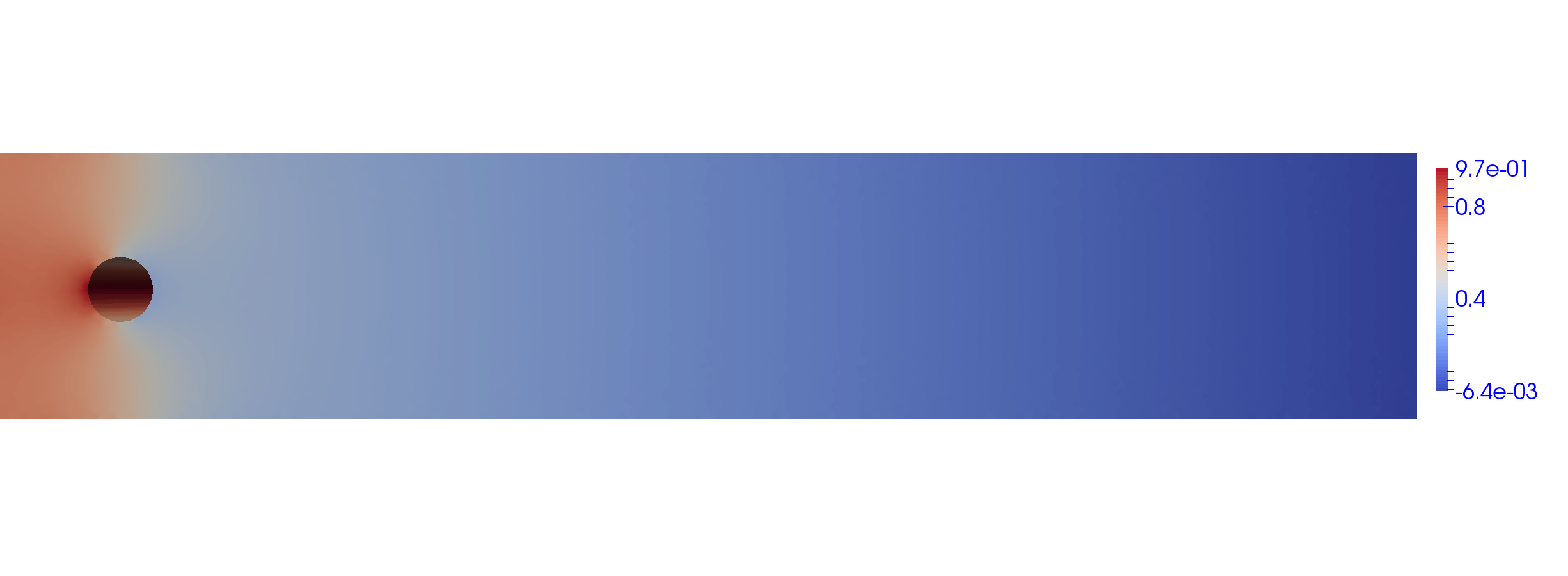}
        \put(20,30){$\overline{q}_{FOM} - \overline{q}_{ROM}$, $t = 5$}
      \end{overpic}
\caption{2D flow past a cylinder: difference between FOM and ROM $\v$ (first row), $\u$ (second row), $q$ (third row),
and $\overline{q}$ (fourth row) at times $t = 1$ (left) and $t = 5$ (right). We consider 2 modes for $\v$, $\u$ and $q$, and 1 mode for $\bar{q}$.
}
\label{fig:comp_t}
\end{figure}

Figs.~\ref{fig:comp_FOM_ROM_t1} and~\ref{fig:comp_FOM_ROM_t5} display the comparison between the computed FOM and ROM fields at two different times, $t = 1$ and $t = 5$. As one can see from the figures, the ROM provides a good global reconstruction 
of both velocities and both pressures. 
Next, we make the comparison more quantitative. Fig.~\ref{fig:comp_t} displays the difference between the computed FOM and ROM 
fields for $t = 1$ and $t = 5$.
The maximum absolute error between the FOM and the ROM for the velocity fields 
is of the order of $10^{-2}$ at $t = 1$, while it reaches lower values (of order $10^{-3}$) at 
$t = 5$. The maximum absolute error for $q$ is of the order of  $10^{-1}$ at $t = 1,5$.
Again, these values are in perfect agreement with those reported 
for the lid driven cavity and Y-junction flow in \cite{Star2019} (no absolute error is reported in
\cite{Stabile2017, Stabile2018}) and they indicate
that our ROM is able to reproduce the main flow features at different times. 
Finally, we note that the maximum absolute error for $\qbar$ is of the order of 
$10^{-1}$ at $t = 1$ and $1$ at $t=5$. However, we note that at $t=5$ 
the order of magnitude of $\qbar$ is $10$, as shown in Fig.~\ref{fig:comp_FOM_ROM_t5}.


\begin{figure}
\centering
 \begin{overpic}[width=0.45\textwidth]{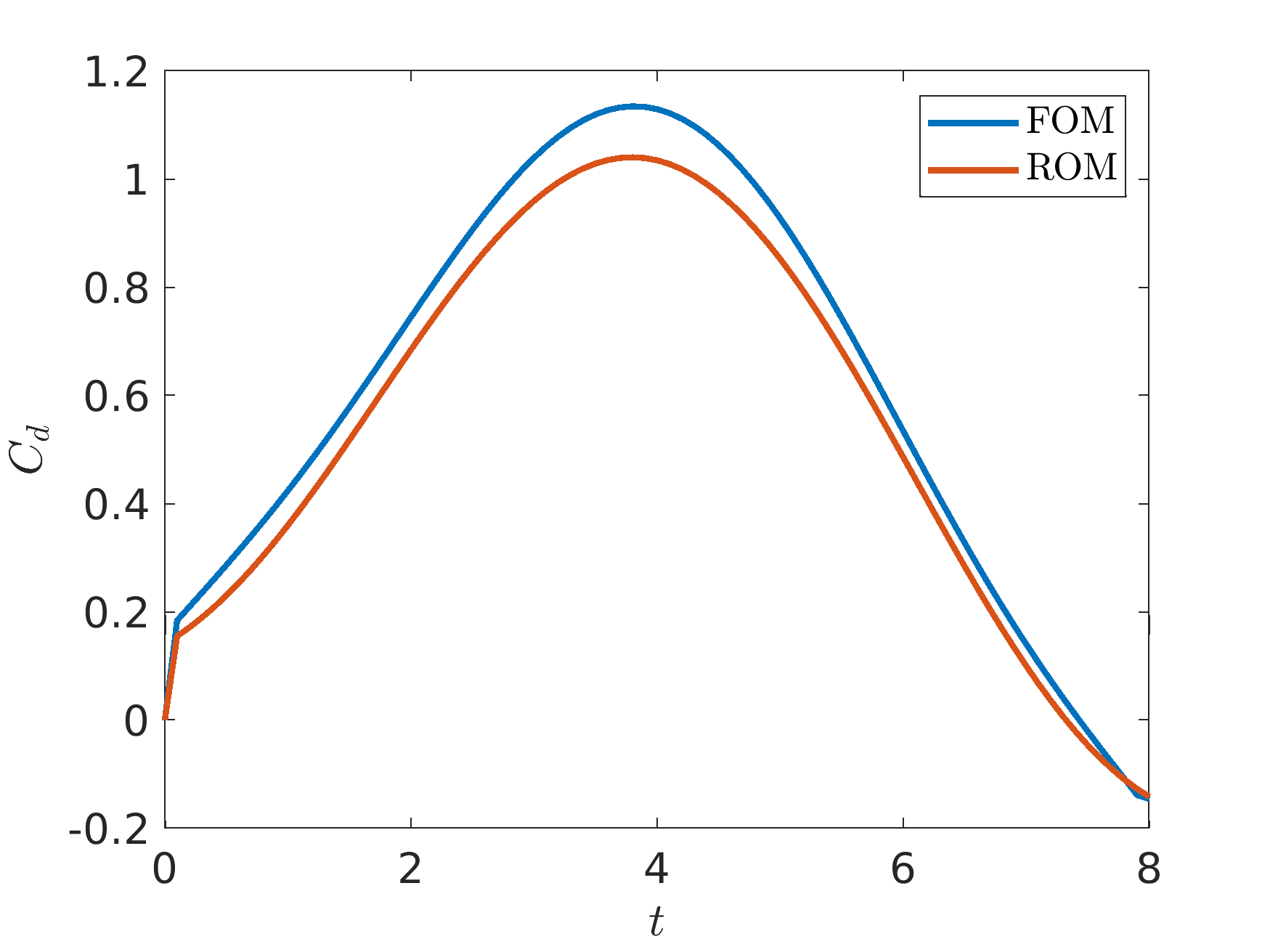}
      \end{overpic}
       \begin{overpic}[width=0.45\textwidth]{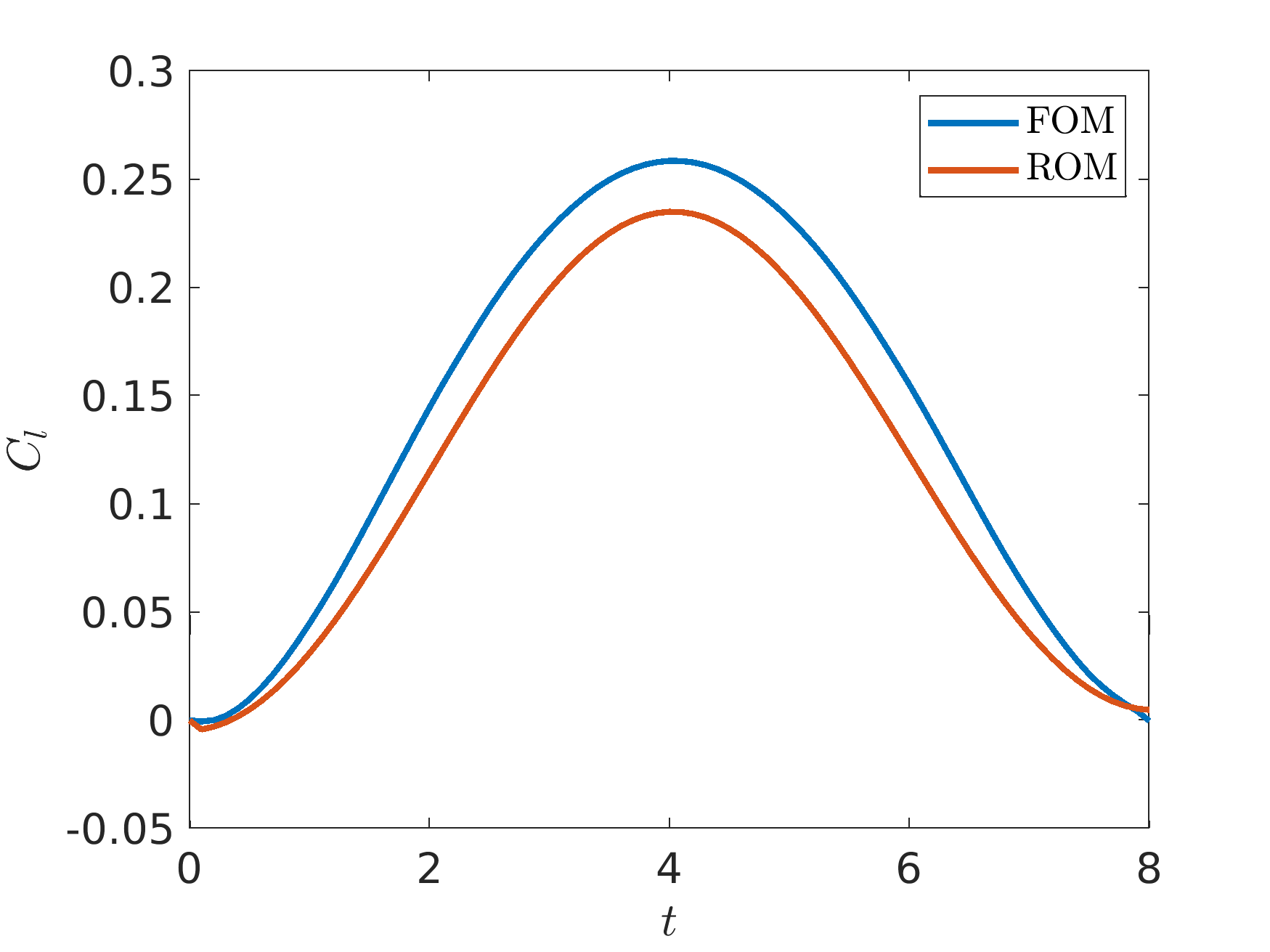}
      \end{overpic}\\
\caption{2D flow past a cylinder: aerodynamic coefficients $C_d$ (left) and $C_l$ (right) computed by FOM and ROM. We consider 2 modes for $\v$, $\u$ and $q$, and 1 mode for $\bar{q}$.}
\label{fig:coeff_t}
\end{figure}

To further quantitate the accuracy of the ROM with respect to the FOM, we consider
the quantities of interest for this benchmark, i.e.~the drag and lift coefficients \cite{John2004,turek1996}:
\begin{align}\label{eq:cd_cl}
c_d(t) = \dfrac{2}{\rho L_{r} {U}^2_{r}} \int_S \left(\left(2 \mu \nabla \u - q\boldsymbol{I}\right)
\cdot \boldsymbol{n}\right) \cdot \boldsymbol{t}~dS, \quad
c_l(t) = \dfrac{2}{\rho L_{r} {U}^2_{r}} \int_S \left(\left(2 \mu \nabla \u - q\boldsymbol{I}\right)
\cdot \boldsymbol{n}\right) \cdot \boldsymbol{n}~dS,
\end{align}
where $U_{r}= 1$ is the maximum velocity at the inlet/outlet, $L_r = 0.1$ is the cylinder diameter, 
$S$ is the cylinder surface, and $\boldsymbol{t}$ and $\boldsymbol{n}$ are the tangential and normal unit vectors
to the cylinder, respectively. The FOM/ROM comparison for the coefficients in \eqref{eq:cd_cl} over time
is reported in Fig.~\ref{fig:coeff_t}.
We observe that the amplitude of the force coefficients are slightly underestimated for all the time instants by ROM.
The ROM reconstruction of the lift coefficient appears to be more critical, especially around the center of the time interval. 
This could be due to the fact that larger errors for pressure $q$ are localized close to the cylinder, 
as one can see in Fig.~\ref{fig:comp_t} (third row). 

Table \ref{tab:coeffs_t} compares the maximum lift and drag coefficients 
and times at which the maxima occur. In addition, we obtain the following errors:
\begin{align}
E_{c_d} = \dfrac{c_{d,max}^{FOM} - c_{d,max}^{ROM}}{c_{d,max}^{FOM}} = 0.089,~
E_{t_{c_d}} = \dfrac{t_{c_{d,max}}^{FOM} - t_{c_{d,max}}^{ROM}}{t_{c_{d,max}}^{FOM}} = 0, \label{eq:E_c_l} \\
E_{c_l} = \dfrac{c_{l,max}^{FOM} - c_{l,max}^{ROM}}{c_{l,max}^{FOM}} = 0.083,~
E_{t_{c_l}} = \dfrac{t_{c_{l,max}}^{FOM} - t_{c_{l,max}}^{ROM}}{t_{c_{l,max}}^{FOM}} = 0. \label{eq:E_c_d}
\end{align}
We see that our ROM approach is able to provide a perfect prediction of the time istants where maxima values of the aerodynamic coefficients occur. 
The errors related to the maximum values are both lower than 9\%. In \cite{Stabile2017}, the authors
use the NSE model with no filtering for the steady flow at $Re = 100$ around a cylinder 
and find errors smaller than 5\% for both coefficients. Additional differences, such 
as steady boundary conditions and homogeneous Neumann boundary condition for pressure
in \cite{Stabile2017}, make it harder to identify the reason why we obtain larger errors. 

\begin{table}
\centering
\begin{tabular}{ccccc}
\multicolumn{5}{c}{} \\
\cline{1-5}
 & $t(c_{l,max})$ & $c_{l,max}$ & $t(c_{d,max})$& $c_{d,max}$   \\
\hline
FOM & 4.1  & 0.258 & 3.9 & 1.135  \\
ROM & 4.1  & 0.235 & 3.9 & 1.041  \\
\hline
\end{tabular}
\caption{2D flow past a cylinder: maximum lift and drag coefficients computed by FOM and ROM and times at which the maxima occur.
We recall that the sampling frequency is 0.1. }
\label{tab:coeffs_t}
\end{table}
%

Regarding the computational cost, the CPU time of the FOM model is $3800$ s. 
The CPU time of the ROM is $30$ s. This corresponds to a speed-up of $\approx 127$.

\subsubsection{Parametrization with respect to the filtering radius $\alpha$}

Filtering radius $\alpha$ plays a crucial role in the success of filtering stabilization. So, after
having investigated the ability of our ROM approach to reconstruct the time evolution of velocity and pressure fields,
we consider $\alpha$ as a parameter. To train the ROM, we choose a uniform sample distribution
in the range $\alpha \in [0.0032, h_{min}]$, where $0.0032$ is $Re^{-3/4} L_r$.
We first consider 11 sampling points and then decrease to 6 sampling points. 
For each value of the filtering radius inside the training set, a simulation is run for the entire time interval 
of interest, i.e.~$(0, 8]$. 
Based on the results presented above, the snapshots are collected every 0.1 in this time window, 
for a total number of 880 (resp., 480) snapshots for the first (resp., second) sampling case. 

\begin{figure}
\centering
 \begin{overpic}[width=0.4\textwidth]{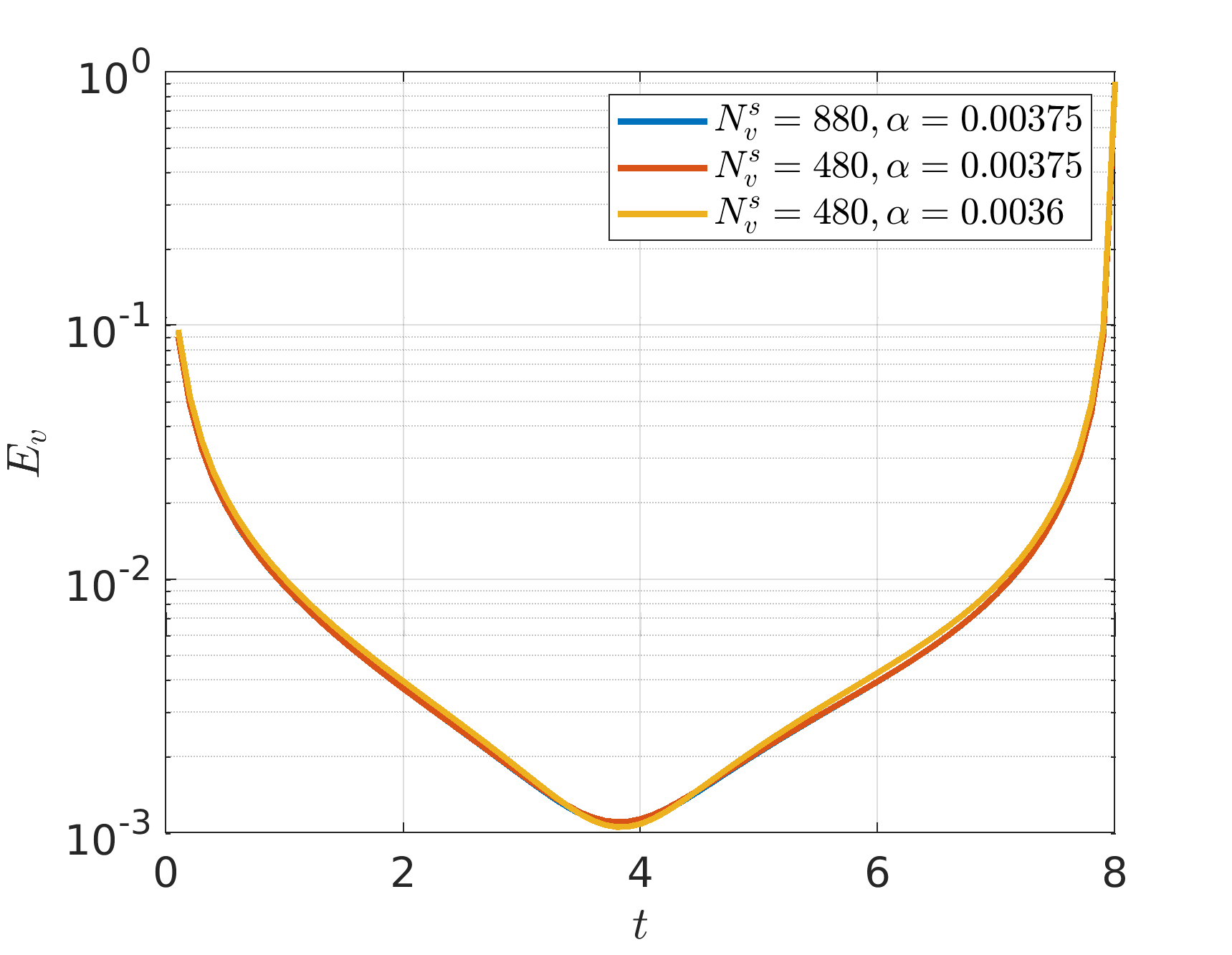}
      \end{overpic}
\begin{overpic}[width=0.45\textwidth]{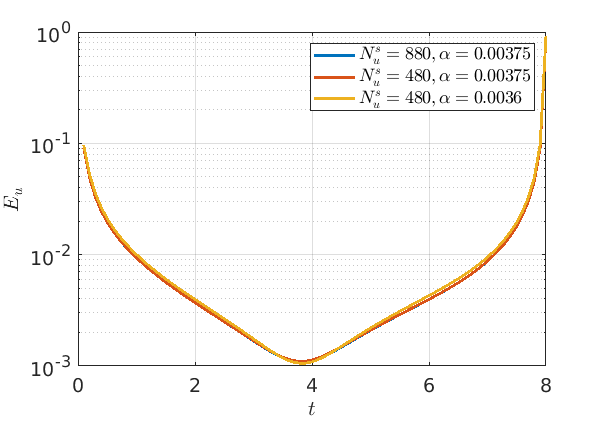}
      \end{overpic}\\
\begin{overpic}[width=0.45\textwidth]{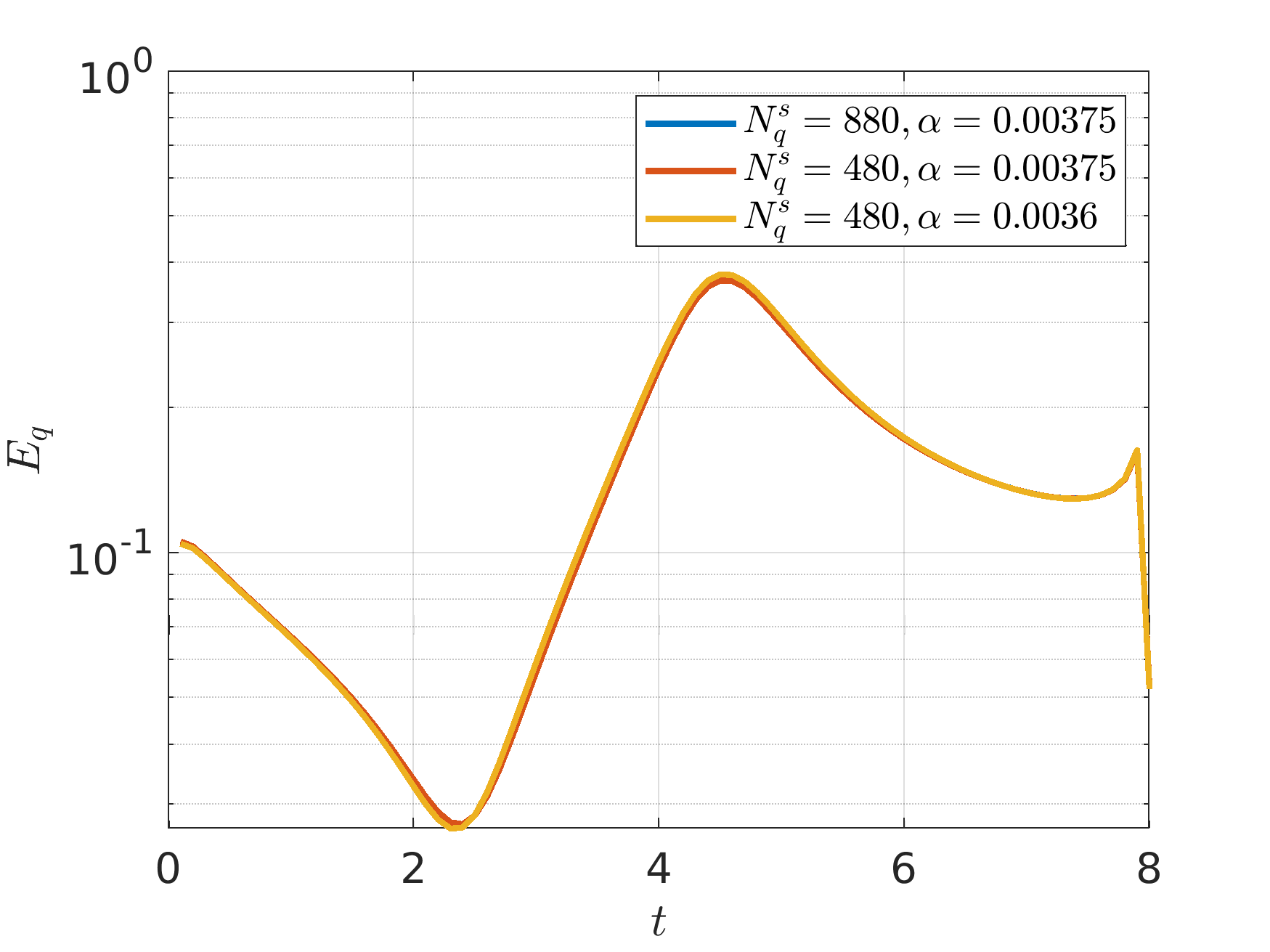}
      \end{overpic}
\begin{overpic}[width=0.45\textwidth]{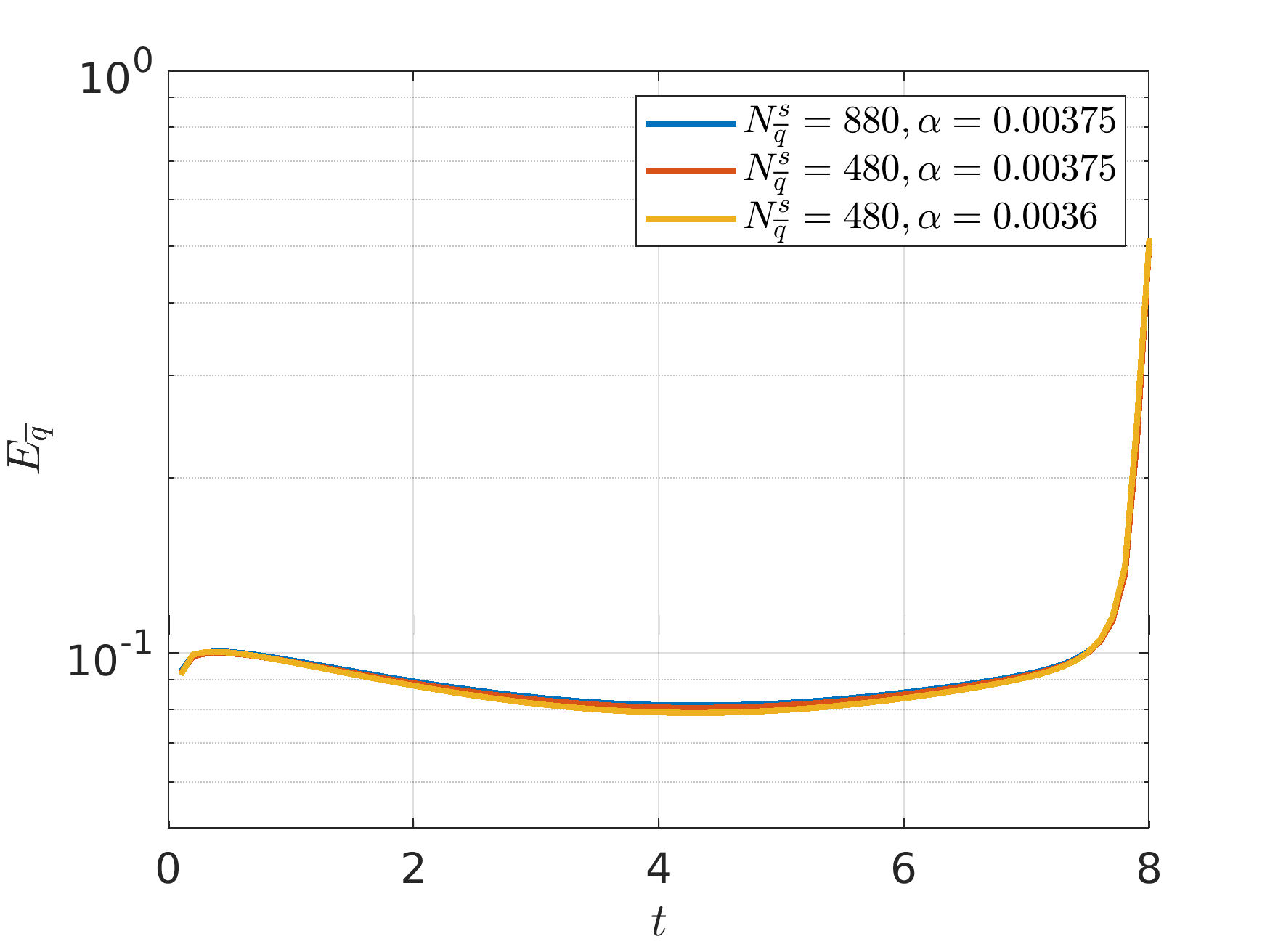}
      \end{overpic}\\
\begin{overpic}[width=0.45\textwidth]{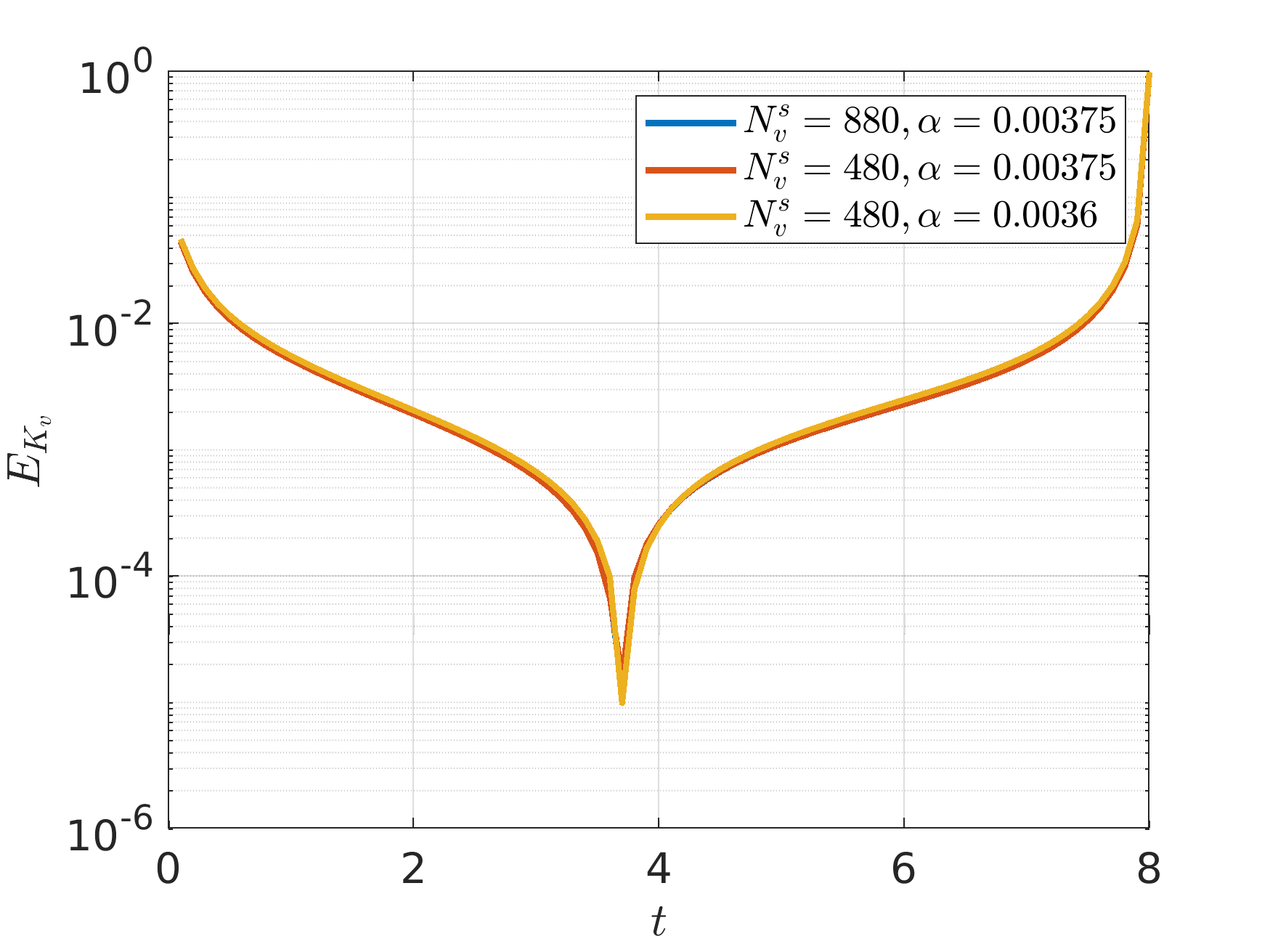}
      \end{overpic}
\begin{overpic}[width=0.45\textwidth]{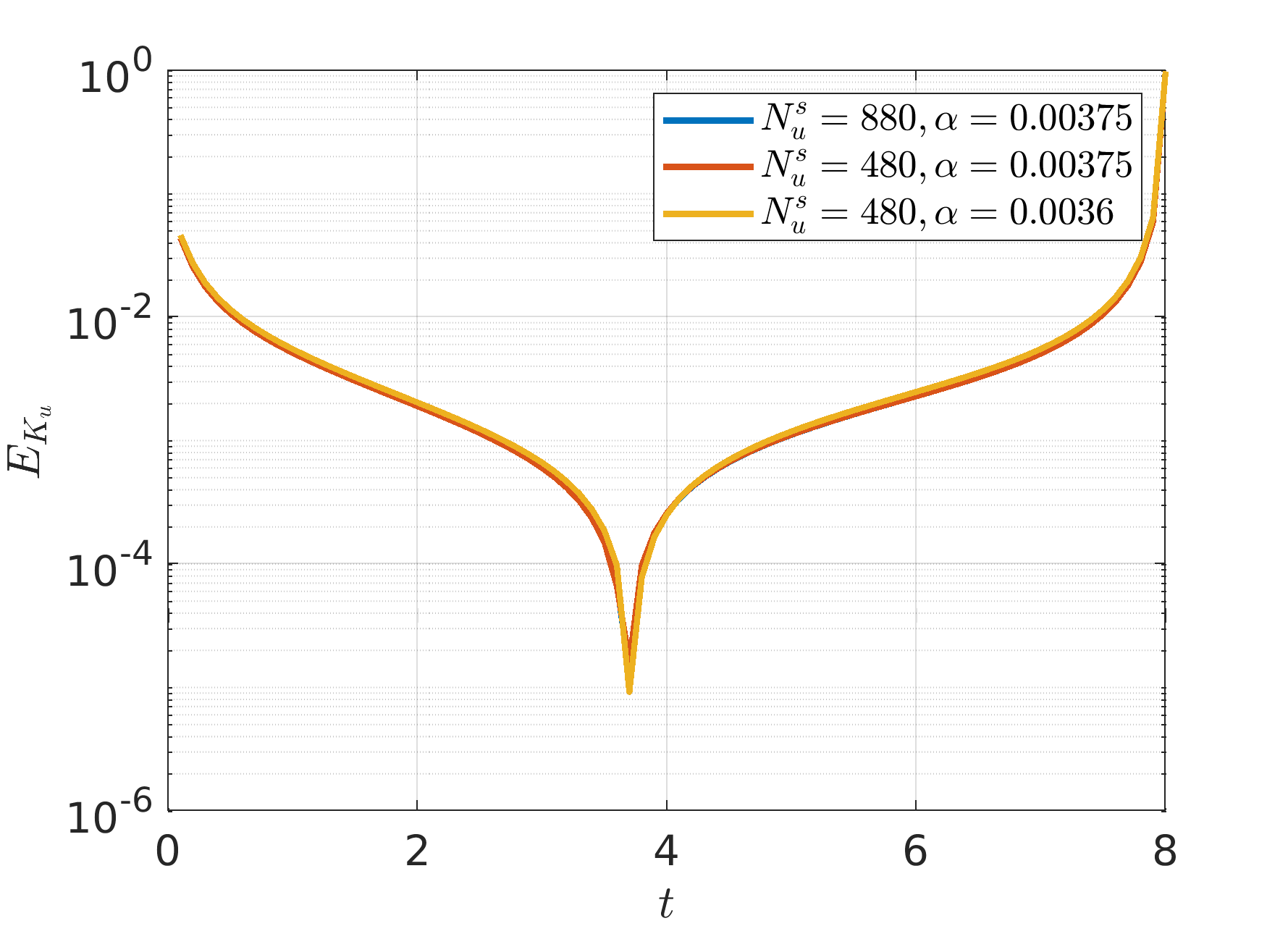}
      \end{overpic}\\
\caption{2D flow past a cylinder: time history of $L^2$ norm of the relative error for velocity fields (top),
pressure fields (center), and for kinetic energies of the system (bottom) for a different number of 
collected snapshots. We consider 2 modes for $\v$, $\u$ and $q$, and 1 mode for $\bar{q}$.}
\label{fig:err_delta}
\end{figure}

We take $\alpha = 0.00375$ (in the range under consideration but not in the training set)
to evaluate the performance of the parametrized ROM. 
A comparison between the two sampling spaces and the FOM is shown in Fig.~\ref{fig:err_delta}.
We observe that the differences are minimal. Thus, from now on we will use 6 sampling points
to reduce the computational time. Fig.~\ref{fig:err_delta} reports also the results for 
$\alpha = 0.0036$, which belongs to the training set. There is no significant change in accuracy
between $\alpha = 0.0036$ and $\alpha = 0.00375$.

Table \ref{tab:errors_delta} reports the maximum, minimum and average relative 
errors for the velocity and pressure fields, and for the kinetic energies of the system 
for filtering radius $\alpha = 0.00375$.
We see that Fig.~\ref{fig:err_delta} and \ref{fig:err_t} are vey similar. Moreover,  
the errors in Table \ref{tab:errors_delta} and \ref{tab:errors_t} are very close.
The POD modes are also very similar to those reported in Fig.~\ref{fig:modes2D}
and thus they are omitted. 

\begin{table}
\centering
\begin{tabular}{lcccccc}
\multicolumn{2}{c}{} \\
\cline{1-7}
 & $\u$ & $\v$ & $q$ & $\overline{q}$ & $K_u$ & $K_v$  \\
\hline
Maximum $E_\Phi$ & 9.1e-1  & 9.1e-1 & 3.7e-1 & 5.1e-1 & 9.8e-1 & 9.8e-1 \\
Miminum $E_\Phi$ & 1.1e-3 & 1.1e-3 & 2.7e-2 & 8.1e-2 & 1.6e-5 & 1.6e-5\\
Average $E_\Phi$ & 2.1e-2 & 2.3e-2 & 1.4e-1 & 9.6e-2 & 1.7e-2 & 1.7e-2\\
\hline
\end{tabular}
\caption{2D flow past a cylinder: maximum, minimum and average relative errors for the velocity and pressure fields, and for the kinetic energy of the system, for $\alpha = 0.00375$. The sampling frequency in time of the snapshots is 0.1 and we take 6 samples in the filtering radius interval.}
\label{tab:errors_delta}
\end{table}

\begin{figure}
\centering
       \begin{overpic}[width=0.45\textwidth]{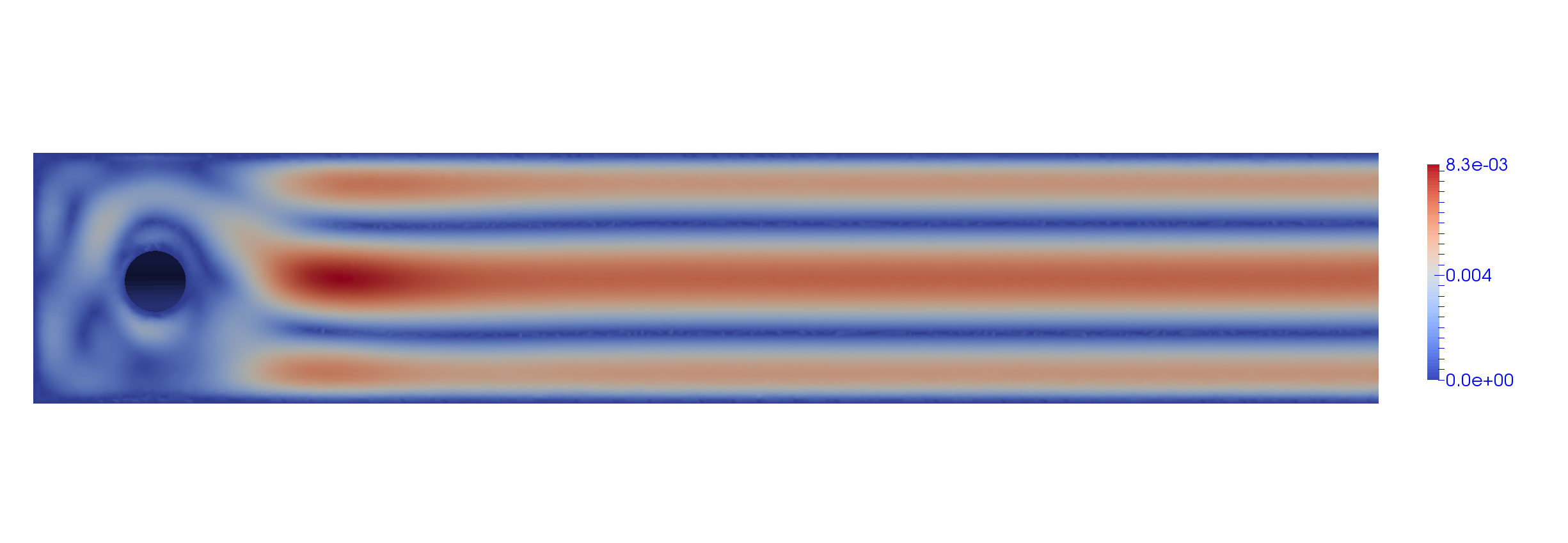}
        \put(20,30){$\v_{FOM} - \v_{ROM}$, $t = 1$}
      \end{overpic}
 \begin{overpic}[width=0.45\textwidth]{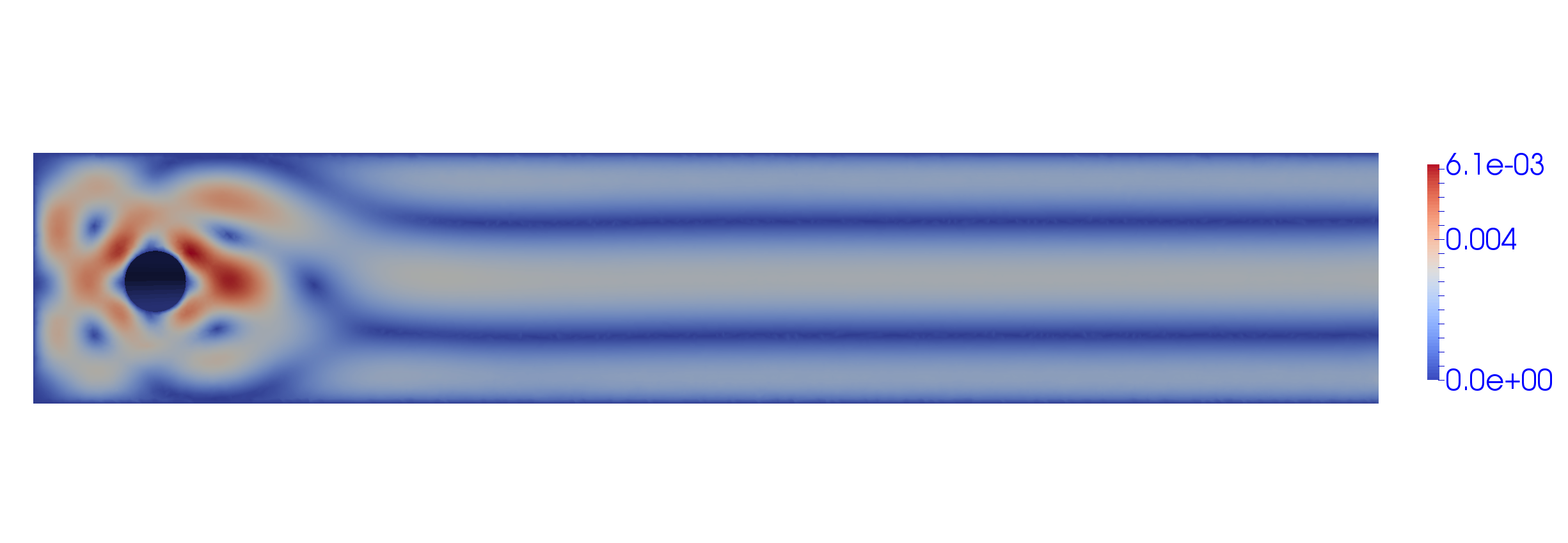}
        \put(20,30){$\v_{FOM} - \v_{ROM}$, $t = 5$}
      \end{overpic}\\
       \begin{overpic}[width=0.45\textwidth]{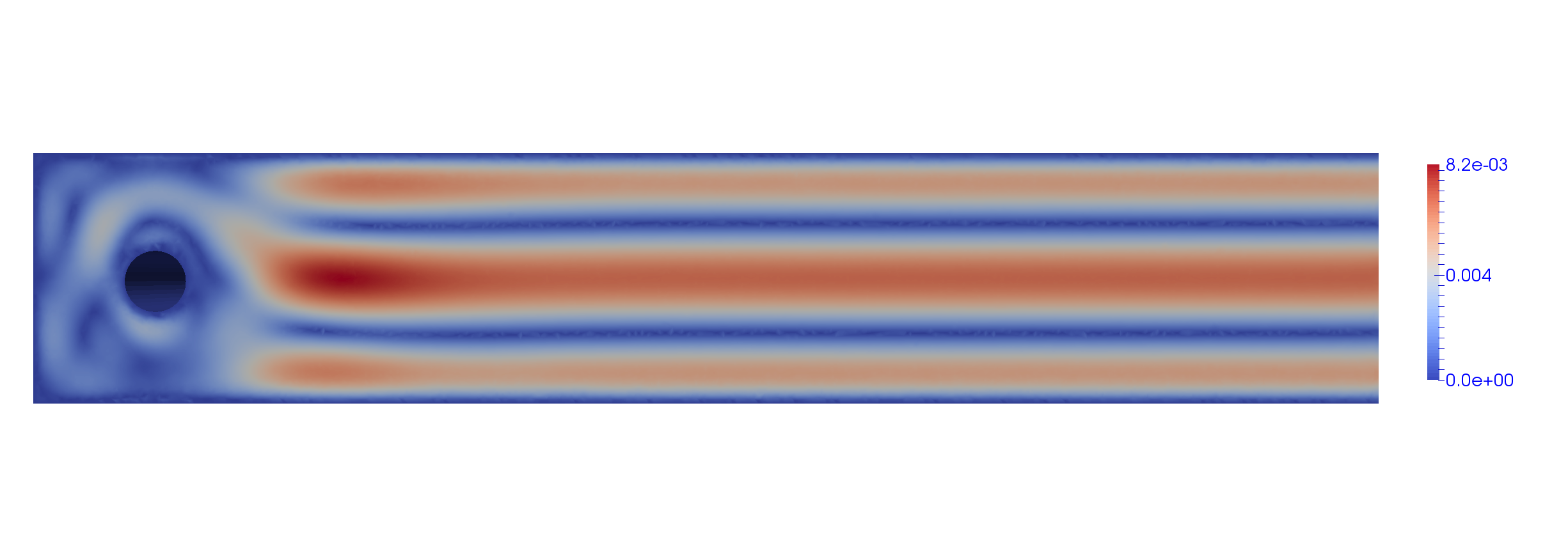}
        \put(20,30){$\u_{FOM} - \u_{ROM}$, $t = 1$}
      \end{overpic}
 \begin{overpic}[width=0.45\textwidth]{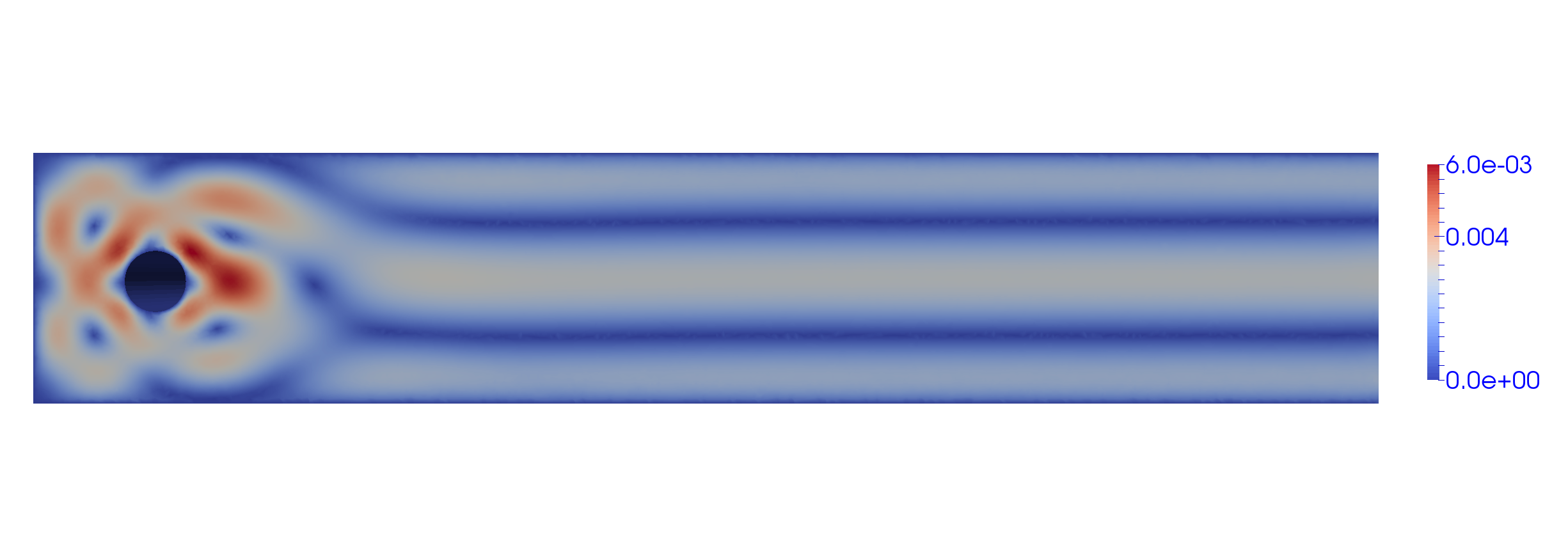}
        \put(20,30){$\u_{FOM} - \u_{ROM}$, $t = 5$}
      \end{overpic}\\
       \begin{overpic}[width=0.45\textwidth]{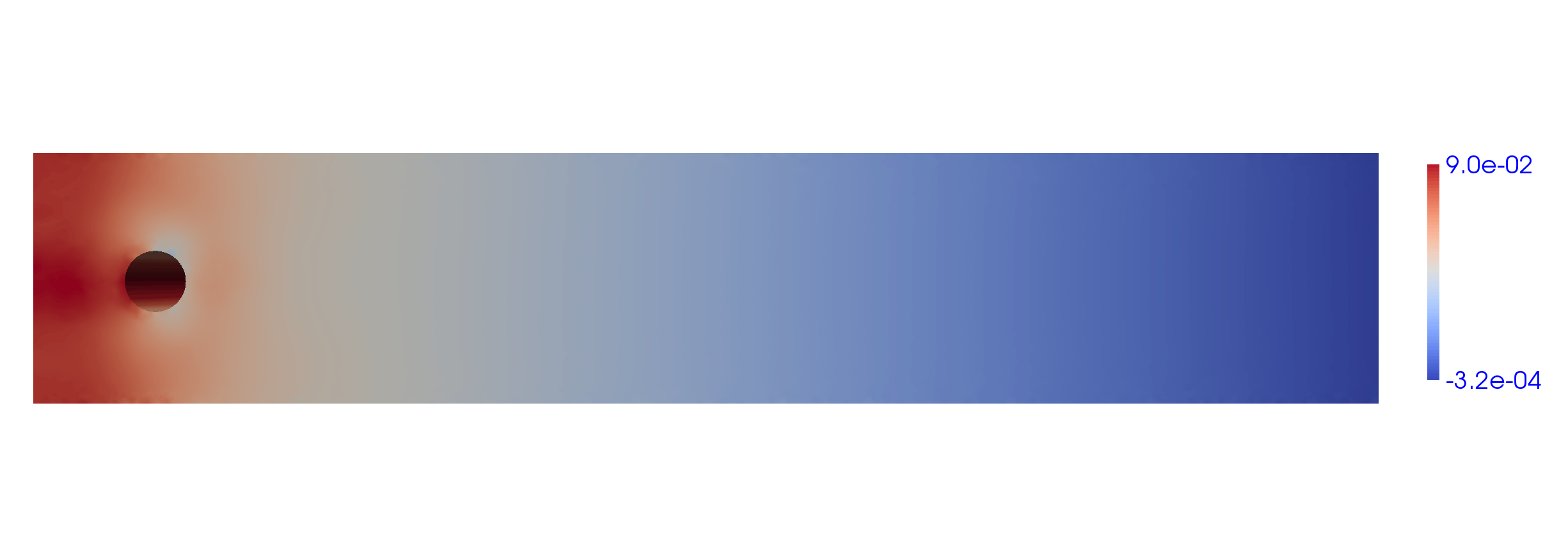}
        \put(20,30){$q_{FOM} - q_{ROM}$, $t = 1$}
      \end{overpic}
 \begin{overpic}[width=0.45\textwidth]{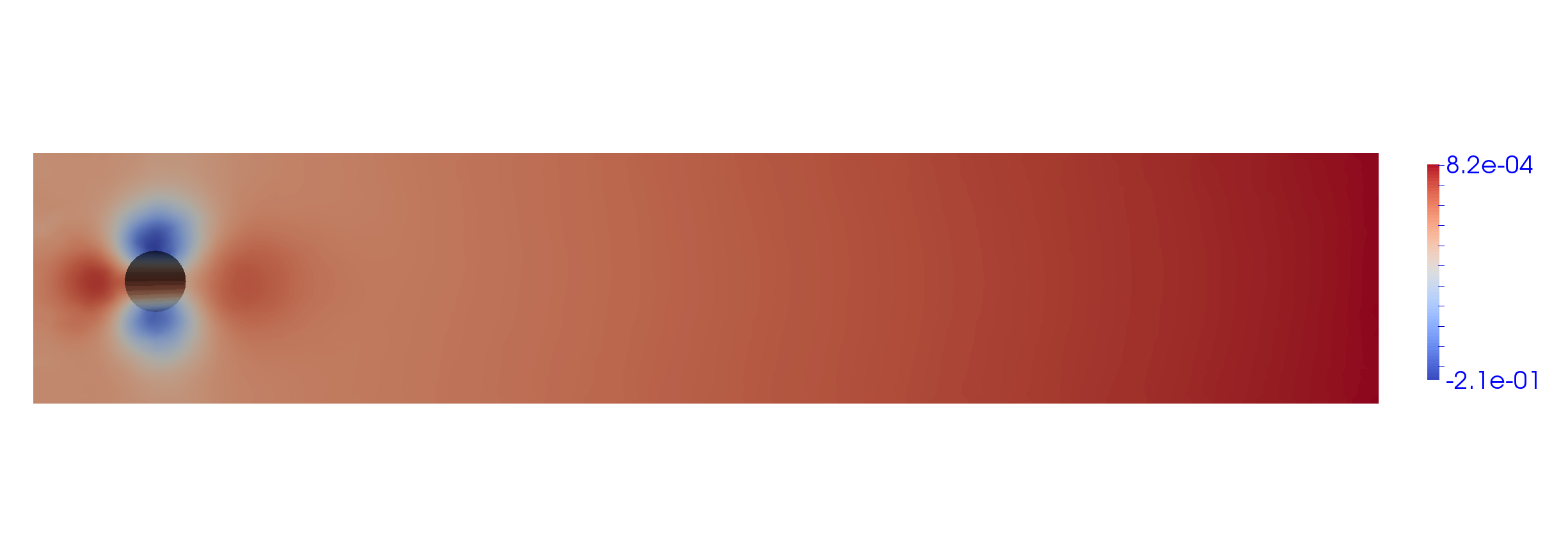}
        \put(20,30){$q_{FOM} - q_{ROM}$, $t = 5$}
      \end{overpic}\\
       \begin{overpic}[width=0.45\textwidth]{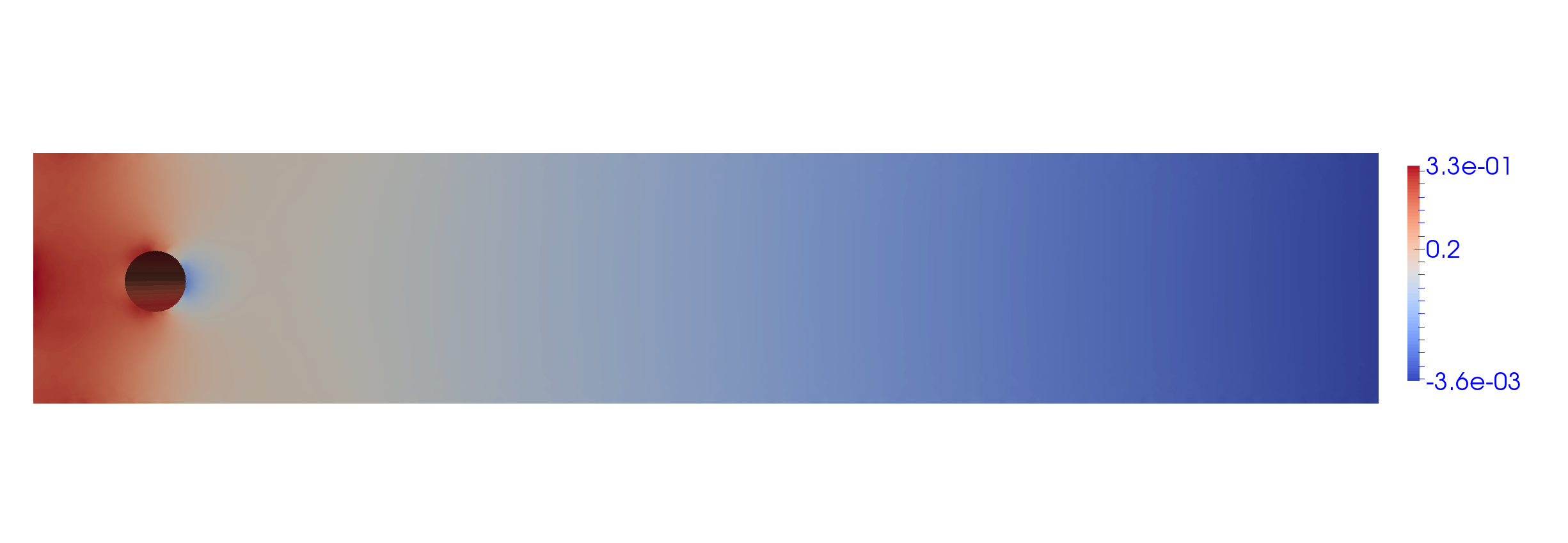}
        \put(20,30){$\overline{q}_{FOM} - \overline{q}_{ROM}$, $t = 1$}
      \end{overpic}
 \begin{overpic}[width=0.45\textwidth]{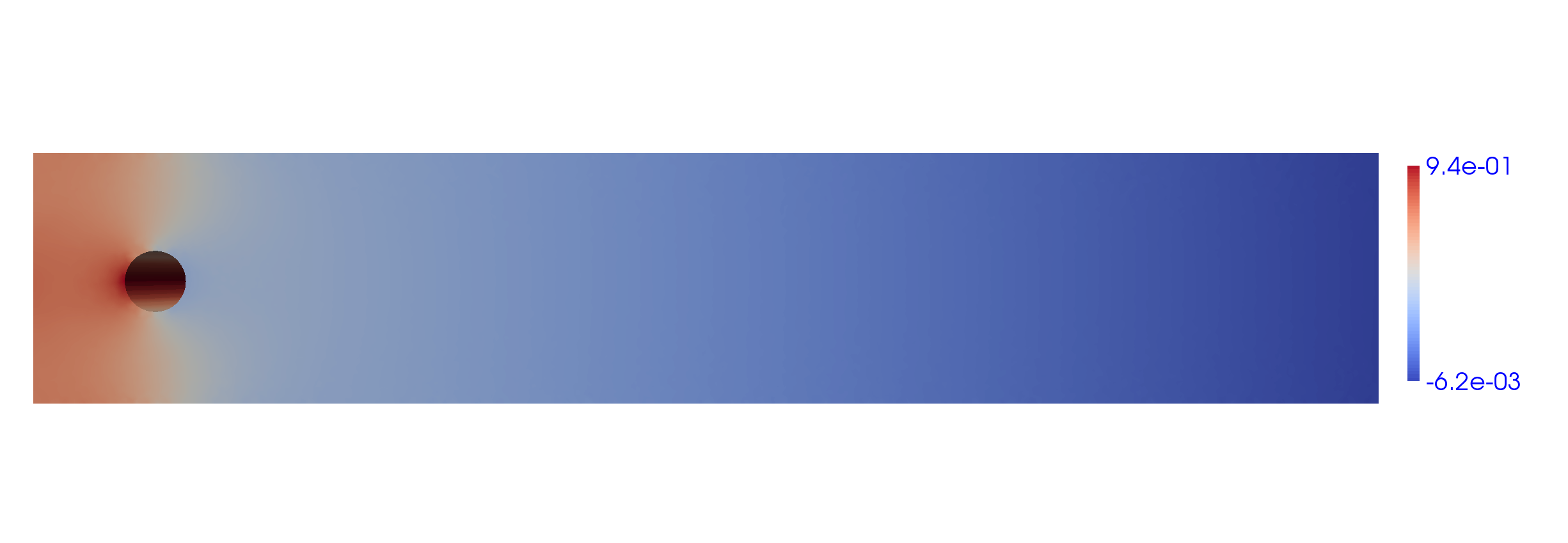}
        \put(20,30){$\overline{q}_{FOM} - \overline{q}_{ROM}$, $t = 5$}
      \end{overpic}
\caption{2D flow past a cylinder: Difference between FOM and ROM $\v$ (first row), $\u$ (second row), $q$ (third row),
and $\overline{q}$ (fourth row) at times $t = 1$ (left) and $t = 5$ (right) for $\alpha = 0.00375$. We consider 2 modes for $\v$, $\u$ and $q$, and 1 mode for $\bar{q}$.}
\label{fig:comp_delta}
\end{figure}

Fig.~\ref{fig:comp_delta} displays the difference between the computed FOM and ROM fields at two
different times: $t = 1$ and $t = 5$. Again, we see that our ROM approach is able to reproduce
the main flow features at different times, with a level of accuracy that is
comparable to the case  where time is the only parameter. 

\begin{figure}
\centering
 \begin{overpic}[width=0.45\textwidth]{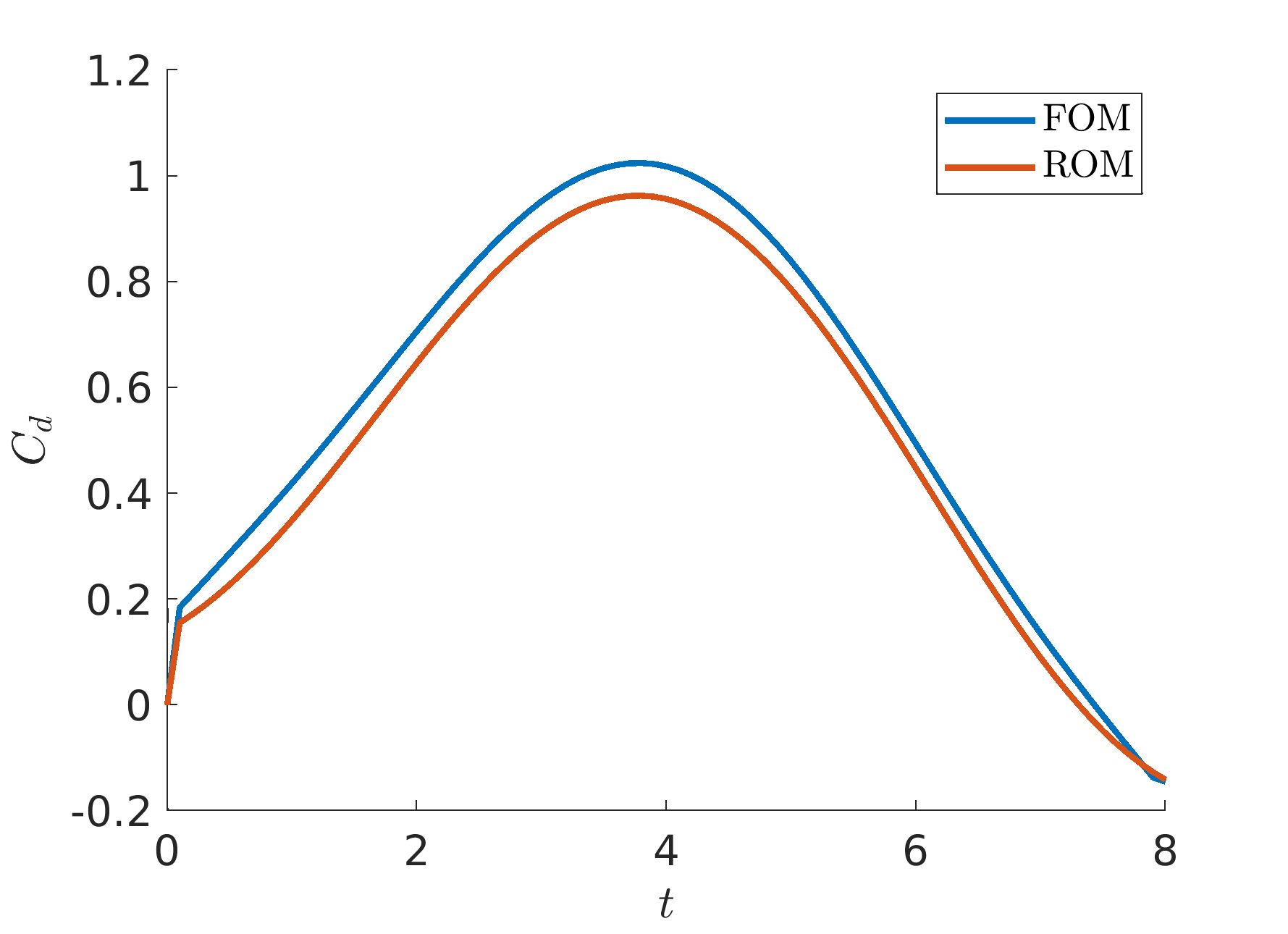}
      \end{overpic}
       \begin{overpic}[width=0.45\textwidth]{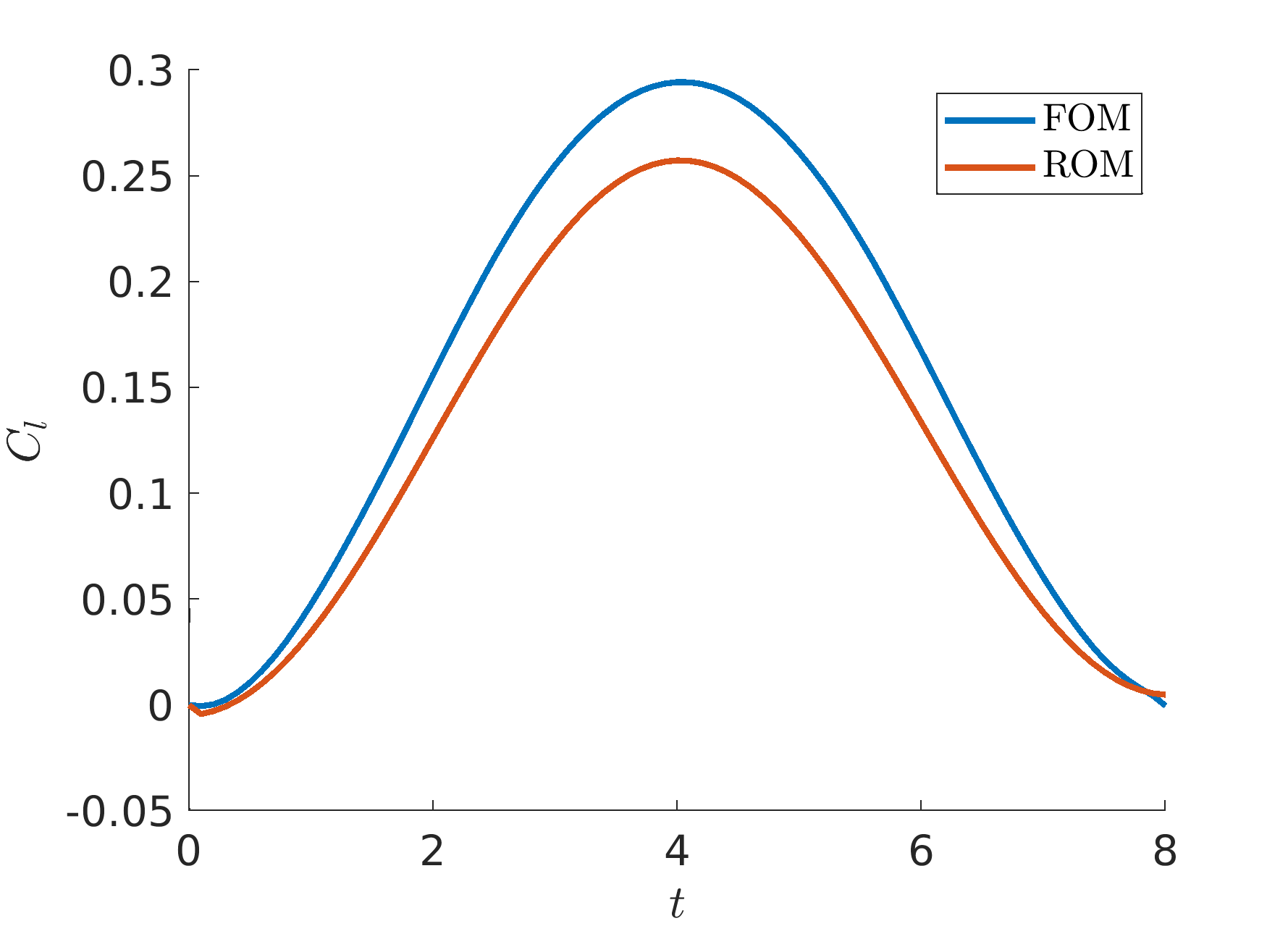}
      \end{overpic}\\
\caption{2D flow past a cylinder: comparison of the FOM aerodynamic coefficients $C_d$ (left) and $C_l$ (right), and the corresponding ROM reconstructions, for $\alpha = 0.00375$. We consider 2 modes for $\v$, $\u$ and $q$, and 1 mode for $\bar{q}$.}
\label{fig:coeff_delta}
\end{figure}

The FOM/ROM comparison for the drag and lift coefficients \eqref{eq:cd_cl} over time
for $\alpha = 0.00375$ is shown in Fig.~\ref{fig:coeff_delta}. We observe that also for this 
value of $\alpha$ the coefficients computed by the ROM underestimate the coefficients
computed by the FOM over almost the entire time interval under consideration, 
while the phase is perfectly reproduced. Table \ref{tab:coeffs_delta} reports
quantitative comparison between ROM and FOM in terms of maximum lift and 
drag coefficients and times at which the maxima occur. The corresponding relative error, as defined in 
\eqref{eq:E_c_l}-\eqref{eq:E_c_d}, are: $E_{c_l} = 0.126$, $E_{t_{c_l}} = 0$, $E_{c_d}= 0.06$, $E_{t_{c_d}} = 0$.
In switching from $\alpha = 0.0032$ to $\alpha = 0.00375$, 
the $C_d$ error decreases for (from 8.3\% to 6\%), while 
the $C_d$ error increases (from 8.9\% to 12.6\%).
The increase in the errore for $C_l$ could be due to the fact that towards
the center of the time interval the absolute error for pressure $q$ around
the cylinder is slightly larger for $\alpha = 0.00375$: compare Fig.~\ref{fig:comp_t} (third row, right) 
with Fig.~\ref{fig:comp_delta} (third row, right). 

\begin{table}[h]
\centering
\begin{tabular}{ccccc}
\multicolumn{5}{c}{} \\
\cline{1-5}
 & $t(c_{l,max})$ & $c_{l,max}$ & $t(c_{d,max})$& $c_{d,max}$   \\
\hline
FOM & 4.1  & 0.294 & 3.9 & 1.025  \\
ROM & 4.1  & 0.257 & 3.9 & 0.963  \\
\hline
\end{tabular}
\caption{2D flow past a cylinder: maximum lift and drag coefficients, and times at which the maxima occur for FOM and ROM, for $\alpha = 0.00375$.}
\label{tab:coeffs_delta}
\end{table}


\subsection{3D flow past a cylinder}

In this section, we aim at showing that our ROM approach can easily handle three-dimensional problems. 
The 3D benchmark we consider has been studied for the first time in \cite{turek1996} 
and further investigated in \cite{Bayraktar2012, John2006}.

We choose to adopt the classical Leray model, instead of the EF algorithm, which is know
to be over-diffusive especially in 3D. As a regularized ROM for convection-dominated problems,
the Leray model has been investigated in \cite{Sabetghadam2012, Iliescu2016, Gunzburger2019b, Xie2018, Wells2017}.
The discrete in time model reads as follows: Given velocities $\u^{n-1}$ and $\u^{n}$, at $t^{n+1}$ find $\u^{n+1}$, $q^{n+1}$,
$\overline{\u}^{n+1}$ and $\qbar^{n+1}$ such that:
\begin{align}
\rho\, \frac{3 \u^{n+1} - 4\u^{n} + \u^{n-1}}{2\Delta t}\,  + \rho\, \div \left(\u^* \otimes \u^{n+1}\right) - 2\mu\Delta\u^{n+1} +\nabla q^{n+1} & = 0,\label{eq:evolve-1.1_leray}\\
\div \u^{n+1} & = 0\label{eq:evolve-1.2_leray}, \\
\frac{\rho}{\Delta t} \overline{\u}^{n+1}  - \mubar \Delta\overline{\u}^{n+1} + \nabla \qbar^{n+1}  &= \frac{\rho}{\Delta t} \u^{n+1}, \label{eq:filter-1.1_leray} \\
\div \overline{\u}^{n+1} & = 0, \label{eq:filter-1.2_leray}
\end{align}
where $\u^* = 2 \overline{\u}^n-\overline{\u}^{n-1}$. Also with this model, 
we compute pressure fields and apply the filter for both FOM and ROM. However, 
for sake of brevity we will show results for $\u$ and $q$ only.

\begin{figure}
\centering
 \begin{overpic}[width=0.4\textwidth]{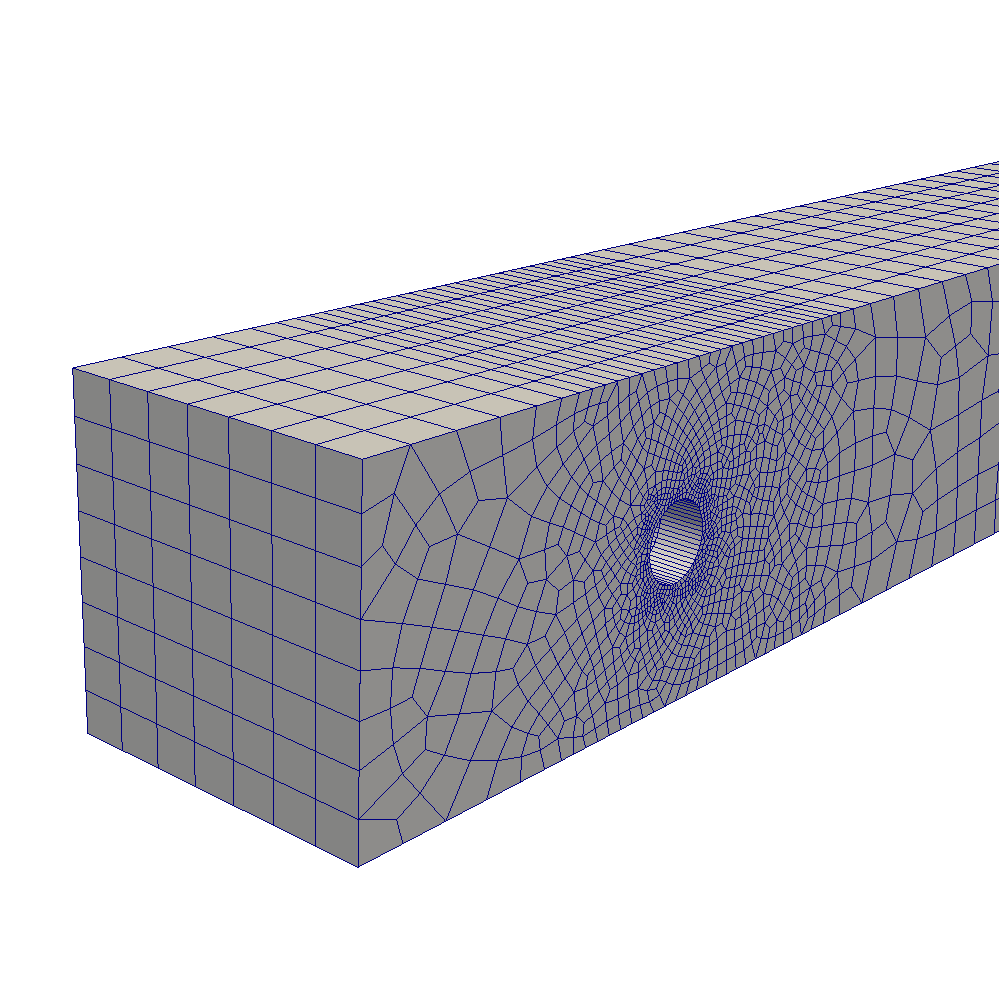}
      \end{overpic}
\begin{overpic}[width=0.5\textwidth]{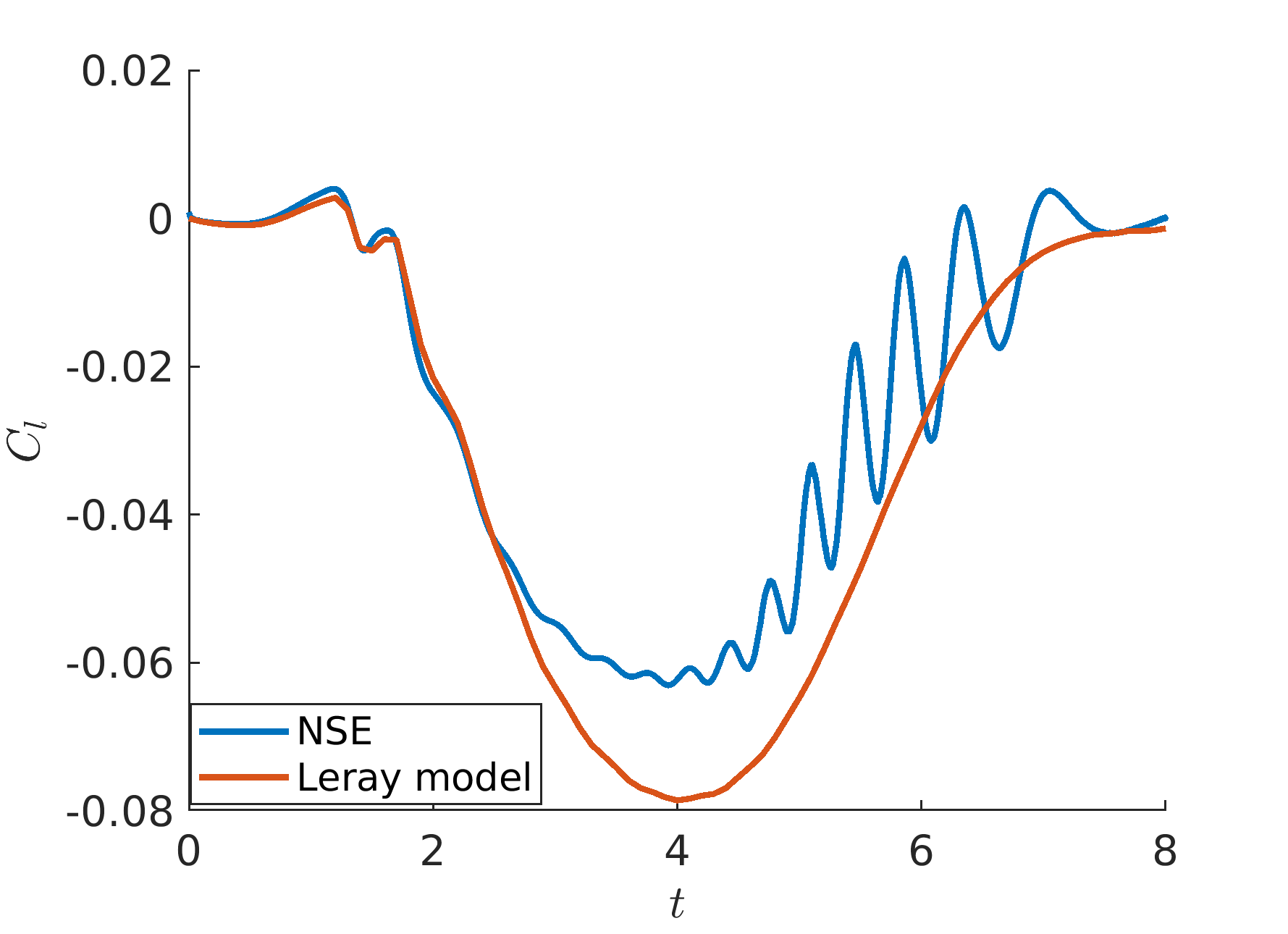}
      \end{overpic}\\
\caption{3D flow past a cylinder: (left) part of the mesh under consideration and (right) 
evolution of $c_l$ computed with the Leray model and NSE model.}
\label{fig:MESH3D}
\end{figure}

The computational domain is a 2.5 $\times$ 0.41 $\times$ 0.41 parallelepiped with a 
cylinder whose axis is parallel to the $z$-axis and center is located at (0.5, 0.2)
when taking the bottom left corner of the channel as the origin of the axes. Fig.~\ref{fig:MESH3D} (left) 
shows part of the computational domain. The channel is filled with fluid with density $\rho = 1$ 
and viscosity $\mu = 0.001$. We impose a no slip boundary condition on the channel
walls and on the cylinder. At the inflow, we prescribe the following velocity profile:
\begin{align}\label{eq:cyl_bc3d}
\u(0,y,z,t) = \left(\dfrac{36}{0.41^4} \sin\left(\pi t/8 \right) y z \left(0.41 - y \right) \left(0.41 - z \right), 0, 0\right), \quad y,z \in [0, 0.41], \quad t \in (0, 8].
\end{align}
In addition, on the channel walls, cylinder, and at the inlet we impose ${\partial q}/{\partial \n} = {\partial \overline{q}}/{\partial \n} = 0$ 
where $\n$ is the outward normal. At the outflow, we prescribe $\nabla \u \cdot  \n = 0$ and $q = \qbar = 0$. 
Note that the Reynolds number is time dependent, with $0 \leq Re \leq 100$ \cite{turek1996, Bayraktar2012, John2006}. 
Like for the 2D benchmark, we start the simulations from fluid at rest.

We consider a hexahedral grid with $h_{min} = 9e-3$ and $h_{max} = 6.6e-2$ and a total of 1.07e4 cells. 
This level of refinement is far from the one required by DNS.
The mesh features very low values of maximum non-orthogonality (34$^\circ$), average non-orthogonality (7$^\circ$), skewnwss (0.6), 
and maximum aspect ratio (25). In addition, the mesh is refined next to the cylinder, like the meshes 
used in \cite{Bayraktar2012, John2006}.

Before applying the ROM, we test the Leray model and compare its results with the results produced by a
NSE solver. Like in the 2D case, we use a second-order accurate Central Differencing (CD) 
scheme \cite{Lax1960} for the discretization of the convective term. For the NSE solver, 
we set $\Delta t = 1e-4$, while for the Leray model we choose $\Delta t = 5e-3$. 
The larger time step for the Leray model is motivated by an attempt to 
contain the artificial viscosity $\bar{\eta}$ defined in eq.~(\ref{eq:filter-1.1}). 
We set $\alpha = 0.0032$ and refer to \cite{Girfoglio2019} for details on this choice. 
Fig.~\ref{fig:MESH3D} (right) shows the evolution of $c_l$ over time computed by the Leray and NSE model.
We observe that the Leray model dampens the unphysical oscillations in the NSE solution 
and reduces the maximum lift coefficient.
We report in Table \ref{tab:coeffs_comp_3D} the computed values of maximum drag and lift coefficients, together with
the reference values fom \cite{turek1996, Bayraktar2012, John2006}. 
From Fig.~\ref{fig:MESH3D} (right), we can conclude that the Leray model outperforms the NSE model on a coarse mesh.  

\begin{table}
\centering
\begin{tabular}{lcc}
\multicolumn{2}{c}{} \\
\cline{1-3}
 & $C_d$ & $C_l$  \\
\hline
 NSE algorithm & 3.21  &  0.004 \\
 Leray model & 3.12  &  0.0028  \\
 Ref. values \cite{turek1996, Bayraktar2012, John2006} & [3.2, 3.3]  &  [0.002, 0.004] \\
\hline
\end{tabular}
\caption{3D flow past a cylinder: Maximum drag and lift coefficients given by NSE and Leray model. The bottom row reports the reference values from \cite{turek1996, Bayraktar2012, John2006}.}
\label{tab:coeffs_comp_3D}
\end{table}

\begin{figure}
\centering
 \begin{overpic}[width=0.45\textwidth]{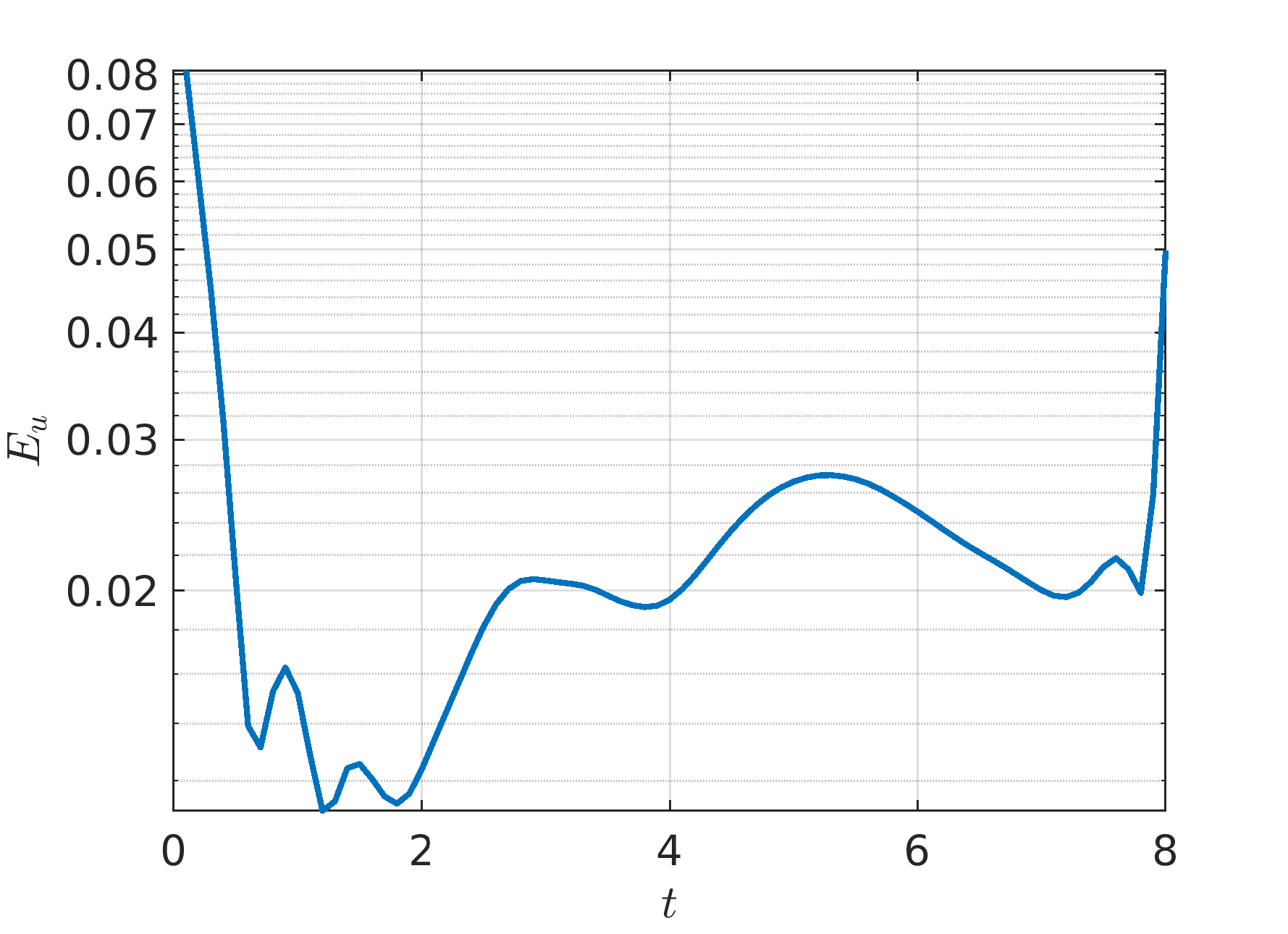}
      \end{overpic}
\begin{overpic}[width=0.45\textwidth]{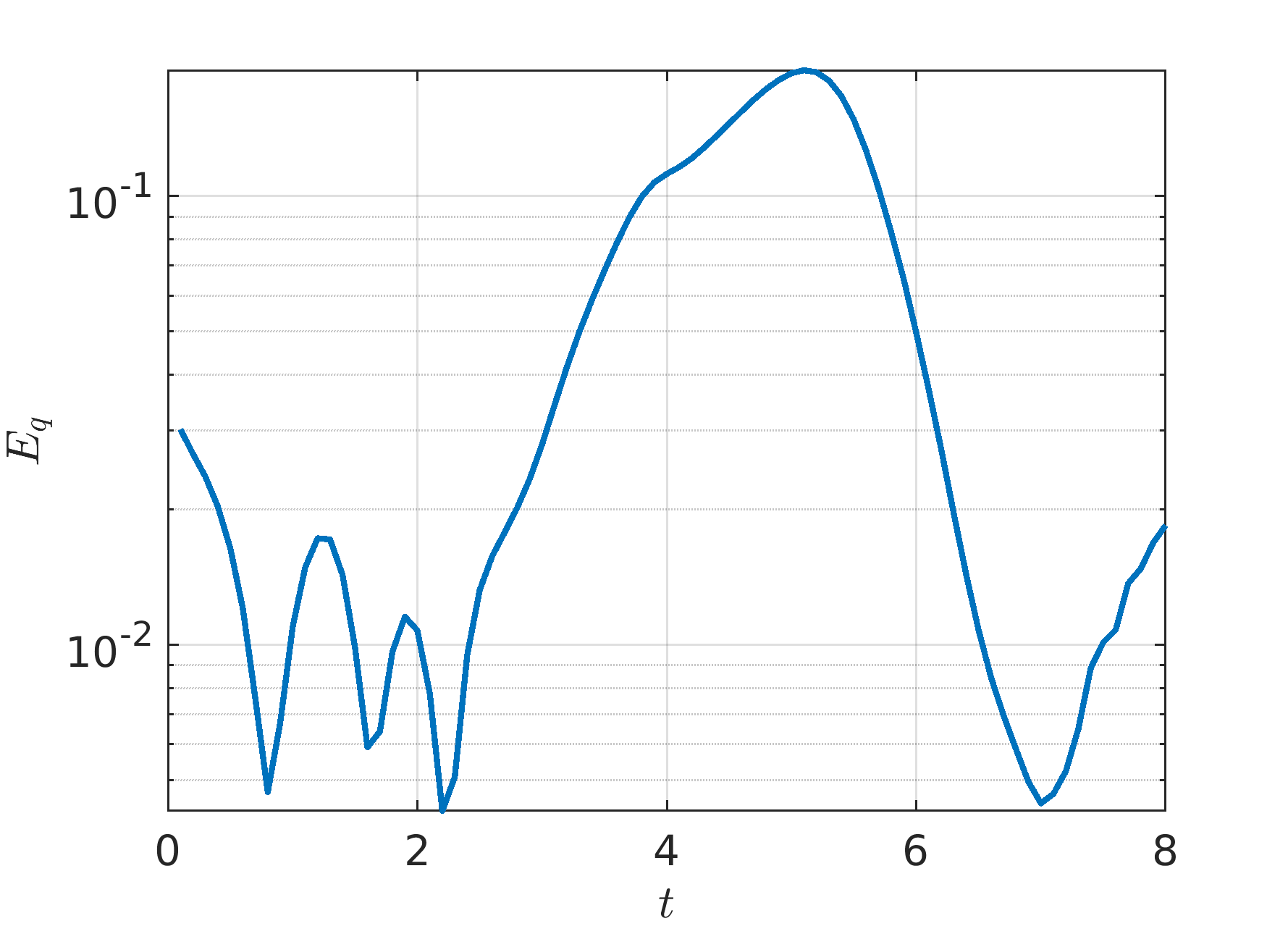}
      \end{overpic}\\
\caption{3D flow past a cylinder: time history of $L^2$-norm of the relative error for velocity $\u$ (left) 
and pressure $q$ (right). We consider 9 modes for $\u$ and 4 modes for $q$.}
\label{fig:err3D}
\end{figure}

The snapshots are collected every 0.1 s using an equispaced grid method in time. Fig.~\ref{fig:err3D} shows error \eqref{eq:error1} 
for the velocity $\u$ and pressure $q$ field over time. 
As for the 2D case, the largest relative errors for $\u$ occur around the beginning 
and end of the simulation, while the relative error for pressure $q$ reaches its maximum value
around $t = 5$. 
The minimum, maximum and average (over time) relative errors are reported in Table \ref{tab:errors3D}. 
Average errors for both the velocity and pressure field are comparable 
with the ones found for the 2D benchmark.


\begin{table}
\centering
\begin{tabular}{lcc}
\multicolumn{2}{c}{} \\
\cline{1-3}
 & $\u$ &  $q$   \\
\hline
Maximum $E_\Phi$ & 8.1e-2  & 1.9e-1  \\
Miminum $E_\Phi$ & 1.1e-2 & 4.3e-3 \\
Average $E_\Phi$ & 2.2e-2 & 5.3e-2 \\
\hline
\end{tabular}
\caption{3D flow past a cylinder: maximum, minimum and average (over time) relative errors for velocity $\u$
and pressure $q$.}
\label{tab:errors3D}
\end{table}

Table \ref{tab:cum3D} lists the first 10 cumulative eigenvalues, based on the first 15 most energetic POD modes. 
We see that 9 modes for $\u$ and 4 modes for $q$ are sufficient to retain 99.99\% of the energy contained in the snapshots. 
With respect to the 2D case, a larger number of modes is necessary. 
This could be due to the fact that the EF algorithm used for the 2D test 
introduces more numerical dissipation and dampens high frequency modes.

\begin{table}
\centering
\begin{tabular}{ccc}
\multicolumn{3}{c}{} \\
\cline{1-3}
Mode number & $\u$ & $q$   \\
\hline
 1 & 0.972257  & 0.957533 \\
 2 & 0.989825  & 0.999779 \\
 3 & 0.996379  & 0.999922 \\
 4 & 0.998582  & 0.999982 \\
 5 & 0.999232  & 0.999989 \\
 6 & 0.999609  & 0.999994 \\
 7 & 0.999760 & 0.999996 \\
 8 & 0.999862  & 0.999997 \\
 9 & 0.999909  & 0.999998 \\
 10 & 0.999944 &  0.999999 \\
\hline
\end{tabular}
\caption{3D flow past a cylinder: first 10 cumulative eigenvalues for $\u$ and $q$.}
\label{tab:cum3D}
\end{table}

\begin{figure}
\centering
 \begin{overpic}[width=0.4\textwidth]{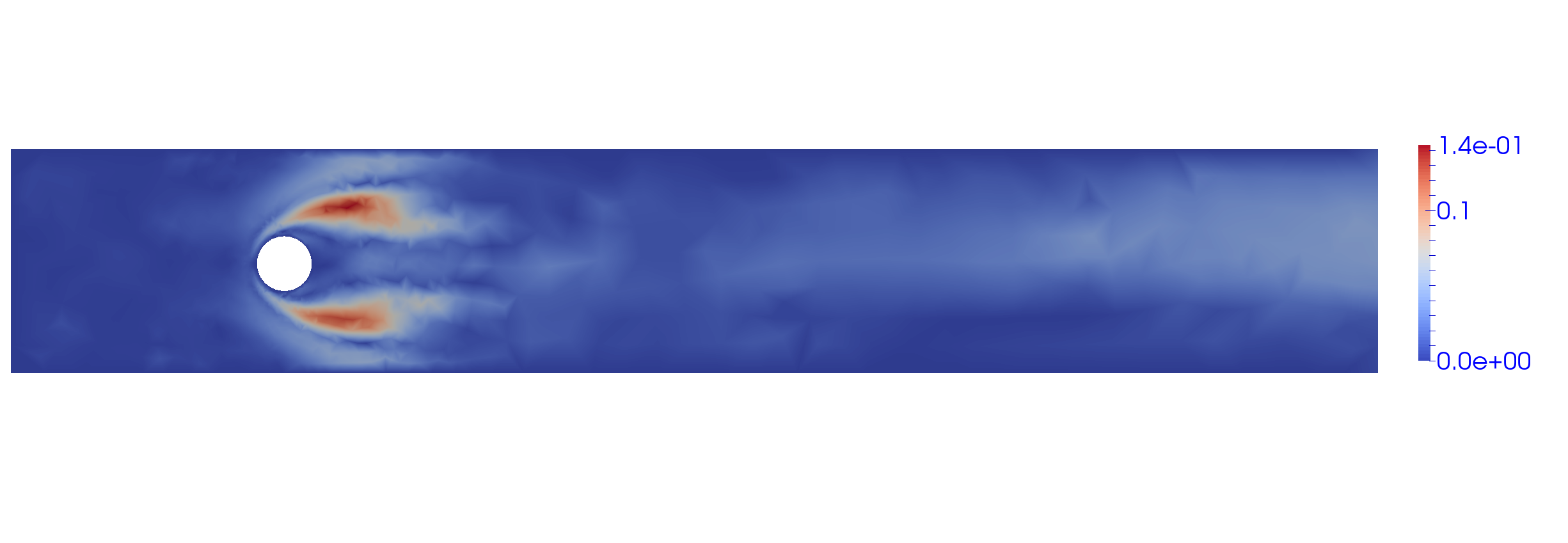}
      \end{overpic}
 \begin{overpic}[width=0.4\textwidth]{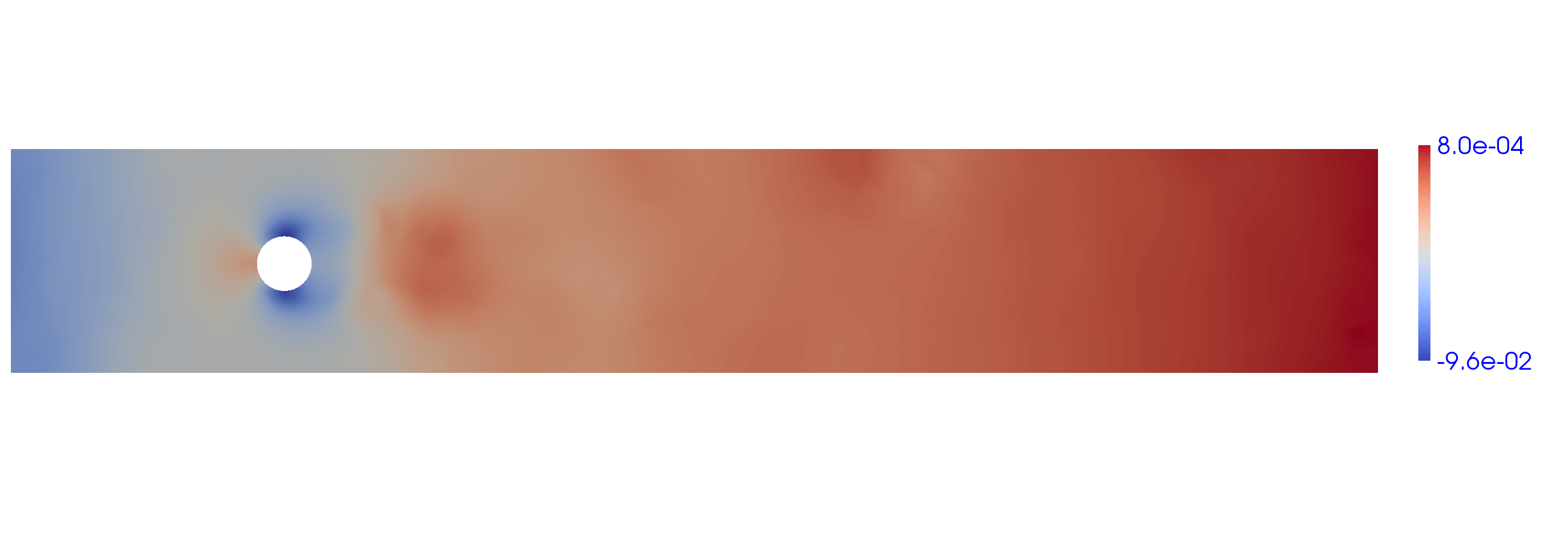}
      \end{overpic}
\caption{3D flow past a cylinder: difference between FOM and ROM $\u$ (left) and $q$ (right) on the midsection at time $t = 5$. We consider 9 modes for $\u$ and 4 modes for $q$.}
\label{fig:comp3D}
\end{figure}

Fig.~\ref{fig:comp3D} displays the difference between the solutions on the midsection 
($z = 0.205$) computed by FOM and ROM
at time $t = 5$. The maximum absolute errors between the FOM and the ROM for both $\u$ and $q$ are of the order $0.1$. In general, we can still conclude that 
our ROM approach is able to reproduce the main flow features. 

\begin{figure}
\centering
 \begin{overpic}[width=0.45\textwidth]{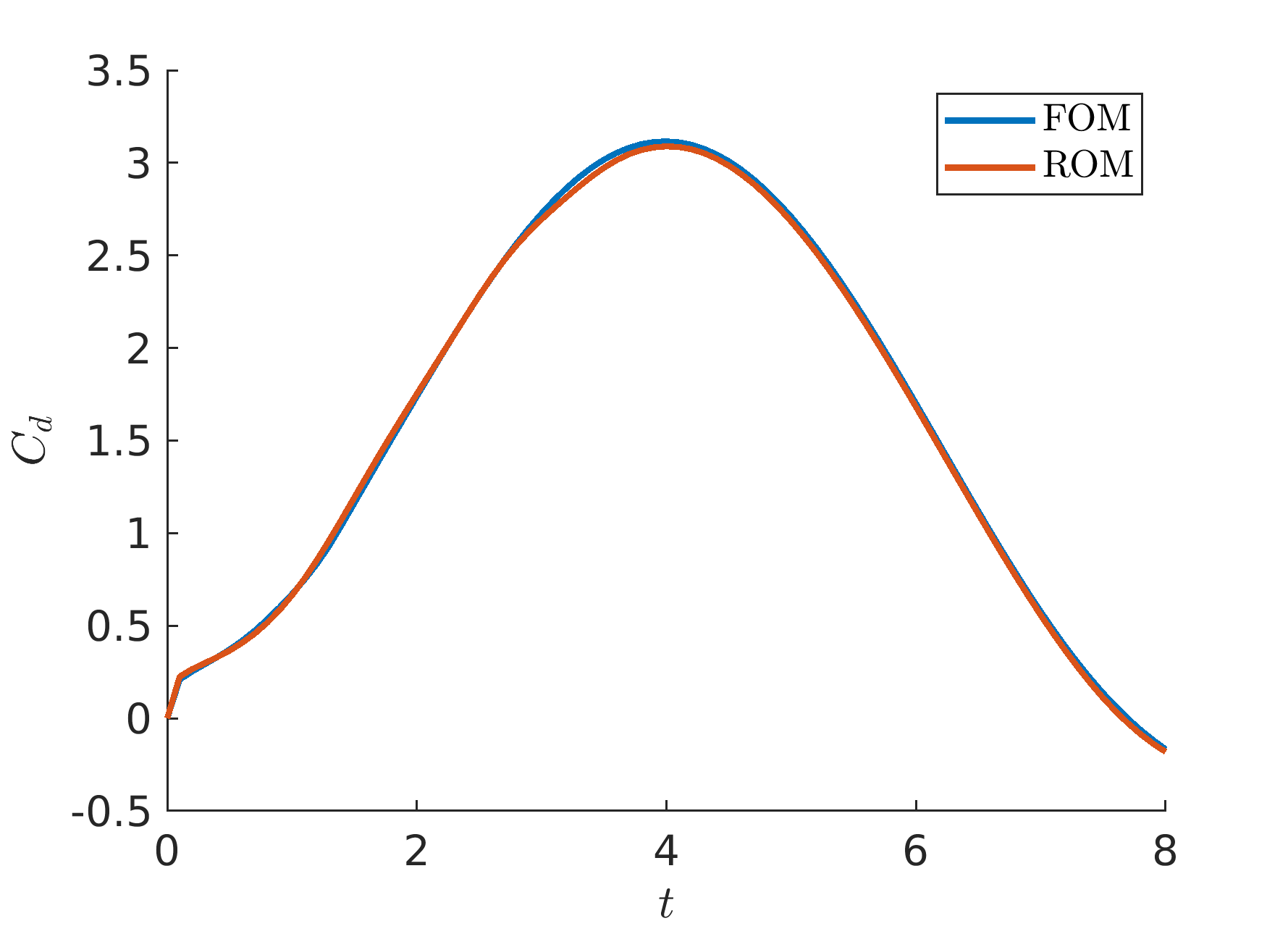}
      \end{overpic}
       \begin{overpic}[width=0.45\textwidth]{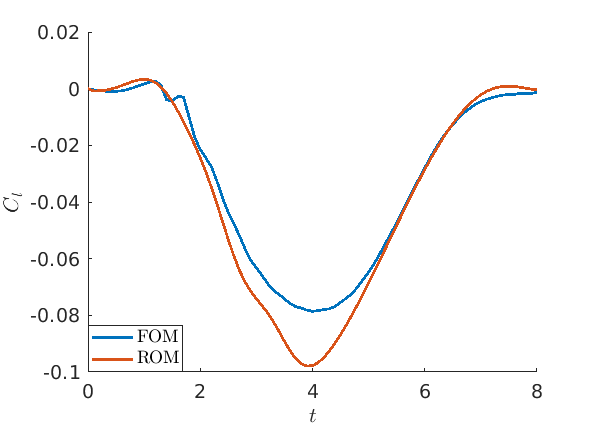}
      \end{overpic}
       \begin{overpic}[width=0.45\textwidth]{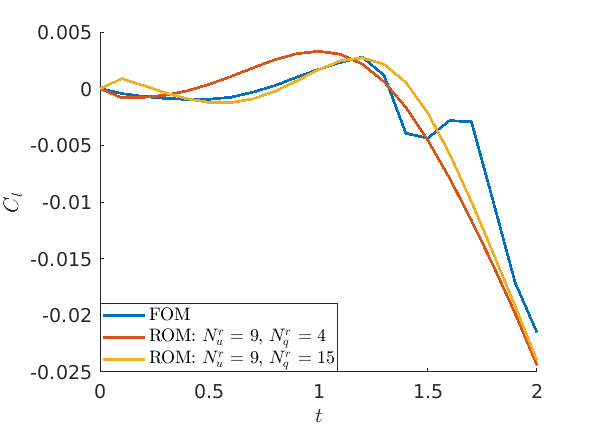}
      \end{overpic}\\
\caption{3D flow past a cylinder: drag coefficient $C_d$ (top left) and lift coefficient $C_l$ (top right)
computed by FOM and ROM for $N_u^r = 9$ and $N_q^r = 4$. A zoomed-in view of the lift coefficient around
the time the maximum is reached is shown in the bottom row for two different values of $N_q^r$.}
\label{fig:coeff3D}
\end{figure}

We report in Fig.~\ref{fig:coeff3D} the evolution of the lift and drag coefficients computed by 
FOM and ROM with two different values of $N_q^r$. 
We see that the amplitude of the drag coefficient is reproduced with excellent accuracy 
over the entire time interval. However, the ROM reconstruction of the lift coefficient
appears to be more critical, as already observed in the 2D test. 
Again, this could be due to the fact that larger errors for pressure $q$ are localized 
close to the cylinder, as one can see in Fig.~\ref{fig:comp3D} (right). 
Table \ref{tab:coeffs3D_1} compares the maximum lift and drag coefficients given by ROM and FOM. 
The relative errors related for the maximum values are $E_{c_d} = 0.3\%$ and $E_{c_l} = 17.9\%$. 
A better reconstruction of the maximum $C_l$ is obtained by using 
a higher number of modes for the pressure (i.e., larger $N_q^r$) 
in order to reconstruct more accurately the low amplitude oscillatory pattern around the time
the peak is reached. See Fig.~\ref{fig:coeff3D} (bottom). 
When switching from $N_q^r = 4$ to $N_q^r = 15$, the $C_l$ error decreases from $17.9\%$ to $3.6\%$.

\begin{table}[h]
\centering
\begin{tabular}{lcc}
\multicolumn{3}{c}{} \\
\cline{1-3}
 & $c_{l,max}$ & $c_{d,max}$   \\
\hline
FOM & 0.0028 & 3.12   \\
ROM & 0.0033 & 3.11  \\
\hline
\end{tabular}
\caption{3D flow past a cylinder: maximum lift and drag coefficients for FOM and ROM.}
\label{tab:coeffs3D_1}
\end{table}

Finally, we comment on the computational costs. 
The total CPU time required by a FOM simulation is 540 s. 
Our ROM approach takes 28 s. The speed-up is about 19, which is much smaller than the speed-up
found for the 2D test. We identified two possible reasons: (i) the value of $\bar{\mu}$ is smaller that one used for the 2D test and the high-fidelity solver converges faster (ii) the mesh used for the 3D test is coarser than one used for the 2D test.

\section{Conclusions and perspectives}\label{sec:conclusions}

This work presents a POD-Galerkin based reduced order method for a Leray model implemented through the
Evolve-Filter algorithm. Unlike the large majority of the works on Leray models, we choose a 
Finite Volume method because of its computational efficiency.  The novelties of the proposed ROM
approach are: (i) spatial filtering applied both for the collection of the snapshots and in the reduced order model,
(ii) the reconstruction of the pressure fields, and
(iii) the use of different POD basis and coefficients to approximate the velocity and pressure fields
in the two steps of the Evolve-Filter algorithm. 
We assessed our ROM approach through two classical benchmarks: 2D and 3D flow past a cylinder.
We found that our ROM can capture the flow features with an accuracy comparable to other ROMs
applied to similar benchmarks in \cite{Stabile2017, Stabile2018, Star2019}. In addition, we quantified
the relative error in the amplitude and phase of drag and lift coefficients computed by ROM and FOM. 
For the 2D test case, we also performed a parametric study with respect to the filtering radius.

A natural extension to the work presented in this manuscript is the development of 
a ROM for a Leray model with a nonlinear differential filter. To this aim, we are working in order to extend the approach used in \cite{Hijazi2020, GeorgakaStabileStarRozzaBluck2020}, based on the idea of merging/combining projection-based techniques with data-driven reduction strategies. In particular, the strategy in \cite{Hijazi2020, GeorgakaStabileStarRozzaBluck2020} exploits a data-driven reduction method to approximate the eddy viscosity solution manifold and a classical POD-Galerkin projection approach for the velocity and the pressure fields, respectively. We are going to use a data-driven reduction method to approximate the indicator function governing the amount of regularization introduced in the model in the nonlinear framework.  

\section{Acknowledgements}\label{sec:acknowledgements}
We acknowledge the support provided by the European Research Council Executive Agency by the Consolidator Grant project AROMA-CFD ``Advanced Reduced Order Methods with Applications in Computational Fluid Dynamics" - GA 681447, H2020-ERC CoG 2015 AROMA-CFD, PI G. Rozza, and INdAM-GNCS 2019-2020 projects.
This work was also partially supported by US National Science Foundation through grant DMS-1620384
and DMS-195353.

\bibliographystyle{plain}
\bibliography{ROM_GQR.bib}

\end{document}

%% file: packages.tex
\usepackage{amsmath}
\usepackage{amssymb}
\usepackage{mathrsfs}
\usepackage{color}
\usepackage{graphicx}
\usepackage[percent]{overpic}
\usepackage{url}
\usepackage{subfig}
\usepackage{enumerate}
\usepackage{overpic}
\usepackage{amsthm}
\usepackage{moreverb}
\usepackage[pdftex,colorlinks,bookmarksopen,bookmarksnumbered,citecolor=red,urlcolor=red]{hyperref}
\usepackage{authblk}

%% file: commands.tex
\def\b{{\bm b}}

\def\n{{\bm n}}
\def\u{{\bm u}}
\def\v{{\bm v}}

\def\0{\boldsymbol{0}}

\def\mubar{\overline{\mu}}
\def\ubar{\overline{\u}}
\def\vbar{\overline{\v}}

\def\qbar{\overline{q}}
\def\dt{\partial_t}

\def\cl {\nonumber \\}
\def\el {\nonumber }

\newtheorem{rem}{Remark}[section]

\newcommand{\bm}[1]{\mbox{\boldmath{$#1$}}}
\def\div{\nabla\cdot}

\newcommand{\anna}[2][cyan]{\emph{\textcolor{#1}{#2}}}